\documentclass[5p]{elsarticle}

\usepackage{setspace}
\usepackage{amsmath,amssymb,amsthm,graphicx,mathtools}
\usepackage[export]{adjustbox}
\graphicspath{{pics/}} % Specifies the directory where pictures are stored
\usepackage[caption=false]{subfig}
\usepackage{algorithm,algpseudocode,color}

\journal{Journal of \LaTeX\ Templates}

\newcommand{\cc}{\mathbf{c}}
\newcommand{\dd}{\mathbf{d}}
\newcommand{\vv}{\mathbf{v}}
\newcommand{\nn}{\mathbf{n}}

\newcommand{\tl}{\mathbf{t}}
\newcommand{\bbs}{\mathbb{S}}
\newcommand{\pp}{\mathbf{p}}
\newcommand{\RR}{\mathbb{R}}
\newcommand{\vcal}{\mathcal{V}}

\def\eem{\text{e$-$}}
\newcommand{\ee}{\mathbf{e}}
\newtheorem{definition}{Definition}
\providecommand{\restr}[2]{
  \left.\kern-\nulldelimiterspace
  #1 % the function
  \right|_{#2} % this is the delimiter
}
\providecommand{\abs}[1]{{\lvert#1\rvert}}
\providecommand{\norm}[1]{{\lVert#1\rVert}}

\begin{document}
\begin{frontmatter}

\title{Region extraction in mesh intersection}

\author[mymainaddress]{Pablo Antolin}
\ead{pablo.antolin@epfl.ch}

\author[mymainaddress,mysecondaryaddress]{Annalisa Buffa}
\ead{annalisa.buffa@epfl.ch}

\author[mymainaddress]{Emiliano Cirillo\corref{mycorrespondingauthor}}
\cortext[mycorrespondingauthor]{Corresponding author}
\ead{emiliano.cirillo@epfl.ch}

\address[mymainaddress]{Institute of Mathematics, \'Ecole Polytechnique F\'ed\'erale de Lausanne (EPFL), Lausanne, 1015,
Switzerland}
\address[mysecondaryaddress]{Istituto di Matematica Applicata e Tecnologie Informatiche
``Enrico Magenes'' del CNR, 27100 Pavia, Italy}

\begin{abstract}
	Region extraction is a very common task in both Computer Science and Engineering with several applications in object recognition and
  motion analysis, among others. Most of the literature focuses on regions delimited by straight lines, often in the special case of
  intersection detection among two unstructured meshes. While classical region extraction algorithms for line drawings and mesh
  intersection algorithms have proved to be able to deal with many applications, the advances in Isogeometric Analysis require a
  generalization of such problem to the case in which the regions to be extracted are bounded by an arbitrary number of curved segments.
  In this work we present a novel region extraction algorithm that allows a precise numerical integration of functions defined in different
  spline spaces. The presented algorithm has several interesting applications in contact problems, mortar methods, and quasi-interpolation
  problems.
\end{abstract}


\end{frontmatter}

% \twocolumn

\section{Introduction}\label{sec:intro}
	\begin{figure*}[t]
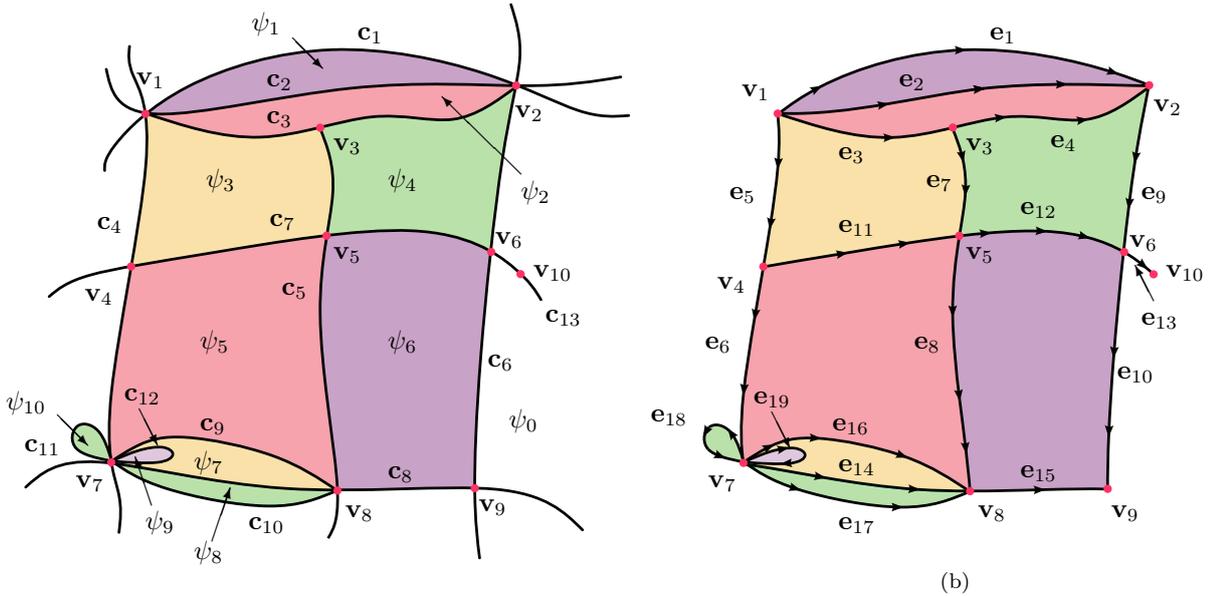
\centering
    	\subfloat {
    		\def\svgwidth{.9\columnwidth}
    		\input{pics/intro_example.eps_tex}
   		}%\quad
   		\subfloat[] {
    		\def\svgwidth{.9\columnwidth}
    		\input{pics/intro_example_paths.pdf_tex}
   		}
  		\caption{
  			An example of a geometric setting for the application of the region extraction algorithm presented in this work. The regions
  			$\psi_0,\dots,\psi_{10}$ to be extracted are marked with different colors in~(a). The edges as in Definition~\ref{def:edge} are
  			represented in~(b).
  		}
  		\label{fig:intro}
	\end{figure*}
	Given $n$ planar curves $\cc_1,\dots,\cc_n$, the region extraction problem consists in finding the regions bounded by the $n$ curves and their intersections, see Figure~\ref{fig:intro}. Such problem is a very common task in both Computer Science and Engineering with several applications in object recognition, motion analysis and stereopsis, among others.

	Most of the existing literature focuses on the case of a \emph{line drawing}, that is a class of pictorial data where the information is
	conveyed by the edges and the vertices of a planar graph, or to the case of intersection detection among two unstructured mesh. In the
	first setting, Jiang and Bunke~\cite{Jiang:1993:AOA} proposed an approach based on the arrangement of groups of edges of a
	planar graph in a counter-clockwise order. The ordered edges are then scanned linearly in order to create \emph{wedges} representing the
	area between two consecutive edges. The list of wedges is again sorted according to their vertices and appended together in order to
	build the different regions in $O(m\log m)$ steps, where $m$ is the number of edges of the line drawing. A similar algorithm
	was proposed by Shih~\cite{Shih:1989:ASA}, with the difference that the region extraction step of the algorithm is performed
	with a tagging process that results in assigning the same tag to the wedges belonging to the same region. This process allows to speed-up
	the number of steps to $O(m)$.
	In~\cite{Fan:1991:SEF} instead, the authors use an \emph{adjacency matrix} $M$ with entries
	\begin{equation}
		\label{eq:adjacency-matrix}
		M_{i,j} = \begin{cases}
			1, &\text{ if there is an edge between vertices $i$ and $j$}\\
			0, &\text{otherwise.}	
		\end{cases}
	\end{equation}
	and the clockwise angles between adjacent edges of the line drawing in order to walk the faces to be extracted in a counter-clockwise
	direction. Many different algorithms have instead been proposed for the identification of solid's faces in a wireframe,
	see~\cite{Dutton:1983:EIT,Brewer:1986:ACO,KUO:2001:AEO} and references therein, but all of these algorithms rely on the fact
	that wireframes can be represented as planar graphs.

	A planar graph is a pair $(V,E)$, where $V = \{\vv_1,\dots,\vv_n\}$ is a set of vertices and $E$ is a set of paired vertices called
	edges. Planar graphs can be represented using an adjacency matrix as in~\eqref{eq:adjacency-matrix} and therefore cannot fully depict a
	situation as the one represented in Figure~\ref{fig:intro}. If the boundary of the regions are indeed curvilinear segments, two vertices
	can be connected by more than one edge, see regions $\psi_1, \psi_2, \psi_7$, and $\psi_8$ in Figure~\ref{fig:intro}, a relation that
	cannot be represented using an adjacency matrix.

	A problem similar to region extraction arises in the context of mesh intersection methods as well. In order to intersect two
	unstructured meshes, Gander and Japhet~\cite{Gander:2013:A9P} proposed to use an approach that can be split in two steps.
	Given two triangular meshes $M_1$ and $M_2$ and two triangles $T_1\in M_1$ and $T_2\in M_2$, the first step consists in the
	identification of the intersections between $T_1$ and $T_2$. During this operation, a list of neighbor triangles of $T_2$ that intersect
	with $M_1$ is retrieved as well. The second step
	makes use of this additional information in order to extract the regions with an advancing front technique, defining the mesh
	intersections of $M_1$ and $M_2$. Another advancing front algorithm was proposed in~\cite{Lee:2004:FAF}. Lee and colleagues
	proposed an algorithm based on the construction of a background quadtree for the first mesh $M_1$ and a
	\emph{self-avoiding walk} for the second mesh $M_2$. Then, following the self-avoiding walk on the triangles of $M_2$, they use the local
	information of a triangle to generate the triangle-intersection set of the next one on the mesh $M_1$.
	In~\cite{Lohner:1988:SUD,Lohner:1989:ARF}, the author uses
	a background quadtree to search for nearby grid elements. Plimpton and colleagues~\cite{Plimpton:2004:APR}, instead proposed
	an approach based on recursive coordinate bisections for searching nearby grid elements.

	While region extraction algorithms for line drawings and mesh intersection algorithms have proved to be able to deal with many
	applications, the advances in Isogeometric Analysis (IGA) and other high-order methods call for a generalization of such
	problem to the case in which the regions to be extracted are not bounded by straight edges but by an arbitrary number of curved segments
	instead. In many applications, indeed, it is required to compute integrals involving spline quantities that are only piecewise
	polynomials with a finite order of continuity. A popular approach in IGA is to ignore the reduced inter-elements continuity of the
	integrand, defining the quadrature points over the domain of one of the splines. While this is an efficient approach, as no information
	about the second spline is required, it can result in large integration errors do to the reduced continuity of the integrand. Such errors
	can become even larger if the integration domain is not fully covered by the image of the two splines, as this can result in a
	discontinuity of the integrand.

	A method to overcome this issue is represented by the so called \emph{segment based integration schemes}. The idea behind these schemes
	is that the restrictions of the splines to their elements, are simple polynomials and therefore can be integrated using standard
	quadrature techniques.
	% While this procedure is trivial when integrating splines belonging to the same spline space, the situation gets
	% more complicated when the integrand is represented by a product of functions defined in different spline spaces. In these cases, it is
	% necessary to identify the regions in which all the splines involved are polynomials. While classical mesh intersection algorithms
	% usually deal with either triangular or quadrilateral regions, in this setting, such regions are bounded by curved B-spline segments.

	Segments based integration schemes are based on the approximation of the image of the regions in which the splines are simple
	polynomials. In the context of isogeometric mortar methods, for example, the product of two B-spline functions needs to be integrated
	over (part of) the boundary of a surface, or a solid, called \emph{interface}. Seitz and
	colleagues~\cite{Seitz:2016:IDM} considered four-sided (linear) quadrilateral approximations of the regions in which the
	splines are polynomials. Quadrilaterals corresponding to the different splines are coupled and projected to a common auxiliary plane,
	where their intersection is identified via a clip polygon algorithm.
	The obtained polygon is then triangulated and the quadrature points are
	defined inside each triangle. The obtained quadrature points are finally mapped back to one of the original curved regions in order
	to evaluate the desired integral. Hesh and Betsch~\cite{Hesch:2012:IAA} used instead a different approach in the context of
	domain decomposition methods. The authors proposed to project the control points of both splines on one of the two surfaces.
	Then the intersection of the edges of the two splines are identified and the set of obtained intersection points are triangulated.
	Quadrature points are again defined in each obtained triangle.

	While both approaches are quite efficient, their main drawback is that they cannot represent exactly the regions in which the involved
	splines are polynomials. If the regions involved are heavily concave a non-negligible amount of quadrature points will still be defined
	outside the desired region, affecting the quality of the numerical integration. In order to overcome this issue, an algorithm that
	automatically recognizes and extracts the curved regions in which the splines have the desired order of continuity is necessary. In
	this work we present a region extraction algorithm that can be applied to this context as well.
	This work is divided as follows. In Section~\ref{sec:the_region_extraction_algorithm} we give the definition of a \emph{curvilinear
	drawing}, a generalization of the line drawing treated in~\cite{Jiang:1993:AOA} and~\cite{Fan:1991:SEF}. Our region extraction algorithm
	is presented in the same section. In Section~\ref{sec:applications} we present several applications of the newly introduced algorithm in
	the context of mesh intersection for B-spline trivariate solids. Finally, some conclusive remarks are given in 
	Section~\ref{sec:conclusion}.

\section{The region extraction algorithm} % (fold)
	\label{sec:the_region_extraction_algorithm}
	Before introducing our algorithm, we first need to set some notations and definitions. The idea of line drawing can be easily generalized
	by the following.
	\begin{definition}\label{def:curve-drawing}
		Let $\cc_1, \dots, \cc_n$ be $n$ curves
		\[
			\cc_j\colon [0, 1]\rightarrow\RR^2, \qquad j = 1,\dots,n,
		\]
		such that $\cc_j\in C^2[0, 1]$ and their intersection points form a discrete set. Then the set of curves $\cc_1, \dots, \cc_n$,
		together with their intersection points, is said to be a curvilinear drawing.
	\end{definition}
	An example of curvilinear drawing with curves $\cc_1,\dots,$ $\cc_{13}$ is shown in Figure~\ref{fig:intro}. The requirement of the curves
	$\cc_1, \dots, \cc_n$ to be $C^2$ is not strictly necessary but is here enforced in order to simplify the discussion. The goal of this
	work is to extract the planar regions delimited by the (restriction of the) curves $\cc_1, \dots, \cc_n$. We remark
	that the condition about the discreteness of the intersection points is enforced in order to avoid curves that are partially coincident.
	For the purpose of this work we assume that we receive in input only the curves and that their restrictions need first to be found via
	curve-curve intersection.
	
	To this end we loop over each curve in order to find its intersections with all the others, that are the \emph{vertices} of the
	curvilinear drawing.
	\begin{definition}\label{def:vertex}
		Let $\vv\in\bar{\Omega}$ be the intersection of at least two curves in a curvilinear drawing. Then we say that $\vv$ is a \emph{vertex}
		of the curvilinear drawing.
	\end{definition}
	Considering again Figure~\ref{fig:intro}, the curvilinear drawing has ten vertices, denoted with $\vv_1,\dots,\vv_{10}$, respectively.
	
	\begin{table}\centering\footnotesize
		\begin{tabular}{c|c}
    		$i$ & $\mathcal{V}_c(\vv_i)$\\\hline
    		$1$ & $\{\cc_1,\cc_2,\cc_3,\cc_4\}$\\
    		$2$ & $\{\cc_1,\cc_2,\cc_3,\cc_6\}$\\
    		$3$ & $\{\cc_3,\cc_5\}$\\
    		$4$ & $\{\cc_4,\cc_7\}$\\
    		$5$ & $\{\cc_5,\cc_7\}$\\
    		$6$ & $\{\cc_6,\cc_7\}$\\
    		$7$ & $\{\cc_4,\cc_8,\cc_9,\cc_{10},\cc_{11},\cc_{12}\}$\\
    		$8$ & $\{\cc_5,\cc_8,\cc_9,\cc_{10}\}$\\
    		$9$ & $\{\cc_6,\cc_8\}$\\
    		$10$& $\{\cc_6,\cc_{13}\}$
  		\end{tabular}
  		\qquad
  		\begin{tabular}{c|cc}
    		$j$ & $\mathcal{C}_v(\cc_j)$\\\hline
    		$1$ & $\{\vv_1,\vv_2\}$\\
    		$2$ & $\{\vv_1,\vv_2\}$\\
    		$3$ & $\{\vv_1,\vv_3,\vv_2\}$\\
    		$4$ & $\{\vv_1,\vv_4,\vv_7\}$\\
    		$5$ & $\{\vv_3,\vv_5,\vv_8\}$\\
    		$6$ & $\{\vv_2,\vv_6,\vv_9, \vv_{10}\}$\\
    		$7$ & $\{\vv_4,\vv_5,\vv_6\}$\\
    		$8$ & $\{\vv_7,\vv_8,\vv_9\}$\\
    		$9$ & $\{\vv_7,\vv_8\}$\\
    		$10$ & $\{\vv_7,\vv_8\}$\\
    		$11$ & $\{\vv_7,\vv_7\}$\\
    		$12$ & $\{\vv_7,\vv_7\}$\\
    		$13$ & $\{\vv_{10}\}$\\
  		\end{tabular}
  		\caption{Maps $\mathcal{V}_c(\vv_i)$, $i = 1,\dots,9$ and $\mathcal{C}_v(\cc_j)$, $j = 1,\dots,13$, for the curves-vertices identification of the curvilinear drawing in Figure~\ref{fig:intro}. The list of vertices in $\mathcal{C}_v(\cc_j)$ follows the order inherited by the parameterization of $\cc_j$.}
  		\label{tab:maps}
	\end{table}
	Let $C$ and $V$ be the lists containing the curves and the currently identified vertices of a line drawing. Every time a new vertex $\vv$
	has	been identified we add it to the list $V$ and update two maps. The first map
	\[
		\mathcal{V}_c\colon V\rightarrow C
	\]
	keeps track of the curves defining each vertex in $V$. The second map
	\[
		\mathcal{C}_v\colon C\rightarrow V
	\]
	helps to identify the vertices lying on the same curve. For each curve $\cc$ we store the vertices lying on $\cc$ following the order
	given by the corresponding parameters. Assuming that all the curves in Figure~\ref{fig:intro}~(a) are parameterized as shown
	in Figure~\ref{fig:intro}~(b), Table~\ref{tab:maps} shows the corresponding maps $\mathcal{V}_c$ and $\mathcal{C}_v$ for each vertex and
	each curve. Notice that $\mathcal{C}_v(\cc_{11})$ and $\mathcal{C}_v(\cc_{12})$ contain twice the same vertex, indicating that $\cc_{11}$
	and $\cc_{12}$ are closed curves that should be traversed in both directions. When all the vertices of the curvilinear drawing
	have been identified, $\mathcal{V}_c$ and $\mathcal{C}_v$ allow us to find the corresponding \emph{edges}.
	\begin{definition}\label{def:edge}
		Let $\vv_i, \vv_j$, and $\cc$ be respectively two vertices and a curve in a curvilinear drawing such that there exist two parameters
		$t_i, t_j, t_i < t_j$, with $\cc(t_i) = \vv_i$ and $\cc(t_j) = \vv_j$. Then the restriction $\ee = \restr{\cc}{[t_i, t_j]}$ is an
		(oriented) edge of the curvilinear drawing.
	\end{definition}

	The edges of the curvilinear drawing in Figure~\ref{fig:intro}~(a), together with their orientation, are shown in
	Figure~\ref{fig:intro}~(b). Using a rather standard notation, here and in the following, we denote with $-\ee_i$ the edge $\ee_i$
	traversed following its inverse parameterization.
	With these definitions at hand, a region of the curvilinear drawing can be represented as a \emph{closed trail}.
	\begin{definition}
		A closed trail is a sequence of pairs of vertices and edges $$((\vv_{i_1}, \ee_{i_1}), (\vv_{i_2}, \ee_{i_2}), \dots, (\vv_{i_m},
		\ee_{i_m}))$$ such that each $\ee_{i_j}$ is an edge between $\vv_{i_j}$ and $\vv_{i_{j+1}}$ and $\ee_{i_m}$ is an edge between
		$\vv_{i_m}$ and $\vv_{i_1}$.
	\end{definition}
	Note that the edges in this definition are oriented. As an example, let us consider again the drawing in
	Figure~\ref{fig:intro}. The region $\psi_3$ can be represented as
	\[
		\psi_3 = ((\vv_1,\ee_5),(\vv_4,\ee_{11}), (\vv_5, -\ee_7), (\vv_3,-\ee_3)).
	\]
	The main advantage of this notation is that one can easily represent regions that would be difficult to represent with, e.g., vertex
	based representations. Examples of these regions are the ones bounded by closed edges as $\psi_7$ and $\psi_9$ in Figure~\ref{fig:intro},
	which are
	\[
		\psi_7 = ((\vv_7,\ee_{14}),(\vv_8,\ee_{16}),(\vv_7,\ee_{19}))
	\]
	and
	\[
		\psi_9 = ((\vv_7,-\ee_{19})),
	\]
	respectively.

	As remarked in Section~\ref{sec:intro}, a line drawing can be efficiently represented as the adjacency matrix $M$
	in~\eqref{eq:adjacency-matrix}. Since this is not the case for a curvilinear drawing with loops (e.g. $\psi_9$ and $\psi_{10}$ in
	Figure~\ref{fig:intro}) and multiedges (e.g. $\ee_1, \ee_2, \ee_{14},\ee_{16}$, and $\ee_{17}$ in Figure~\ref{fig:intro}), we store its
	connectivity as lists of unvisited paths, one for each vertex of the curvilinear drawing.
	Given a vertex $\vv$, its unvisited path list $\Pi$ contains the (oriented) edges that originate in $\vv$.
	For instance, the initial list of unvisited paths for vertex $\vv_5$ in Figure~\ref{fig:intro} is given by
	\begin{equation}
		\label{eq:unvisited-list-5}
		\Pi_5 = \{-\ee_7,\ee_{12},-\ee_{11},\ee_8\}.
	\end{equation}
	Using the maps $\mathcal{V}_c$ and $\mathcal{C}_v$ in Table~\ref{tab:maps}~(a) and~(b), we can easily build the initial lists of unvisited paths for the drawing in Figure~\ref{fig:intro}, see Table~\ref{tab:unvisited_paths}.
	By Definition~\ref{def:edge}, the list corresponding to vertex $\vv_{10}$ contains only the path $-\ee_{13}$, as curve $\cc_{13}$ in Figure~\ref{fig:intro} has no further vertices. Vertices like $\vv_{10}$ are called \emph{dangling nodes} and are easily recognizable as vertices having only one open edge among their unvisited paths.

	Our region extraction algorithm is described in Algorithm~\ref{alg:region-extraction}.
	\begin{algorithm}[t]
  		\caption{ExtractRegions : Extracts the regions from the curvilinear drawing}
  		\label{alg:region-extraction}
  		\hspace*{\algorithmicindent} \textbf{Input} Vertices list $V = \{\vv_1,\dots,\vv_m\}$;\\
  		\hspace*{\algorithmicindent} \textbf{Input} Unvisited paths lists, $\Pi = \{\Pi_1,\dots,\Pi_m\}$;\\
  		\hspace*{\algorithmicindent} \textbf{Output} List $\Psi$ of the extracted regions;
  		\begin{algorithmic}[1]
    		\State $\Psi\coloneqq\emptyset$;    /* Initialize list of regions. */ 
    		\State PurgeDanglingNodes($V$, $\Pi$);\label{line:purge_dangling_nodes}
    		\For{$\vv_{i_1}\in V$}\label{line:starting-vertex}
    			\While{$\Pi_{i_1}\ne\emptyset$}
    				\State $\ee_{i_1} \gets$ path in $\Pi_{i_1}$; \label{line:first-edge}
    				\State $\vv_{i_2} \gets$ endpoint of $\ee_{i_1}$;
    				\State $\ee_{i_2} = \text{arg}\displaystyle\max_{\ee\in\Pi_{i_2}}\measuredangle[\ee\ee_{i_1}]$;\label{line:pick-a-path}
    				\State $\psi \coloneqq \{(\vv_{i_1}, \ee_{i_1})\}$;\label{line:initialize-current-face}    /* Initialize region. */
    					\State $j \coloneqq 2$;
    					\While{$\ee_{i_j} \ne \ee_{i_0}$}
	    					\State $\psi = \psi \cup \{(\vv_{i_j}, \ee_{i_j})\}$;
    						\State $\Pi_{i_j} = \Pi_{i_j} \setminus \{\ee_{i_j}\}$;    /* Update $\Pi_{i_j}$. */\label{line:update-Pi_i}
    						\State $\vv_{i_{j + 1}} \gets$ endpoint of $\ee_{i_j}$;    /* Find vertex. */\label{line:set-next-seed-1}
    						\State $\ee_{i_{j + 1}} = \text{arg}\displaystyle\max_{\ee\in\Pi_{i_{j + 1}}}\measuredangle[\ee\ee_{i_j}]$;\label{line:set-next-seed-2}
    						\State $j = j + 1$;
    					\EndWhile
    					\State $\Pi_{i_1} = \Pi_{i_1} \setminus \{\ee_{i_1}\}$;    /* Update $\Pi_{i_0}$. */\label{line:update-Pi_0}
    					\State $\Psi = \Psi\cup\psi$;
    				\EndWhile
    		\EndFor
    		\State \Return $\Psi$;
  		\end{algorithmic}
	\end{algorithm}
	Line~\ref{line:purge_dangling_nodes} calls a simple routine that takes care of recognizing and discarding the dangling nodes in the
	curvilinear drawing, see Algorithm~\ref{alg:purge_dangling_nodes}. We remark that, assigning an orientation to the edges, we implicitly 
	duplicate the edges as each edge belongs to exactly two vertices and is hence traversed once in each direction in
	Algorithm~\ref{alg:region-extraction}.

	Given a curvilinear drawing, the main idea behind Algorithm~\ref{alg:region-extraction} is to construct a corresponding \emph{rotation
	system}~\cite{Heffter:1891:UDB} and to use the obtained topological information to extract the regions bounded by its halfedges. For any 
	vertex $\vv$ in a given multigraph, a rotation system associates an ordering to the edges originating in $\vv$. Such ordering
	is implicitly defined by the orientation of the surface in which the multigraph lies. Assigning a counter-clockwise ordering to all the
	edges that originate in $\vv$, it is then possible to extract the regions around that vertex, by simply walking all closed
	trails starting and ending in $\vv$.

	\begin{table}\centering\footnotesize\onehalfspacing
		\begin{tabular}{c|c}
    	$i$ & $\Pi_i$\\\hline
    	$1$ & $\{\ee_1,\ee_2,\ee_3,\ee_5\}$\\
    	$2$ & $\{-\ee_1,-\ee_2,-\ee_4,\ee_9\}$\\
    	$3$ & $\{-\ee_3,\ee_4,\ee_7\}$\\
    	$4$ & $\{-\ee_5,\ee_6,\ee_{11}\}$\\
    	$5$ & $\{-\ee_7,\ee_8,-\ee_{11},\ee_{12}\}$\\
    	$6$ & $\{-\ee_9,\ee_{10},\ee_{13},-\ee_{12}\}$\\
    	$7$ & $\{-\ee_6,\ee_{14},\ee_{16},\ee_{17},\pm\ee_{18},\pm\ee_{19}\}$\\
    	$8$ & $\{-\ee_8,-\ee_{14},\ee_{15},-\ee_{16},-\ee_{17}\}$\\
    	$9$ & $\{-\ee_{10},-\ee_{15}\}$\\
    	$10$ & $\{-\ee_{13}\}$\\
  	\end{tabular}
  	\caption{List of unvisited paths for the vertices $\vv_1,\dots,\vv_{10}$ in Figure~\ref{fig:intro}.}
  	\label{tab:unvisited_paths}
	\end{table}
	
	In order to clarify how the region extraction algorithm works, let us consider again the drawing in Figure~\ref{fig:intro}. Without loss
	of generality, we assume that $\vv_{i_1} = \vv_5$ in Line~\ref{line:starting-vertex} of Algorithm~\ref{alg:region-extraction}. At the
	beginning of the algorithm the list of unvisited paths for $\vv_5$ is given in~\eqref{eq:unvisited-list-5}.
	The algorithm assigns $-\ee_7$ to $\ee_{i_1}$ in Line~\ref{line:first-edge} and picks its end-point $\vv_3$. Therefore it sets $\vv_{i_2} = \vv_3$ and looks in the corresponding list $\Pi_3$, see Table~\ref{tab:unvisited_paths}, for the edge that forms the maximum counter-clockwise angle with $\ee_{i_1}$. To this end we use the following definition of angle between two edges.
	\begin{definition}
		Let $\vv$ be a vertex and let $\bar\ee$ and $\tilde\ee$ be two edges originating in $\vv$. Then the angle in $\vv$ between edges $\bar\ee$ and $\tilde\ee$ is given by the counterclockwise angle between $\bar\tl(\vv)$ and $\tilde\tl(\vv)$, where $\bar\tl(\pp)$ and $\tilde\tl(\pp)$ denote the tangents of $\bar\ee$ and $\tilde\ee$ in a point $\pp\in\bar\Omega$, respectively. The orientation of the tangents are chosen so that they always point toward the interior of the corresponding edges.
	\end{definition}
	Therefore we set $\ee_{i_2} = -\ee_3$ and we initialize the region $\psi = \{(\vv_5,-\ee_7)\}$. 
	Since $\ee_{i_2}\ne\ee_{i_0}$, we add the pair $(\vv_3,-\ee_3)$ to $\psi$ and we remove $-\ee_3$ from the list of unvisited paths of $\vv_3$. The algorithm then sets $\vv_{i_3} = \vv_1$ and $\ee_{i_3} = \ee_5$, since $\ee_5$ is the edge forming again the maximum angle with $\ee_{i_2}$. The pair $(\vv_1,\ee_5)$ is added to $\psi$ and $\ee_5$ is removed from the list of unvisited paths of $\vv_1$. Since the end-point of $\ee_5$ is $\vv_4$, Lines~\ref{line:set-next-seed-1} and~\ref{line:set-next-seed-2} set $\vv_{i_4} = \vv_4$ and $\ee_{i_4} = \ee_{11}$. 
	The list $\psi$ is then updated with the pair $(\vv_4,\ee_{11})$ and $\ee_{11}$ is removed by $\Pi_4$. Now the end-point of $\ee_{11}$ is the initial vertex $\vv_5$ and the edge with maximal angle is $-\ee_7$ that is still in the list of unvisited paths $\Pi_{i_1} = \Pi_5$.
	Therefore, we set $\vv_{i_5} = \vv_5$ and $\ee_{i_5} = -\ee_7$. Since $\ee_{i_5} = \ee_{i_0}$, the algorithm terminates the inner loop, deletes $-\ee_7$ from $\Pi_5$ and adds the identified face $\psi_3$ to the list $\Psi$ of extracted regions. After the identification of $\psi_3$, the list of unvisited paths for $\vv_5, \vv_3, \vv_1$, and $\vv_4$ are
	\begin{equation*}
		\begin{aligned}
			\Pi_5 &= \{\ee_8,-\ee_{11},\ee_{12}\},\\
			\Pi_3 &= \{\ee_4,\ee_7\},\\
			\Pi_1 &= \{\ee_1,\ee_2,\ee_3\},
		\end{aligned}
	\end{equation*}
	and
	\[
		\Pi_4 = \{-\ee_5,\ee_6\},
	\]
	respectively.
	Algorithm~\ref{alg:region-extraction} proceeds then with the extraction of the regions $\psi_6$, $\psi_4$ and $\psi_5$ corresponding to $\vv_{i_1} = \vv_5$ and $\ee_{i_1} = \ee_8, \ee_{12}$ and $-\ee_{11}$, respectively.
	Once all the regions surrounding $\vv_5$ have been extracted, Algorithm~\ref{alg:region-extraction} finds a different vertex with a non empty unvisited paths list, until all regions in the curvilinear graph have been recognized.
	
	Let us remark some features about Algorithm~\ref{alg:region-extraction}. An edge in a curvilinear drawing is adjacent to
	exactly two regions, one for each direction in which we can traverse it, while any region adjacent to a vertex is always bounded by at
	least one of the edges originating at said vertex. Therefore, Lines~\ref{line:pick-a-path} and~\ref{line:set-next-seed-2} in
	Algorithm~\ref{alg:region-extraction} are always guaranteed to succeed. On the other hand, for any closed region in a curvilinear graph
	there exists at least one edge that is adjacent to it and hence the algorithm is guaranteed to extract all of them. Finally, the
	algorithm is guaranteed to terminate when all the edges of the drawing have been traversed in both directions.

	The detection of the maximal angle extensively used in Algorithm~\ref{alg:region-extraction} can be complicated by the fact that two or 
	more edges can have same tangent directions in a vertex. In order to simplify the discussion, we did not consider this case in the 
	description of the region extraction algorithm. Nevertheless, these corner cases can be treated by considering beforehand the curvature
	of the edges having same tangents in a single vertex.
	% as follow. Let $\vv$ be a vertex obtained as the intersection of two curves with the same tangent, $\vv_1,\dots,\vv_s$ be the end-points of the edges starting in $\vv$ and
	% \[
	% 	d = \min_{i = 1,\dots,s}\norm{\vv - \vv_i}.
	% \]
	% Then we consider a circle centered in $\vv$ with radius $d/2$ and its intersections with the edges originating at $\vv$, $\pp_1,\dots,\pp_s$.
	% Since each point $\pp_i$ corresponds to a different edge $\ee_i$ starting in $\vv$, we use the counterclockwise angle between the segments $\vv\pp_i$, $i = 1,\dots,s$, in order to identify the edge to be traversed, see Figure~\ref{fig:cotangent-curves}.
	% \begin{figure}[t]\centering
 %    	\mbox {
 %    		\def\svgwidth{.9\columnwidth}
 %    		\input{pics/coplanar_curves.eps_tex}
 %   		}
 %  		\caption{Construction for the tangent identification for curves with the same tangents at a vertex $\vv$. The arrows denote the directions to be used as tangents for the corresponding edges.}
 %  		\label{fig:cotangent-curves}
	% \end{figure}
	%
	\begin{algorithm}[t]
  		\caption{PurgeDanglingNodes : Recognizes and delete the dangling nodes and the relative edges.}
  		\label{alg:purge_dangling_nodes}
  		\hspace*{\algorithmicindent} \textbf{Input} Vertices list $V = \{\vv_1,\dots,\vv_m\}$;\\
  		\hspace*{\algorithmicindent} \textbf{Input} Unvisited paths lists, $\Pi = \{\Pi_1,\dots,\Pi_m\}$;\\
  		\hspace*{\algorithmicindent} \textbf{Output} Updated $V$ and $\Pi$;
  		\begin{algorithmic}[1]
    		\State $L = m$;    /* Initial number of vertices. */
    		\For{$\vv_{i}\in V$}
    			\If{$\norm{\Pi_i} = 1$}    /* $\vv_i$ is a dangling node. */
    				\State $\ee \gets$ edge in $\Pi_i$;
    				\State delete $\ee$ from all unvisited paths lists in $\Pi$;
    				\State $V = V\setminus\{\vv_i\}$;
    				\State $\Pi = \Pi\setminus{\Pi_i}$
    			\EndIf
    		\EndFor
    		\If{$\norm{V} < L$}
    			\State PurgeDanglingNodes($V$,$\Pi$);
    		\EndIf
  		\end{algorithmic}
	\end{algorithm}

	Among the regions extracted by Algorithm~\ref{alg:region-extraction} there is also the external unbounded region $\psi_0$ in 
	Figure~\ref{fig:intro}. Whether similar regions need to be purged or not can be application dependent but it is possible to easily 
	identify them as the regions for which the list of consecutive vertices and edges follows a clockwise direction.
	In order to avoid further operations, it is possible to compute the angles between the edges in the interval $[-\pi,\pi]$ and to sum them
	up. Once a region has been extracted, it can be purged if the sum of the angles is positive.
% section the_region_extraction_algorithm (end)
\section{Applications} % (fold)
	\label{sec:applications}
	As mentioned in Section~\ref{sec:intro}, Algorithm~\ref{alg:region-extraction} can be used as a tool to numerically compute integrals
	involving spline functions with a finite order of continuity and their product. Usually, standard quadrature rules provide accurate
	approximations of the
	integrals only if the integrand can be well approximated by polynomials. In order to numerically compute the integral of splines
	functions it is therefore necessary to identify the sub-regions of their domain in which they are polynomials, and apply the quadrature
	rule separately in each sub-region. Algorithm~\ref{alg:region-extraction} represents a robust tool to automatically recognize and extract
	such regions in the case of an integrand represented by a product of spline functions.

	In this section we further develop this idea in the context of trivariate B-spline functions defined over three-dimensional solids.
	Despite we here restrict ourselves to the case of $n = 2$ intersecting solids, the methods described in the rest of this section can be
	easily generalized to the case $n \ge 3$, using the domain decomposition presented in~\cite{Antolin:2020:OMI}.

	Let us denote with $B_{i,\dd} = B_{i, \dd, \tl}$ the $i$-th tensor-product trivariate B-spline basis function of degrees
	$\mathbf{d} = (d_u, d_v, d_w)$ with knot vectors $\tl = (\tl_u, \tl_v, \tl_w)$ and with $\bbs_{\dd,\tl}$ the corresponding linear space
	spanned by these basis functions. Moreover, for every bounded domain $\Omega\in\RR^s$, let us denote with $\bar\Omega$ its closure.
	\begin{figure*}[t]\centering
  	\mbox {
   		\subfloat[] {
     		\includegraphics[height = 4.6cm, valign = c]{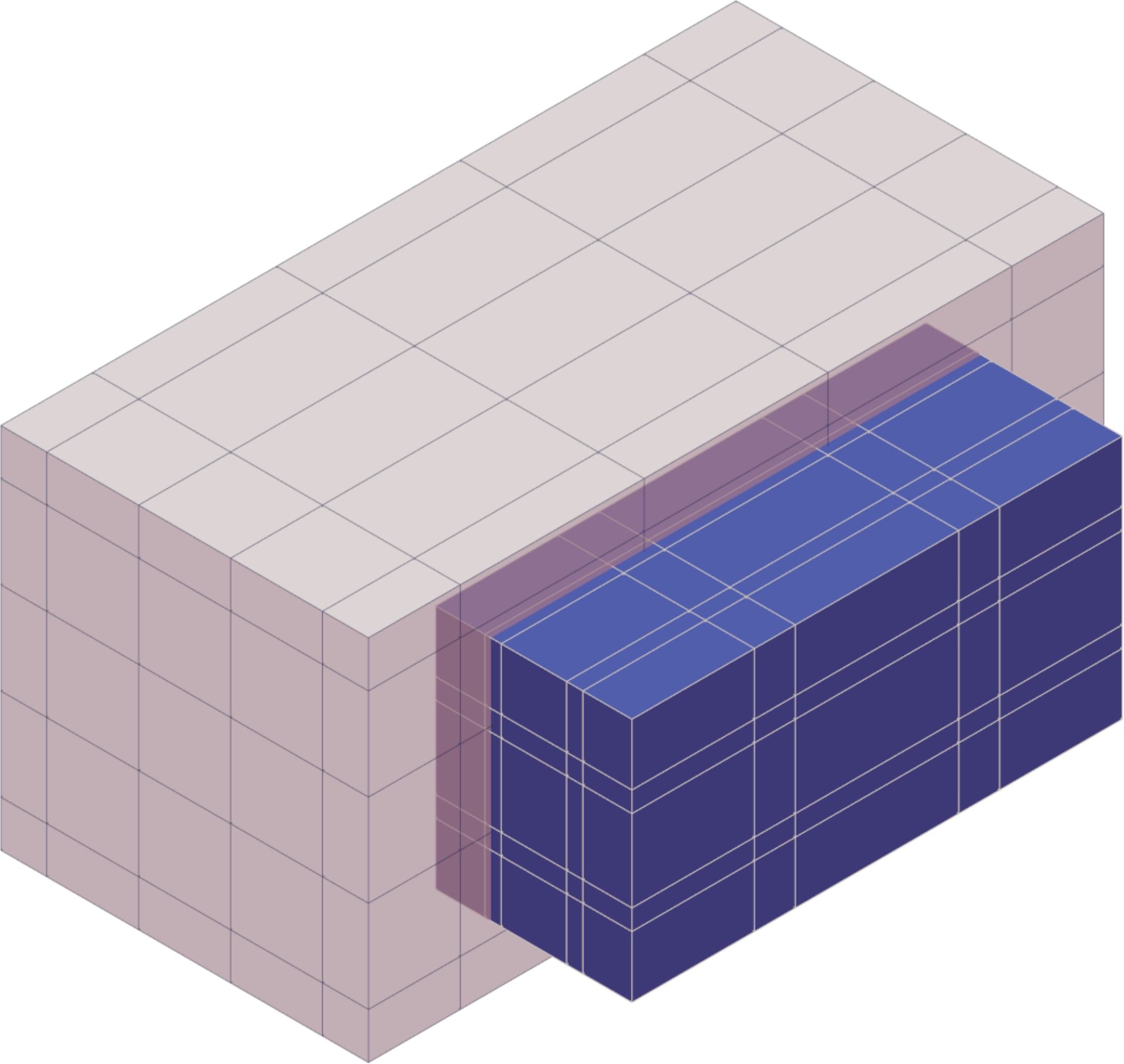}
     		\begin{picture}(0,0)
       			\put(-120, 45){$\Omega^*_1$}
       			\put(-40, -55){$\Omega^*_2$}
   			\end{picture}
   		}\qquad
   		\subfloat[] {
   			\includegraphics[height = 4.6cm, valign = c]{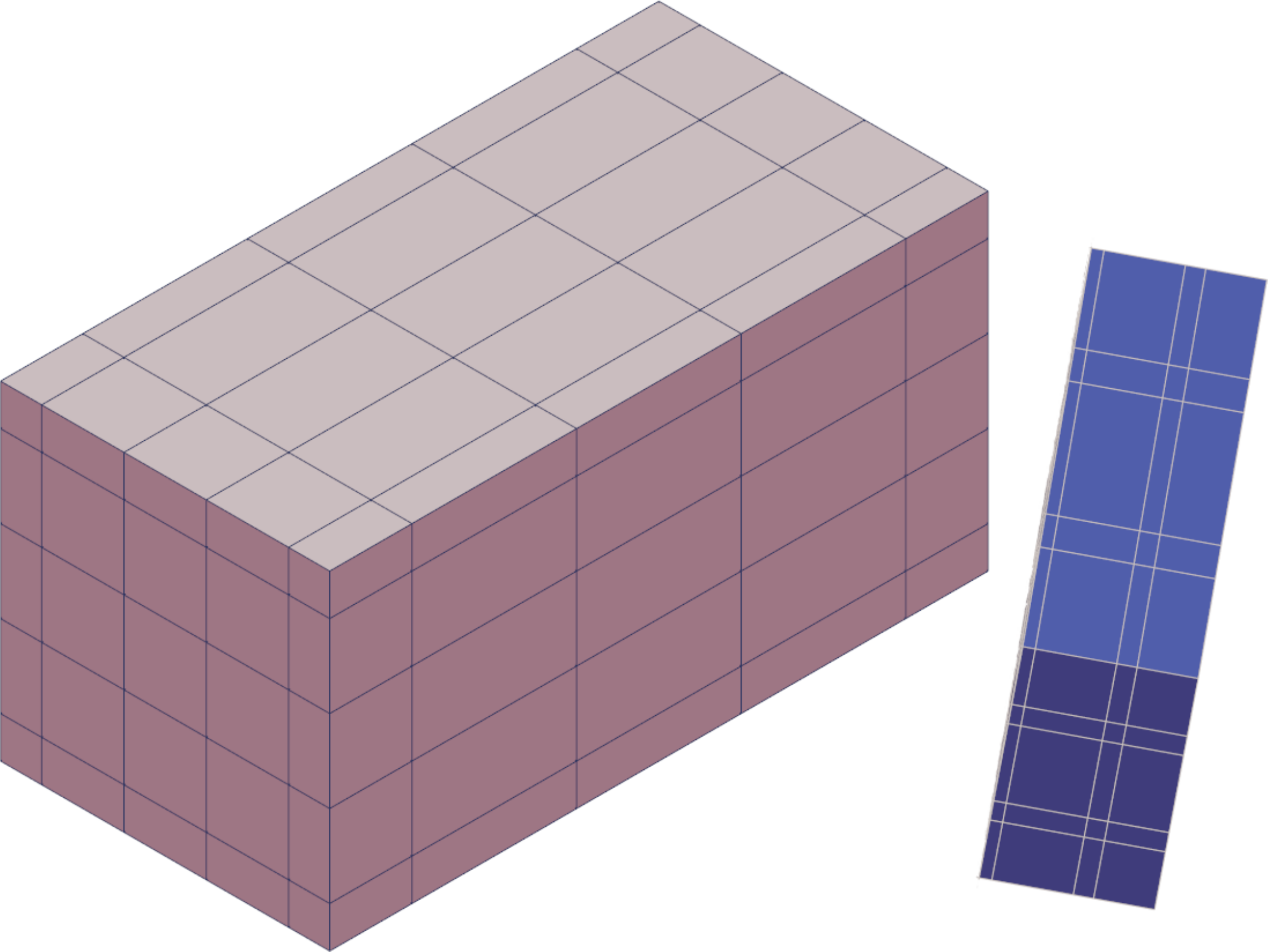}
   			\begin{picture}(0,0)
       			\put(-155, 45){$\Omega_1$}
       			\put(-60, -55){$\Omega_2$}
   			\end{picture}
   		}\qquad
   		\subfloat[] {
   			\includegraphics[height = 4.6cm, valign = c]{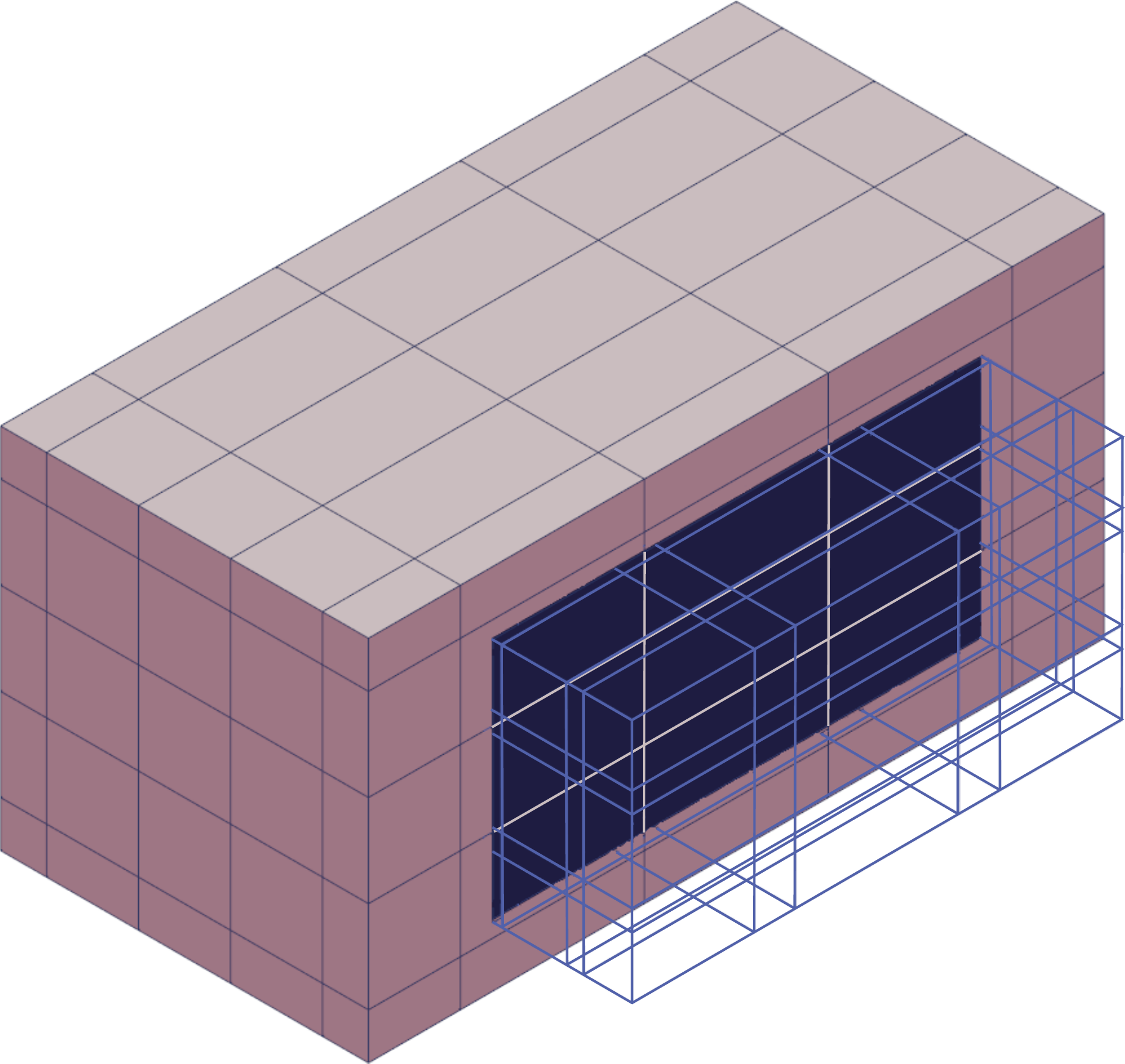}
     	  	\begin{picture}(0,0)
       			\put(-120, 45){$\Omega_1$}
       			\put(-40, -55){$\Omega_2$}
 	     			\put(-17, -45){$\Gamma_{1,2}$}
       			\put(-18, -37){\vector(-1, 1){22}}
   			\end{picture}
   		}
 		}
 		\caption{
 			Two domains $\Omega^*_1$ and $\Omega^*_2$ are shown in (a), while their partitions $\Omega_1$ and $\Omega_2$ as in
 			Equation~\eqref{eq:data-structure} are visible in (b). The local interface $\Gamma_{1,2}$ defined in Equation~\eqref{eq:interface},
 			together with the control meshes of $\Omega_1$ and $\Omega_2$, is shown in (c).
 		}
 		\label{fig:interface}
	\end{figure*}
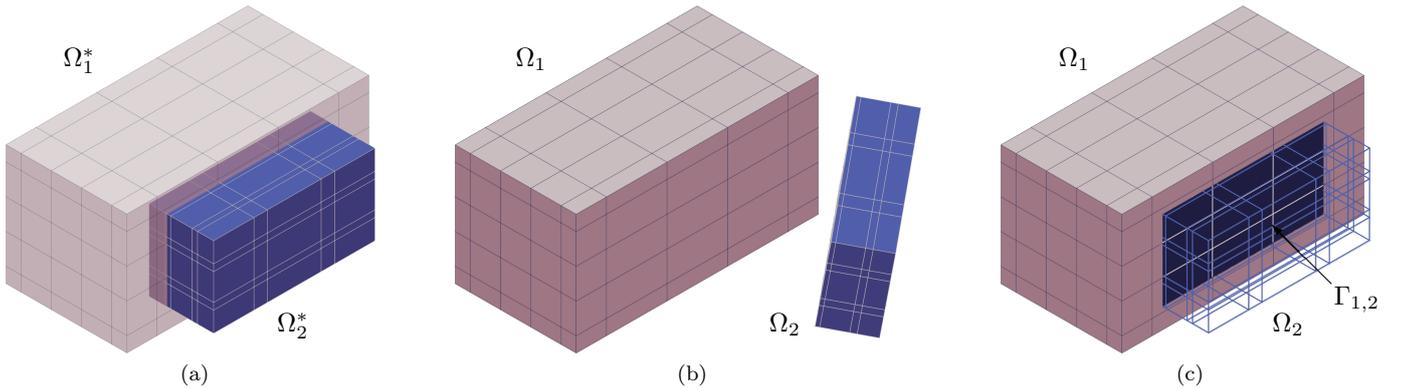
	
\begin{figure*}\centering
		\def\svgwidth{1.2\columnwidth}
		\input{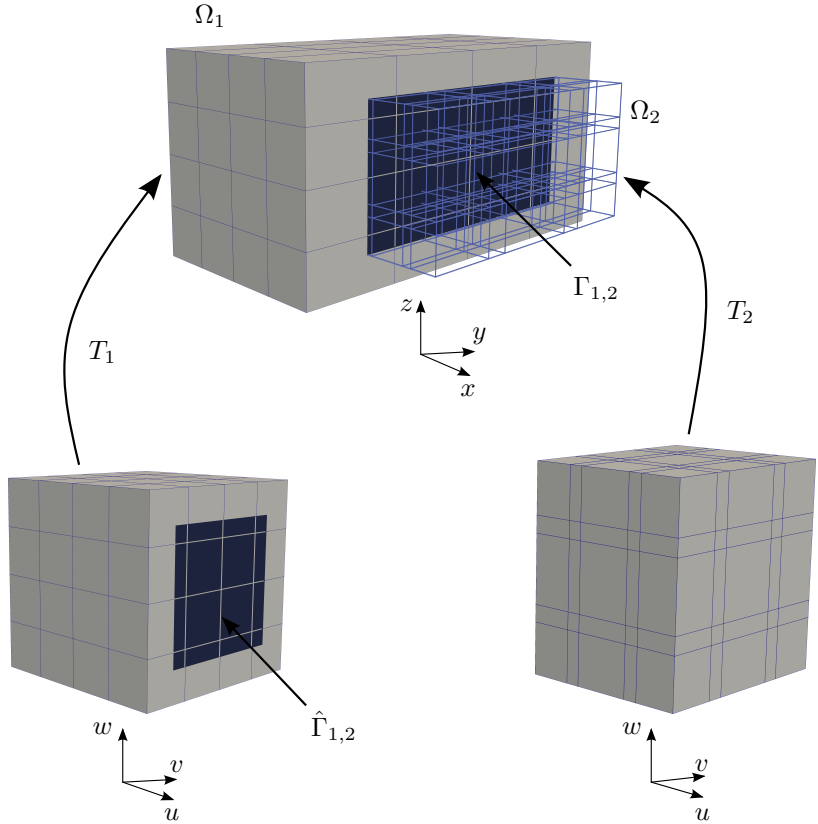}
 		\caption{A partition of the union of two domains $\Omega^*_1$ and $\Omega^*_2$. The interface $\Gamma_{1,2}$ and its preimage
 			$\hat\Gamma_{1,2}$ are highlighted in blue in both the Euclidean space and the parametric domain of $T_1$.}
 		\label{fig:maps}
	\end{figure*}

	Let us formalize our setting. Let $\Omega\subset\RR^3$ be a connected, bounded domain such that there exist two possibly
	overlapping domains $\Omega^{*}_1$ and $\Omega^{*}_2$ such that $\Omega = \Omega^{*}_1\cup\Omega^{*}_2$. The domain $\Omega$ can be
	partitioned as $\{\Omega_1,\Omega_2\}$, where
	\begin{equation}
		\label{eq:data-structure}
		\begin{aligned}
			\Omega_1 &= \Omega^{*}_1,\\
			\Omega_2 &= \Omega^*_2 \setminus \Omega^{*}_1.
		\end{aligned}
	\end{equation}
	We define the \emph{interface} $\Gamma_{1, 2}$ as
	\begin{equation}
		\label{eq:interface}
  	\Gamma_{1, 2} = \partial \Omega^{*}_1\cap \bar\Omega_2.
	\end{equation}
	Figure~\ref{fig:interface} shows an example of such a domain decomposition. In Figure~\ref{fig:interface}~(a) the two overlapping domains
	$\Omega^{*}_1$ and $\Omega^{*}_2$ are shown. Figure~\ref{fig:interface}~(b) shows the partition $\{\Omega^*_1,\Omega^*_2\}$ as defined
	in~\eqref{eq:data-structure}, while (c) shows the interface $\Gamma_{1,2}$ as defined in~\eqref{eq:interface}. Let
	$T_1\in\bbs_{\dd_1,\tl_1}$ and $T_2\in\bbs_{\dd_2,\tl_2}$ be two trivariate B-spline parameterizations of $\bar{\Omega}^{*}_1$ and
	$\bar{\Omega}^{*}_2$, respectively, see Figure~\ref{fig:maps}. Here and in the rest of this section we finally denote with
	$\hat\Gamma_{1,2}$ the preimage of $\Gamma_{1,2}$ in the parametric space of $T_1$.
	In this section we are going to face three different problems. In Section~\ref{sub:computation_of_integrals_over_planar_domains} we are
	going to use Algorithm~\ref{alg:region-extraction} in order to create a suitable quadrature rule for the approximate computation of
	integrals defined over the interface of two solids. Two numerical experiments are carried out, one concerning the integration of a
	smooth function and one regarding the integration of the product of splines defined in different spline spaces.

	In Sections~\ref{sub:weak_continuity} and~\ref{ssub:weak_continuity_not_conf} we apply the quadrature rule described in 
	Section~\ref{sub:computation_of_integrals_over_planar_domains} in order to enforce weak continuity constraints to volumetric objects.
	The constraints are imposed differently, depending on the relative position of the two objects. In one of the cases we are able to
	reproduce the deformation of $\Omega^*_2$ in the spline space $\mathbb{S}_{\dd_1}$, provided that the latter is sufficiently refined.
	In the second case, the appearance of oscillations does not allow us to obtain the same results but we approximate the deformation of
	$\Omega^*_2$ using a convolution based strategy.

	Finally, in Section~\ref{sub:poisson}, we use our algorithm to solve the Poisson's problem for two bodies in a contact position, using
	a mortar-like approach as described in~\cite{Brivadis:2015:IMM}.
	\subsection{Precise computation of integrals over the interface}
		\label{sub:computation_of_integrals_over_planar_domains}
\begin{figure*}[t]\centering
      \mbox {
    		\subfloat[] {
    			\includegraphics[height = 4.55cm, valign = c]{pics/img/quadrature_mesh_1.pdf}
    			\begin{picture}(0,0)
      			\put(-140, 78){$\Omega_1$}
  	  			\put(-33, -50){$\Omega_2$}
    			\end{picture}
    		}\quad
    		\subfloat[] {
    			\includegraphics[height = 4.55cm, valign = c]{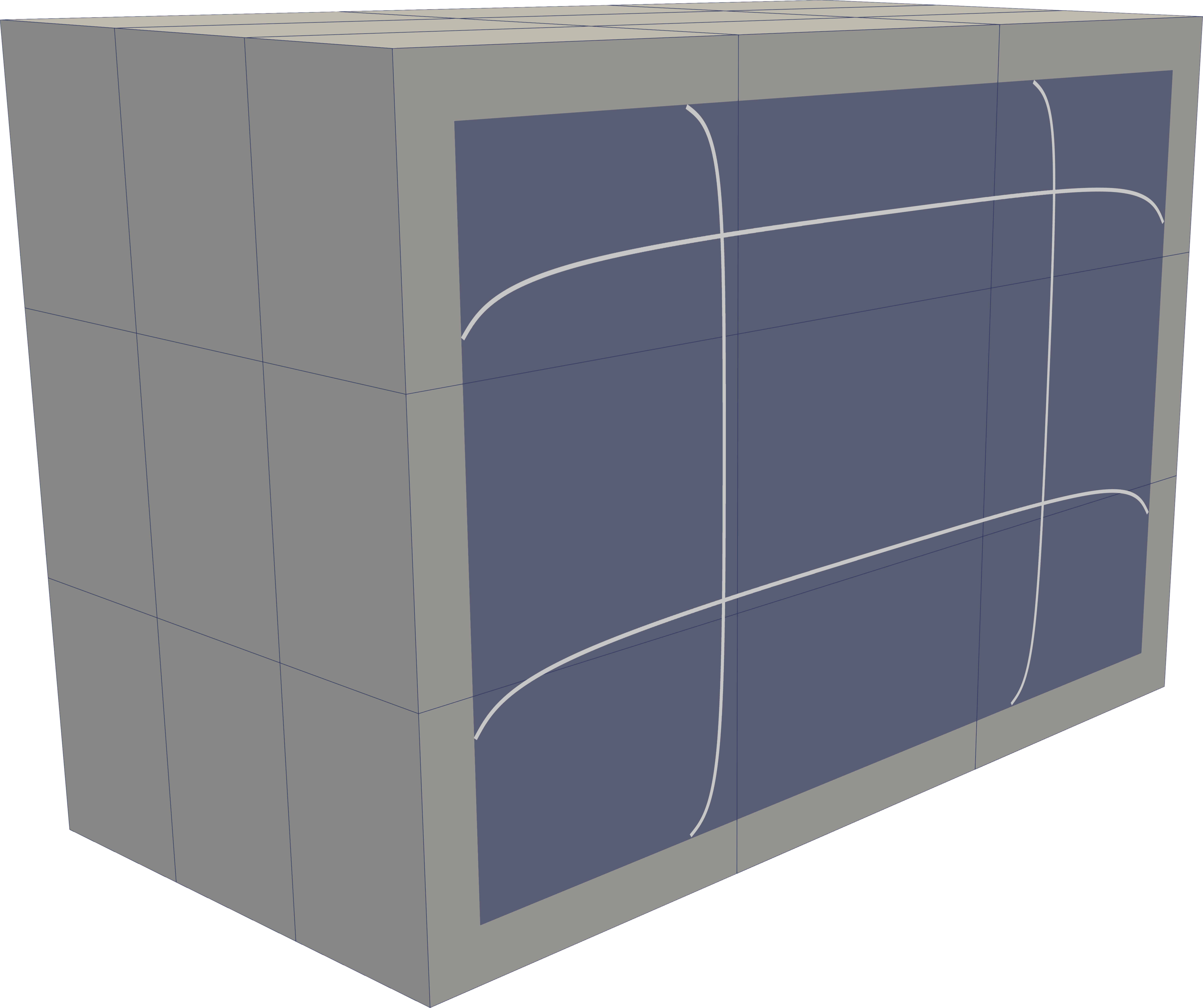}
    			\begin{picture}(0,0)
      			\put(-120, 78){$\Omega_1$}
  	  			\put(-30, -57){$\Gamma_{1,2}$}
      			\put(-27, -47){\vector(-1, 2){15}}
    			\end{picture}
    		}\quad
    		\subfloat[] {
    			\includegraphics[height = 4.55cm, valign = c]{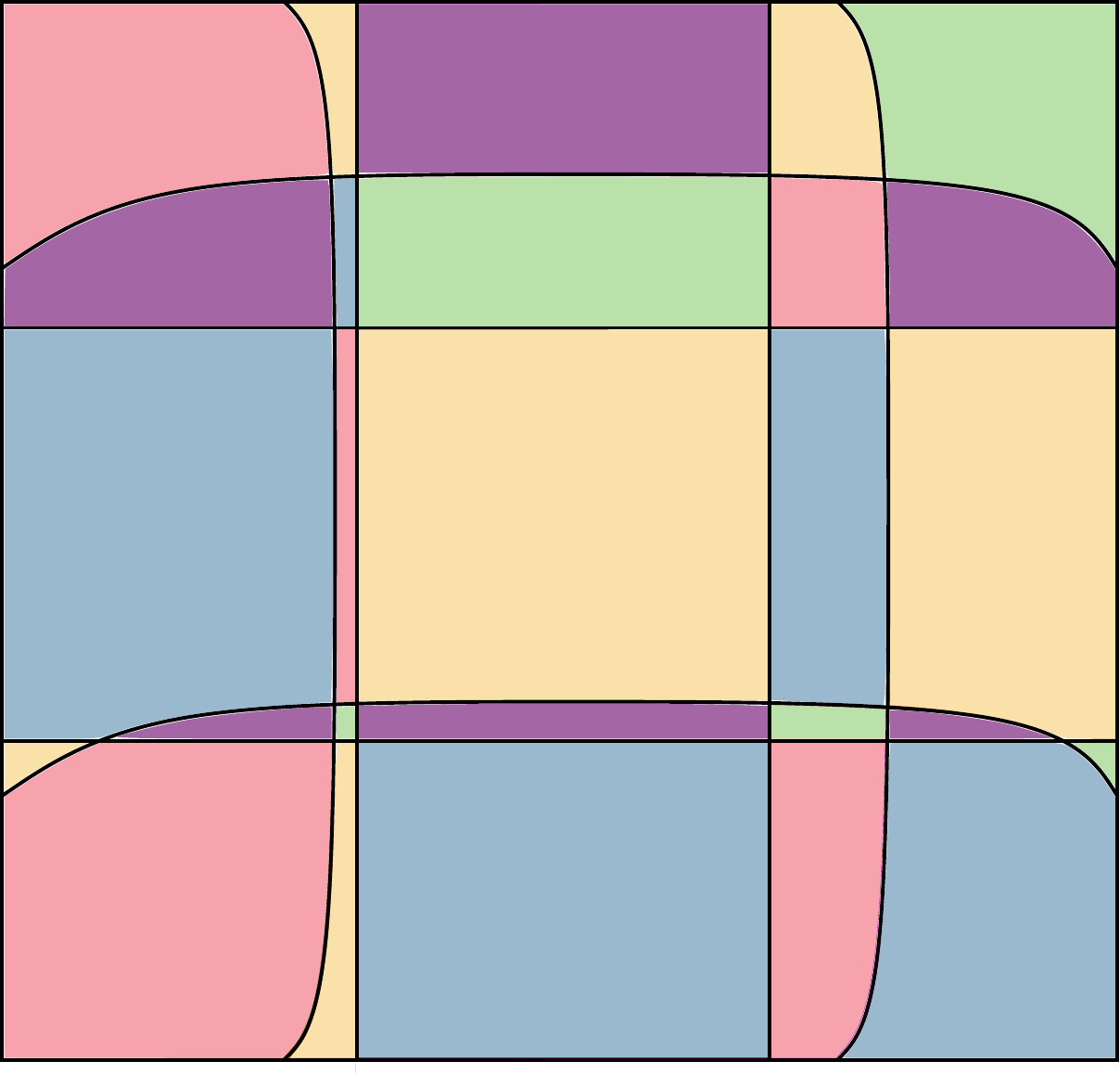}
    			\begin{picture}(0,0)
    			\put(-90,75){$\hat\Gamma_{1,2}$}
    			\end{picture}
   			}
  		}
  		\caption{Mesh intersection between the mesh inherited by $\Gamma_{1,2}$ and the one inherited by $\Omega_2$. In (a), the knots
  						 isoparametric surfaces of $T_2$ are extracted. Their intersection with $\Gamma_{1,2}$ produce the curves represented in
  						 white in (b). Finally, (c) shows the pull-back of the intersection curves in (b) via $T_1^{-1}$, together with the knot 
  						 lines of $\Gamma_{1,2}$ in the parametric space $\hat\Gamma_{1,2}$ and the regions extracted by
  						 Algorithm~\ref{alg:region-extraction} in different colors. See also Figure~\ref{fig:maps}.}
  		\label{fig:pull-back}
		\end{figure*}

		In this section we are going to compute integrals over the interface $\Gamma_{1,2}$. We propose two different examples, one with a
		smooth function and one with a product of splines belonging to the spline spaces of $T_1$ and $T_2$, respectively.
		In order to numerically compute the integral in the latter case it is necessary to find the mesh intersection between the mesh
		inherited by the two splines. By construction, see Equation~\eqref{eq:interface}, $\Gamma_{1,2}$ is always part of the boundary of
		$\Omega^{*}_1$ and therefore it inherits from $T_1$ its mesh information. Our goal is therefore to find out how the mesh of $T_2$
		intersects with the natural one of $\Gamma_{1,2}$. To perform this operation we follow three steps.
		\begin{enumerate}
			\item Extract the knots isoparametric surfaces of $T_2$ in Euclidean space, see Figure~\ref{fig:pull-back}~(a);
						\label{enum:surface-extraction}
			\item Find the intersection curves of each isoparametric surface with the interface $\Gamma_{1,2}$. Note that not all the
						isoparametric surfaces in Step~\ref{enum:surface-extraction} have necessarily an intersection with $\Gamma_{1,2}$, see
						Figure~\ref{fig:pull-back}~(b);
			\item Pull-back the obtained curves in $\hat\Gamma_{1,2}$ via $T_1^{-1}$, see Figure~\ref{fig:pull-back}~(c).
		\end{enumerate}
		The list of pull-back curves, the parametric grid of $\hat\Gamma_{1,2}$, and the boundary curve of $\hat\Gamma_{1,2}$, together with
		their intersection points compose now a curvilinear drawing as in Definition~\ref{def:curve-drawing} and we can hence apply
		Algorithm~\ref{alg:region-extraction} in order to extract the list of regions $\Psi$ that represents the mesh intersection of the
		interface $\Gamma_{1,2}$, see colored regions in Figure~\ref{fig:pull-back}~(c). We remark that several of the vertices of
		the curvilinear drawing are represented by the intersections of curves with straight lines that are parallel to the $u$ and $v$
		directions and therefore can be retrieved very efficiently.

\begin{figure*}[t]\centering
    	\subfloat[] {    	
    		\def\svgwidth{0.8\columnwidth}
    		%!TEX root = ../ANTOLIN2021REI.tex
%% Creator: Inkscape 1.1 (c68e22c387, 2021-05-23), www.inkscape.org
%% Accompanies image file 'numerical_integration_models.pdf' (pdf, eps, ps)
\begingroup%
  \makeatletter%
  \providecommand\color[2][]{%
    \errmessage{(Inkscape) Color is used for the text in Inkscape, but the package 'color.sty' is not loaded}%
    \renewcommand\color[2][]{}%
  }%
  \providecommand\transparent[1]{%
    \errmessage{(Inkscape) Transparency is used (non-zero) for the text in Inkscape, but the package 'transparent.sty' is not loaded}%
    \renewcommand\transparent[1]{}%
  }%
  \providecommand\rotatebox[2]{#2}%
  \newcommand*\fsize{\dimexpr\f@size pt\relax}%
  \newcommand*\lineheight[1]{\fontsize{\fsize}{#1\fsize}\selectfont}%
  \ifx\svgwidth\undefined%
    \setlength{\unitlength}{5661bp}%
    \ifx\svgscale\undefined%
      \relax%
    \else%
      \setlength{\unitlength}{\unitlength * \real{\svgscale}}%
    \fi%
  \else%
    \setlength{\unitlength}{\svgwidth}%
  \fi%
  \global\let\svgwidth\undefined%
  \global\let\svgscale\undefined%
  \makeatother%
  \begin{picture}(1,0.75382535)%
    \lineheight{1}%
    \setlength\tabcolsep{0pt}%
    \put(0,0){\includegraphics[width=\unitlength,page=1]{img/numerical_integration_models.pdf}}%
    \put(0.14138577,0.72){\color[rgb]{0,0,0}\makebox(0,0)[lt]{\lineheight{1.25}\smash{\begin{tabular}[t]{l}$\Omega_1^*$\end{tabular}}}}%
    \put(0.75,0.04){\color[rgb]{0,0,0}\makebox(0,0)[lt]{\lineheight{1.25}\smash{\begin{tabular}[t]{l}$\Omega_2^*$\end{tabular}}}}%
  \end{picture}%
\endgroup%

    	}
    	\qquad
    	\subfloat[] {    	
    		\def\svgwidth{0.57\columnwidth}
    		%!TEX root = ../ANTOLIN2021REI.tex
%% Creator: Inkscape 1.1 (c68e22c387, 2021-05-23), www.inkscape.org
%% Accompanies image file 'numerical_integration.pdf' (pdf, eps, ps)
\begingroup%
  \makeatletter%
  \providecommand\color[2][]{%
    \errmessage{(Inkscape) Color is used for the text in Inkscape, but the package 'color.sty' is not loaded}%
    \renewcommand\color[2][]{}%
  }%
  \providecommand\transparent[1]{%
    \errmessage{(Inkscape) Transparency is used (non-zero) for the text in Inkscape, but the package 'transparent.sty' is not loaded}%
    \renewcommand\transparent[1]{}%
  }%
  \providecommand\rotatebox[2]{#2}%
  \newcommand*\fsize{\dimexpr\f@size pt\relax}%
  \newcommand*\lineheight[1]{\fontsize{\fsize}{#1\fsize}\selectfont}%
  \ifx\svgwidth\undefined%
    \setlength{\unitlength}{467.03762704bp}%
    \ifx\svgscale\undefined%
      \relax%
    \else%
      \setlength{\unitlength}{\unitlength * \real{\svgscale}}%
    \fi%
  \else%
    \setlength{\unitlength}{\svgwidth}%
  \fi%
  \global\let\svgwidth\undefined%
  \global\let\svgscale\undefined%
  \makeatother%
  \begin{picture}(1,0.95098331)%
    \lineheight{1}%
    \setlength\tabcolsep{0pt}%
    \put(0,0){\includegraphics[width=\unitlength,page=1]{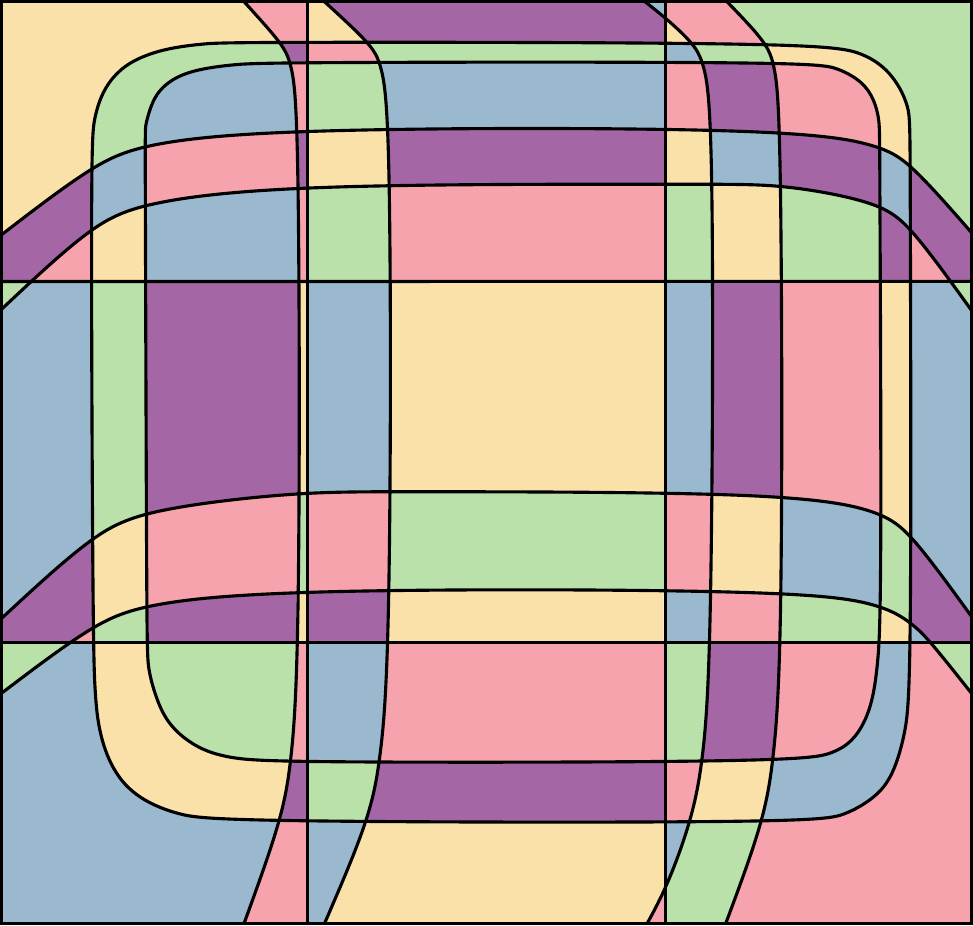}}%
    \put(0.8,0.99){\color[rgb]{0,0,0}\makebox(0,0)[lt]{\lineheight{1.25}\smash{\begin{tabular}[t]{l}$\hat\Gamma_{1,2}$\end{tabular}}}}%
  \end{picture}%
\endgroup%

    	}
  		\caption{Example described in Section~\ref{sub:computation_of_integrals_over_planar_domains}. A volumetric box $\Omega^*_1$ and a
  		volumetric sphere $\Omega^*_2$ with non empty intersection are shown in (a). The intersection curves between the knot surfaces of
  		$T_2$ and $\Gamma_{1,2}$ are pulled-back in the parametric space $\hat\Gamma_{1,2}$ and, together with the knot lines of the
  		interface, form the curvilinear drawing in~(b). The regions to be extracted are marked with different colors.}
  		\label{fig:example-0}
		\end{figure*}

		In both the examples of this section we consider the geometric setting in Figure~\ref{fig:example-0}~(a).
		The solids $\Omega^*_1$ and $\Omega_2^*$ are parameterized by two trivariate B-splines $T_1\in\bbs_{\dd_1,\tl_1}$ and
		$T_2\in\bbs_{\dd_2,\tl_2}$, respectively. Figure~\ref{fig:example-0}~(b) shows the pull-back of the intersection curves between the 
		isoparametric surfaces of $\Omega^*_2$ and $\Gamma_{1,2}$ in the parametric space of $T_1$.

		In the first example we compute the integral
		\begin{equation}
			\label{eq:num_integration_1}
			\int_{\Gamma_{1,2}}\sin\Big(\frac{\pi}{2}x\Big)\cos(\pi y)e^x,
		\end{equation}
		using our algorithm and comparing the result with the one obtained using standard quadrature rules.
		Denoting with
		\[
			f(x,y) = \sin\Big(\frac{\pi}{2}x\Big)\cos(\pi y)e^x,
		\]
		Equation~\eqref{eq:num_integration_1} can be rewritten as
		\begin{equation}
			\label{eq:num_integration_2}
			\int_{\Gamma_{1,2}}f = \int_{\hat\Gamma_{1,2}}f\circ S_{1,2} \det(\nabla S_{1,2}),
		\end{equation}
		where $S_{1,2}\colon \hat\Gamma_{1,2}\rightarrow\Gamma_{1,2}$ is a parameterization of the interface $\Gamma_{1,2}$, that is
		$S_{1,2} = \restr{T_1}{\hat\Gamma_{1,2}}$.
		Since $f$ is an analytic function, Equation~\eqref{eq:num_integration_2} can be numerically computed easily with standard quadrature
		rules.
		We denote with $I_f$ the value of~\eqref{eq:num_integration_1} obtained in such a way, using an overkill number of quadrature points.
		In order to test our algorithm, we further write~\eqref{eq:num_integration_2} as
		\begin{equation}
			\label{eq:num_integration_3}
			\begin{aligned}
				\int_{\hat\Gamma_{1,2}}f\circ & S_{1,2} \det(\nabla S_{1,2}) = \\
					&\sum_{\psi\in\Psi}\int_{\psi}f\circ S_{1,2} \det(\nabla S_{1,2}),
			\end{aligned}
		\end{equation}
		where $\Psi$ contains all the regions extracted by Algorithm~\ref{alg:region-extraction}, see Figure~\ref{fig:example-0}~(b). 
		Finally, in order to compute numerically the integrals in~\eqref{eq:num_integration_3}, we need to create a suitable quadrature
		rule for each region $\psi$. To this end, we apply the untrimming algorithm proposed in~\cite{Wei:2021:IBC} in order to split $\psi$ in
		a list of four-sided, non overlapping, free-form quadrilaterals parameterized as planar parametric patches
		$\pi_1, \dots,\pi_{L_\psi}$
		\[
			\pi_j\colon [0,1]^2\rightarrow \psi, \qquad j = 1,\dots,L_\psi.
		\]
		Each $\pi_j$ is guaranteed to be a B\'ezier patch and the union of their images is a partition of $\psi$ in $\hat\Gamma_{1,2}$.
		Therefore, we can finally compute~\eqref{eq:num_integration_1} as
		\begin{equation}
			\label{eq:num_integration_4}
			\begin{aligned}
				\int_{\hat\Gamma_{1,2}}f\circ & S_{1,2} \det(\nabla S_{1,2}) = \\
					&\sum_{\psi\in\Psi}\sum_{j=1}^{L_\psi}\int_{[0,1]^2}f\circ S_{1,2}\circ\pi_j \det(\nabla \pi_j) \det(\nabla S_{1,2}),	
			\end{aligned}
		\end{equation}
		which can be computed with standard quadrature techniques. We remark that the algorithm proposed in~\cite{Wei:2021:IBC} minimizes the
		number of quadrilaterals necessary to partition each region $\psi$ with a greedy algorithm and therefore $\psi$ is rarely split in
		more than two patches.

		The integrals in~\eqref{eq:num_integration_4} are computed separately using $2^j$ quadrature points per tile
		direction and we denote with $I_f^{(j)}$ the obtained value of the integral in~\eqref{eq:num_integration_1}. These results
		are compared with $I_f$ by considering the error
		\begin{equation}\label{eq:ej}
			E_f^{(j)} = \abs{I_f - I^{(j)}_f}.
		\end{equation}

		In the second test we want to show that Algorithm~\ref{alg:region-extraction} is a suitable method for computing the integral of
		piecewise polynomials defined over $\hat\Gamma_{1,2}$. To this end, we consider the function
		\[
			s(u,v) = s_1(u,v)\tilde s_2(u,v),
		\] 
		where $s_1\in\bbs_{\dd_1,\tl_1}$ and
		\begin{equation}
			\label{eq:push-pull}
			\tilde s_2 = s_2 \circ T_2^{-1} \circ T_1,
		\end{equation}
		for some $s_2\in\bbs_{\dd_2,\tl_2}$.

		As remarked in Section~\ref{sec:applications}, in order to integrate $s$ over $\hat\Gamma_{1,2}$ it is necessary to identify the regions
		in which $s_1$ and $\tilde s_2$ have maximum order of continuity and the algorithm presented in this work allows us to easily recognize
		such regions.
		We therefore express the integral of $s$ as
			\begin{equation}
				\label{eq:int_s}
				\int_{\hat\Gamma_{1,2}} s = \sum_{\psi\in\Psi}\sum_{j=1}^{L_\psi}\int_{[0,1]^2}(s_1\tilde s_2)\circ \pi_j \det(\nabla \pi_j)
			\end{equation}
		and we numerically compute each integral in~\eqref{eq:int_s} using an overkill number of quadrature points for each tile.
		Denoting such value with $I_s$, we define
		\begin{equation}
			\label{eq:ej-spline}
			E^{(j)}_s = \abs{I_s - I_s^{(j)}}
		\end{equation}
		where $I_s^{(j)}$ is the value of the integral approximated using $2^j$ quadrature points per tile direction.

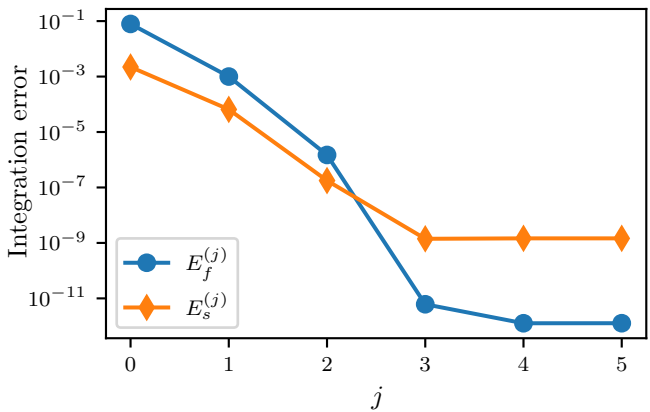
\begin{figure}\centering
    	\includegraphics{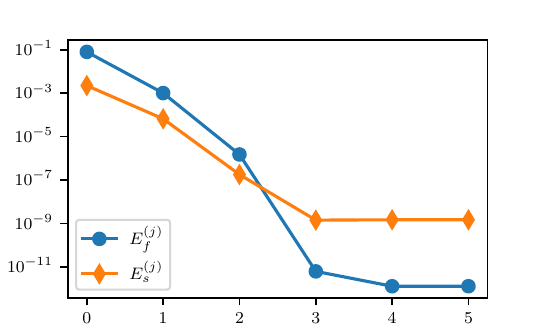}
    	\begin{picture}(0,0)
    		\put(6,5){$j$}
    		\put(-130,60){\rotatebox{90}{Integration error}}
  		\end{picture}%
  		\vspace{-0.2cm}
  		\caption{Semi-log plot of $E_f^{(j)}$ and $E_s^{(j)}$in~\eqref{eq:ej} and~\eqref{eq:ej-spline}, respectively, for $j = 0,\dots,5$. 			See also Table~\ref{tab:numerical_results}.}
  		\label{fig:plot_convergence}
		\end{figure}

		In both experiments, we stop the computation of the quadrature rule if two subsequent approximated values of the integrals are close
		enough, that is when
		\[
			\abs{I^{(j)}_\alpha - I^{(j-1)}_\alpha} < 10^{-12}, \qquad \alpha = f, s.
		\]
		The results of the integrations for both numerical tests are visible in Table~\ref{tab:numerical_results} and
		Figure~\ref{fig:plot_convergence}.
	
		\begin{table*}\centering\footnotesize\onehalfspacing
			\begin{tabular}{c|c|c|c|c|c|c}
				\hline
				$j$ & $0$ & $1$ & $2$ & $3$ & $4$ & $5$\\
				\hline
				\rule{0pt}{3ex}
				$E_f^{(j)}$ & $7.90\eem02$ & $1.00\eem03$ & $1.49\eem06$ & $6.13\eem12$ & $1.27\eem12$ & $1.27\eem12$\\
				\hline
				\rule{0pt}{3ex}
				$E_s^{(j)}$ & $2.21\eem03$ & $6.53\eem05$ & $1.77\eem07$ & $1.41\eem09$ & $1.46\eem09$ & $1.46\eem09$\\
				\hline
  		\end{tabular}
  		\caption{Error values for the integrations in Section~\ref{sub:computation_of_integrals_over_planar_domains}. The table reports the
  		values of $E^{(j)}_f$ and $E^{(j)}_s$ in~\eqref{eq:ej} and~\eqref{eq:ej-spline}, respectively, for $j = 0,\dots,5$.}
  		\label{tab:numerical_results}
		\end{table*}
		We remark that in the computation of $I^{(j)}_s$ geometric operations such as surface-surface intersections, pull-back of curves and
		curve-curve intersections are of utmost importance.
		If a pull-back curve is computed coarsely, the image of a quadrature point through $T_2^{-1}\circ T_1$ is not guaranteed to be in the
		right knot element of $s_2$.
		If this happen, the quadrature rule of $s$ is not computed in a region corresponding to its maximum order of continuity and the standard
		order of convergence of numerical integration cannot be guaranteed anymore.

		In our implementation of the surface-surface intersection, we utilize the algorithms provided by Open CASCADE
		Technology~\cite{OCC:2021:OCS}, an open source C{}\verb!++! library that allows a minimal tolerance of $10^{-7}$.
		Other operations such as the pull-back of curves and the curve-curve intersections are instead performed using Irit geometric
		modeler~\cite{Elber:2010:I1U}, which allows to modify the involved tolerances according to our needs. Nevertheless, the pull-backs are
		obtained by inverting pointwise each curve in Euclidean space in the parametric space of $\Gamma_{1,2}$ and then the obtained points are
		approximated using a least squares approach.
		All these approximation operations pollute the geometric setting in which our algorithm works and are the reason for the
		plateau reached by $E_f^{(j)}$ and $E_s^{(j)}$ visible in Figure~\ref{fig:plot_convergence}. Similar phenomena were previously
		reported in~\cite{Antolin:2019:IAO} and~\cite{Antolin:2021:QFI}.

	\subsection{Weak continuity}
		\label{sub:weak_continuity}		
		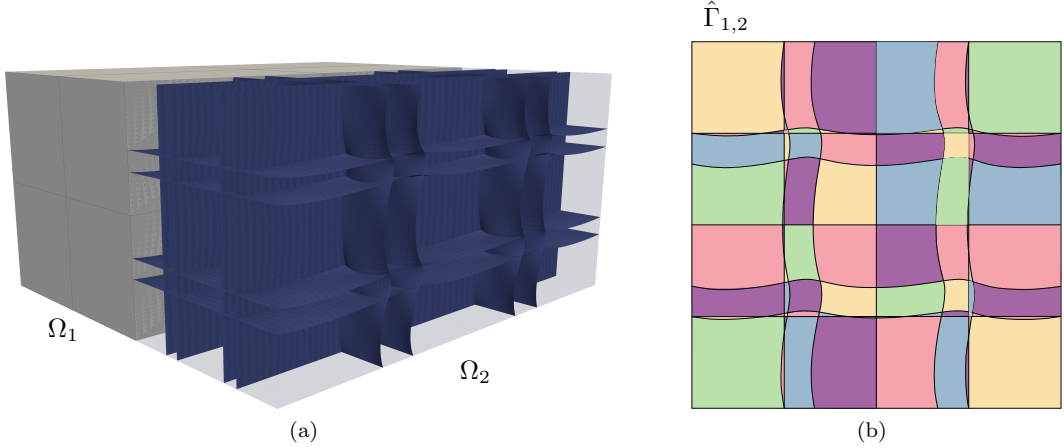
\begin{figure*}[t]\centering
			\subfloat[] {
    		\def\svgwidth{0.9\columnwidth}
    		%!TEX root = ../paper.tex
%% Creator: Inkscape 1.1 (c68e22c387, 2021-05-23), www.inkscape.org
%% Accompanies image file 'weak_continuity_model.pdf' (pdf, eps, ps)
\begingroup%
  \makeatletter%
  \providecommand\color[2][]{%
    \errmessage{(Inkscape) Color is used for the text in Inkscape, but the package 'color.sty' is not loaded}%
    \renewcommand\color[2][]{}%
  }%
  \providecommand\transparent[1]{%
    \errmessage{(Inkscape) Transparency is used (non-zero) for the text in Inkscape, but the package 'transparent.sty' is not loaded}%
    \renewcommand\transparent[1]{}%
  }%
  \providecommand\rotatebox[2]{#2}%
  \newcommand*\fsize{\dimexpr\f@size pt\relax}%
  \newcommand*\lineheight[1]{\fontsize{\fsize}{#1\fsize}\selectfont}%
  \ifx\svgwidth\undefined%
    \setlength{\unitlength}{4455bp}%
    \ifx\svgscale\undefined%
      \relax%
    \else%
      \setlength{\unitlength}{\unitlength * \real{\svgscale}}%
    \fi%
  \else%
    \setlength{\unitlength}{\svgwidth}%
  \fi%
  \global\let\svgwidth\undefined%
  \global\let\svgscale\undefined%
  \makeatother%
  \begin{picture}(1,0.58804714)%
    \lineheight{1}%
    \setlength\tabcolsep{0pt}%
    \put(0,0){\includegraphics[width=\unitlength,page=1]{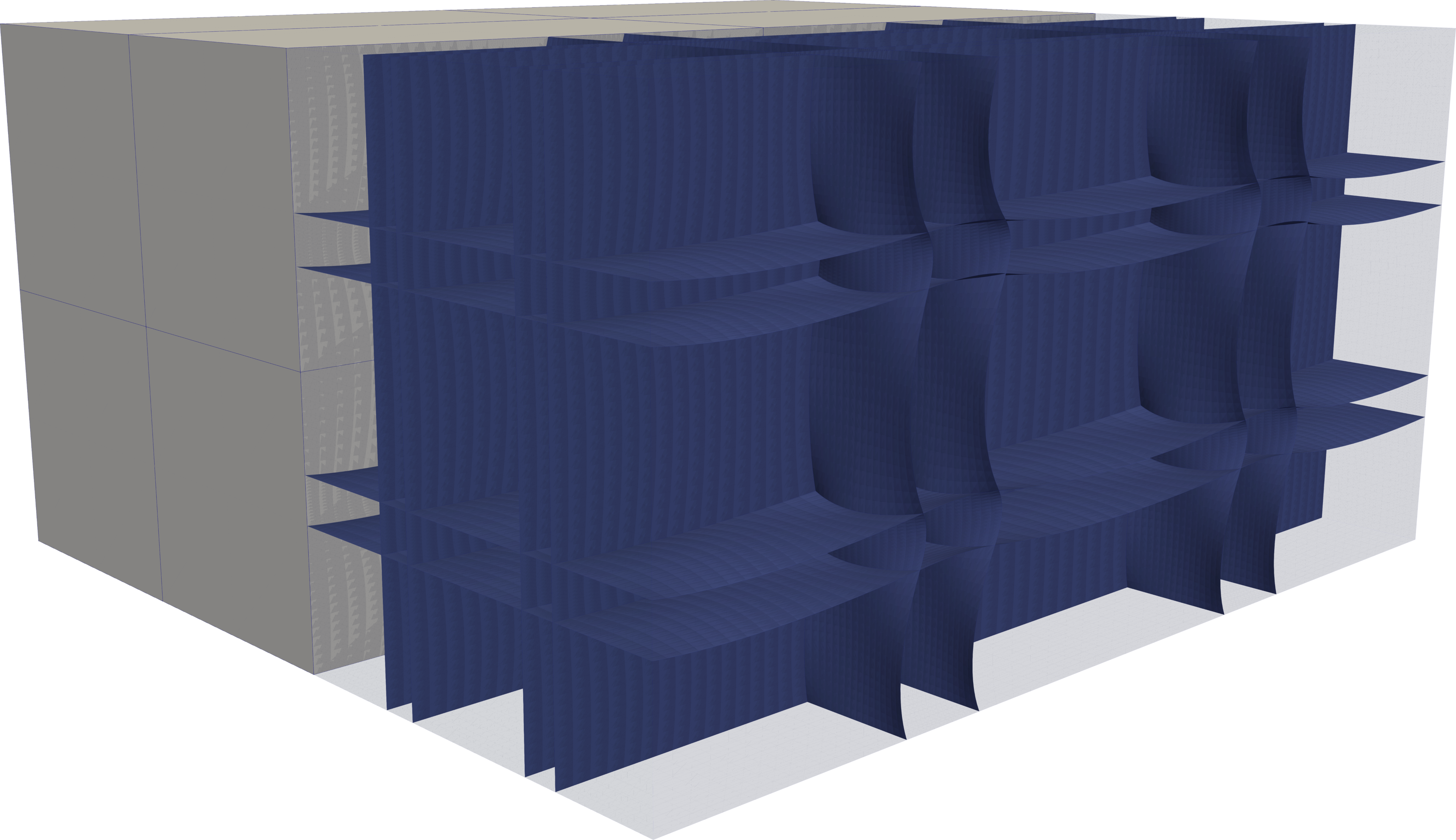}}%
    \put(0.75,0.05){\color[rgb]{0,0,0}\makebox(0,0)[lt]{\lineheight{1.25}\smash{\begin{tabular}[t]{l}$\Omega_2$\end{tabular}}}}%
    \put(0.07,0.12){\color[rgb]{0,0,0}\makebox(0,0)[lt]{\lineheight{1.25}\smash{\begin{tabular}[t]{l}$\Omega_1$\end{tabular}}}}%
  \end{picture}%
\endgroup%

   		}\qquad
   		\subfloat[] {
    		\def\svgwidth{0.55\columnwidth}
    		%!TEX root = ../paper.tex
%% Creator: Inkscape 1.1 (c68e22c387, 2021-05-23), www.inkscape.org
%% Accompanies image file 'drawing.pdf' (pdf, eps, ps)
\begingroup%
  \makeatletter%
  \providecommand\color[2][]{%
    \errmessage{(Inkscape) Color is used for the text in Inkscape, but the package 'color.sty' is not loaded}%
    \renewcommand\color[2][]{}%
  }%
  \providecommand\transparent[1]{%
    \errmessage{(Inkscape) Transparency is used (non-zero) for the text in Inkscape, but the package 'transparent.sty' is not loaded}%
    \renewcommand\transparent[1]{}%
  }%
  \providecommand\rotatebox[2]{#2}%
  \newcommand*\fsize{\dimexpr\f@size pt\relax}%
  \newcommand*\lineheight[1]{\fontsize{\fsize}{#1\fsize}\selectfont}%
  \ifx\svgwidth\undefined%
    \setlength{\unitlength}{566.49452332bp}%
    \ifx\svgscale\undefined%
      \relax%
    \else%
      \setlength{\unitlength}{\unitlength * \real{\svgscale}}%
    \fi%
  \else%
    \setlength{\unitlength}{\svgwidth}%
  \fi%
  \global\let\svgwidth\undefined%
  \global\let\svgscale\undefined%
  \makeatother%
  \begin{picture}(1,1.08434311)%
    \lineheight{1}%
    \setlength\tabcolsep{0pt}%
    \put(0,0){\includegraphics[width=\unitlength,page=1]{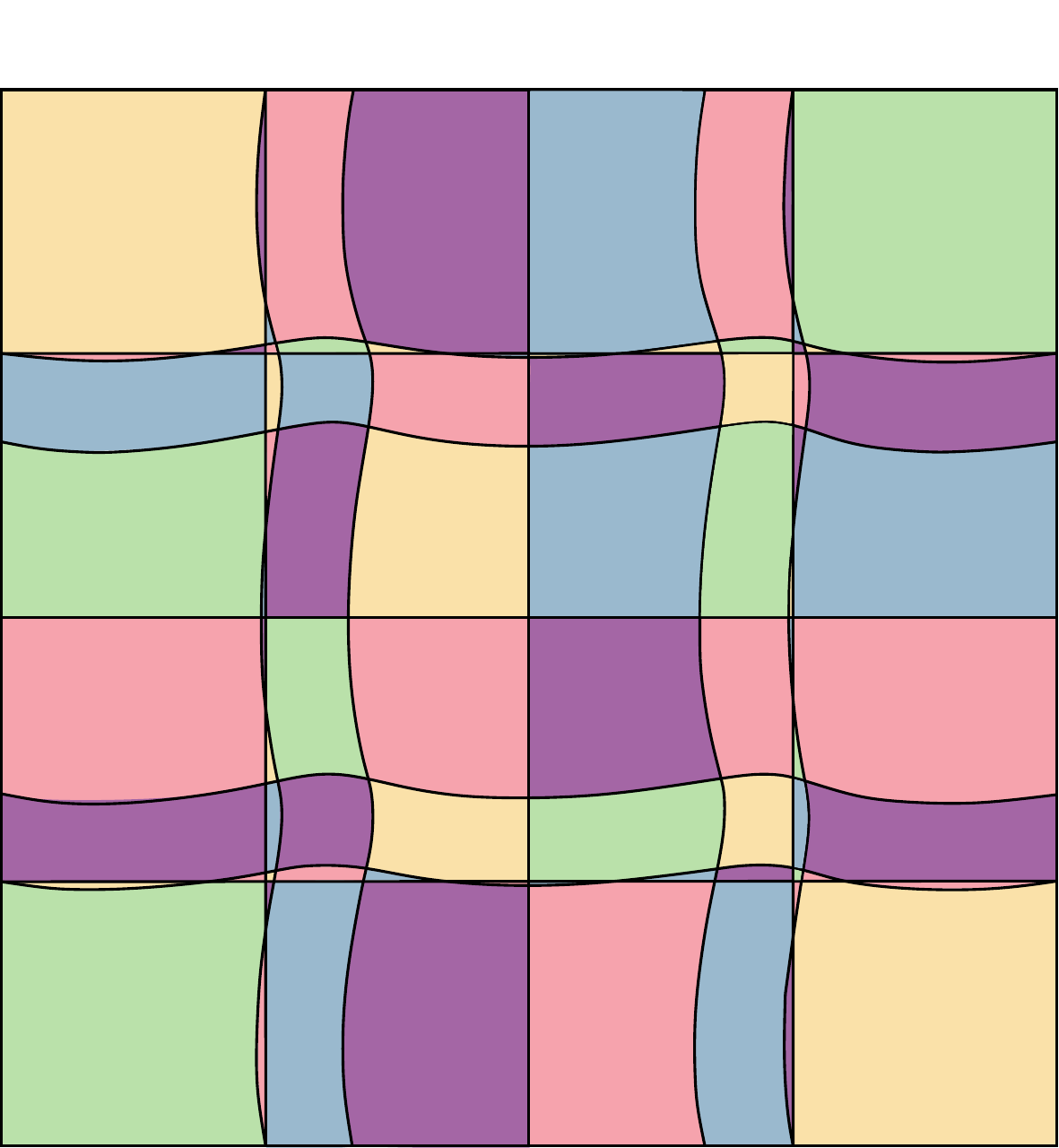}}%
    \put(0.03252826,1.04359097){\color[rgb]{0,0,0}\makebox(0,0)[lt]{\lineheight{1.25}\smash{\begin{tabular}[t]{l}$\hat\Gamma_{1,2}$\end{tabular}}}}%
  \end{picture}%
\endgroup%

   		}
  		\caption{Example described in Section~\ref{sub:weak_continuity}. The original geometries $\Omega_1$ and $\Omega_2$ are shown in~(a).
  						 The pull-back curves of the intersection between the isosurfaces of $\Omega_2$ and the interface $\Gamma_{1,2}$ are shown
  						 in~(b). The regions to be extracted are marked with different colors.}
  		\label{fig:example-wc-res}
		\end{figure*}
		The algorithm proposed in this work can be used to enforce weak continuity constraints to solids geometries. Weak continuity
		constraints have been imposed to bidimensional geometries by Zou and colleagues in~\cite{Zou:2018:IBD}. Their approach is based on
		a newly defined B{\'e}zier projector. A similar result can be achieved with a Lee-Lyche-M{\o}rken 
		quasi-interpolant~\cite{Lee:2000:SEO}, taking advantage of the presented algorithm. We remark that also the method presented
		in~\cite{Zou:2018:IBD} can be implemented using our algorithm for improved accuracy.

		Let us denote with $\Omega^*_1$ and $\Omega^*_2$ two solids parameterized by two trivariate B-splines $T_1$ and $T_2$ and with $\Omega$
		their Boolean union. Let $\hat\Omega^*_2$ be a deformation of $\Omega^*_2$ such that there exists a trivariate \emph{displacement}
		B-spline $\delta T_2$ such that $T_2 + \delta T_2$ is a parameterization of the closure of $\hat\Omega^*_2$. Our goal is to find
		a corresponding trivariate displacement $\delta T_1$ such that $T_2 + \delta T_2$ and $T_1 + \delta T_1$ form a weakly continuous
		piecewise parameterization of
		\[
			\hat\Omega = \hat\Omega^*_1 \cup \hat\Omega^*_2,
		\]
		where $\hat\Omega^*_1$ is the domain parameterized by $T_1 + \delta T_1$, representing the corresponding deformation of $\Omega^*_1$.

		To this end we find $\delta T_1$ as the Lee--Lyche--M{\o}rken quasi-interpolant of $\delta T_2$, that is
		$\delta T_1 = I_\text{LLM}\delta T_2$. This family of quasi-interpolants is a widely used method for the local projection a function
		$f\in\mathcal{L}^2(D)$, for some domain $D$, into a given spline space $\mathbb{S}_{\mathbf{d},\tl}$
		\[
			I_{\text{LLM}}\colon \mathcal{L}^2(D)\rightarrow \mathbb{S}_{\mathbf{d},\tl}
		\]
		In their work~\cite{Lee:2000:SEO}, Lee, \emph{Lyche,} and M{\o}rken proposed a procedure to build such quasi-interpolants that is based on the
		use of local spline projectors, see Algorithm~\ref{alg:LLM}.
		There are two main advantages in using a Lee-Lyche-M{\o}rken approach in this setting. On the one hand, being based on local spline 
		projectors, only the elements of $\Omega_1$ close to the interface are to be influenced and therefore a refinement step can help us at
		controlling the influence that $\delta T_2$ has on $\Omega^*_1$. On the other hand, we are sure of exactly reproducing
		$\delta T_2$, provided that the spline space $\mathbb{S}_{\dd_1,\tl_1}$ is large enough.		

		\begin{algorithm}
  		\caption{LLM quasi-interpolant : Computes a quasi-interpolant for a given function $f$; see~\cite{Lee:2000:SEO}.}
  		\label{alg:LLM}
  		\hspace*{\algorithmicindent} \textbf{Input} Degrees $\dd = (d_u, d_v, d_w)$;\\
  		\hspace*{\algorithmicindent} \textbf{Input} Knots vectors $\tl = (\tl_u, \tl_v, \tl_w)$;\\
  		\hspace*{\algorithmicindent} \textbf{Input} Function $f$;\\
  		\hspace*{\algorithmicindent} \textbf{Output} $I_\text{LLM} f$;
  		\begin{algorithmic}[1]
  			\State $N \gets$ number of degrees of freedom in $\mathbb{S}_{\mathbf{d},\tl}$;
    		\For{$\ell = 0, \dots, N$}
    			\State $K_\ell \gets$ knot interval such that $K_\ell\cap \text{supp} (B_{\ell,\dd, \tl}) \ne\emptyset$;\label{line:K_ell-definition}
    			\State $\mathbb{S}_{\mathbf{d},\tl, K_\ell} \gets$ restriction of $\mathbb{S}_{\mathbf{d},\tl}$ to interval $K_\ell$;\label{line:spline-space-restriction}
    			\State $P_{K_\ell} \gets$ local projector on the space $\mathbb{S}_{\mathbf{d},\tl, K_\ell}$;\label{line:projector-definition}
    			\State $f_{K_\ell} \gets$ restriction of $f$ to interval $K_\ell$;
    			\State $c_\ell \gets $ corresponding coefficient in $P_{K_\ell}f_{K_\ell}$;\label{line:c_ell-def}
    		\EndFor
  			\State \Return $\sum_{\ell = 0}^N c_\ell B_{\ell, \dd, \tl}$;
  		\end{algorithmic}
		\end{algorithm}

		Figure~\ref{fig:example-wc-res} describes the geometric setting of this numerical experiment. In Figure~\ref{fig:example-wc-res}~(a),
		the domains $\Omega^*_1$ and $\Omega^*_2$ are represented together with the isoparametric knot surfaces of $\Omega_2$. In this case
		there is no intersection between the bodies $\Omega_1$ and $\Omega_2$ and therefore $\Omega^*_i = \Omega_i$, $i = 1, 2$, while the
		interface $\Gamma_{1,2}$ is simply the intersection of the boundaries of the solids. The
		intersections of these surfaces with the interface $\Gamma_{1,2}$ are pulled-back in $T_1$'s parametric space and form
		the regions shown in Figure~\ref{fig:example-wc-res}~(b).

		Algorithm~\ref{alg:LLM} is a general procedure and allows us to choose freely both the knot interval $K_\ell$ in
		Line~\ref{line:K_ell-definition} and the local projector $P_{K_\ell}$ in Line~\ref{line:projector-definition}.
		Since the only condition about $K_\ell$ is to have a non-empty intersection with the support of $B_{\ell,\dd_1}$, in our work we simply
		set
		\[
			K_\ell = \text{supp} (B_{\ell,\dd_1}),
		\]
		while we use as the local projector $P_{K_\ell}$ in Line~\ref{line:projector-definition} the usual $\mathcal{L}^2$ projector that
		satisfies
		\begin{equation}\label{eq:projector}
			\int_{K_\ell} P_{K_\ell}[f] B_{l, \dd_1} = \int_{K_\ell} f B_{l, \dd_1}, \qquad l\in\Lambda_\ell,
		\end{equation}
		for any function $f$, where
		\[
			 \Lambda_\ell = \{l \colon K_\ell\cap\text{supp}(B_{l,\dd_1})\ne\emptyset\}.
		\]
		We remark that, thanks to Equation~\eqref{eq:data-structure}, the knot interval $K_\ell$ in Algorithm~\ref{alg:LLM},
		Line~\ref{line:K_ell-definition}, is a knot interval of the boundary surface containing the local interface $\Gamma_{1,2}$ and
		therefore both integrals in Equation~\eqref{eq:projector} are integrals of trivariate functions over a planar domain.

		Let us denote with $L$ the cardinality of $\Lambda_\ell$ and with $\hat\imath$ the index of the basis function in $\Lambda_\ell$
		corresponding to $i$, for any $i$ such that $1\le i\le L$. By writing explicitly
		\[
			P_{K_\ell}[\delta T_2] = \sum_{i = 1}^L \lambda_i B_{\hat\imath, \dd_1},
		\]
		Equation~\eqref{eq:projector} reads
		\begin{equation}\label{eq:projector-explicit}
			\sum_{i = 1}^L \lambda_i \int_{K_\ell} B_{\hat\imath, \dd_1}B_{l, \dd_1} = \int_{K_\ell} \delta T_2 B_{l, \dd_1}, \qquad l\in\Lambda_\ell,	
		\end{equation}
		for some unknown coefficients $\lambda_i$. In order to keep the notation as simple as possible, here and in the rest of this section we
		denote with $\delta T_2$  both the function $\delta T_2\in\mathbb{S}_{\dd_2,\tl_2}$ and the corresponding function as defined
		in~\eqref{eq:push-pull}. While this is formally incorrect, we think that no confusion is likely to arise as all integrations are
		performed in the parametric space of $T_1$.

		In order to determine the $\lambda_i$s it is necessary to solve the (sparse) linear system
		\[
			G_\ell\lambda = P,
		\]
		where $\lambda = (\lambda_1,\dots,\lambda_L)^T$,
		\begin{equation}
			\label{eq:right-hand-side}
			P = \bigg(\int_{K_\ell} \delta T_2 B_{\hat 1, \dd_1},\dots,\int_{K_\ell} \delta T_2 B_{\hat L, \dd_1}\bigg)^T
		\end{equation}
		and $G_\ell$ is the Gibbs matrix with entries
		\[
			(G_\ell)_{i, j} = \int_{K_\ell} B_{\hat \imath, \dd_1}B_{\hat \jmath, \dd_1}, \qquad i,j = 1,\dots, L.
		\]
		While the entries of the Gibbs matrix $G_\ell$ are trivial to compute numerically, as $K_\ell$ is an element for the spline
		space in which the basis functions $B_{\hat \imath, \dd_1}$ are defined, the components of the vector $P$ require a special
		treatment. $\delta T_2$ indeed does not belong to the same spline space of the B-spline basis functions $B_{\hat 1, \dd_1}, \dots,
		B_{\hat L, \dd_1}$. We are therefore left with the key problem of the creation of the quadrature mesh along the preimage of the
		interface $\Gamma_{1,2}$. In order to apply the same strategy as in Section~\ref{sub:computation_of_integrals_over_planar_domains}, the
		isoparametric surfaces of $T_2$ are extended in the three parametric directions before computing their intersection with the interface
		$\Gamma_{1,2}$. The regions representing the mesh intersection of $T_1$ and $T_2$ are shown in Figure~\ref{fig:example-wc-res}~(b).
		Therefore, we proceed to write the $i$-th entry of $P$ as
		\begin{equation}
			\label{eq:Pi-weak-continuity}
			P_i = \int_{K_\ell} \delta T_2 B_{\hat \imath, \dd_1} = \sum_{\psi\in\Psi_\ell} \int_{K_\ell} \delta T_2 B_{\hat \imath, \dd_1},
		\end{equation}
		where $\Psi_\ell = \{\psi\in\Psi\colon \psi\subset K_\ell\}$. In each region $\psi\in\Psi_\ell$, the restriction of
		$B_{\hat \imath,\dd_1}$ is a polynomial of degree $\dd_1$, since $\psi$ is guaranteed to be contained in a single knot span of
		$\bbs_{\tl_1,\dd_1}$, while $\delta T_2$ is a polynomial of degree $\dd_2$ on $T_2^{-1}\circ T_1(\psi)$, as the latter is fully
		contained in a knot span of $\bbs_{\tl_2,\dd_2}$.
	
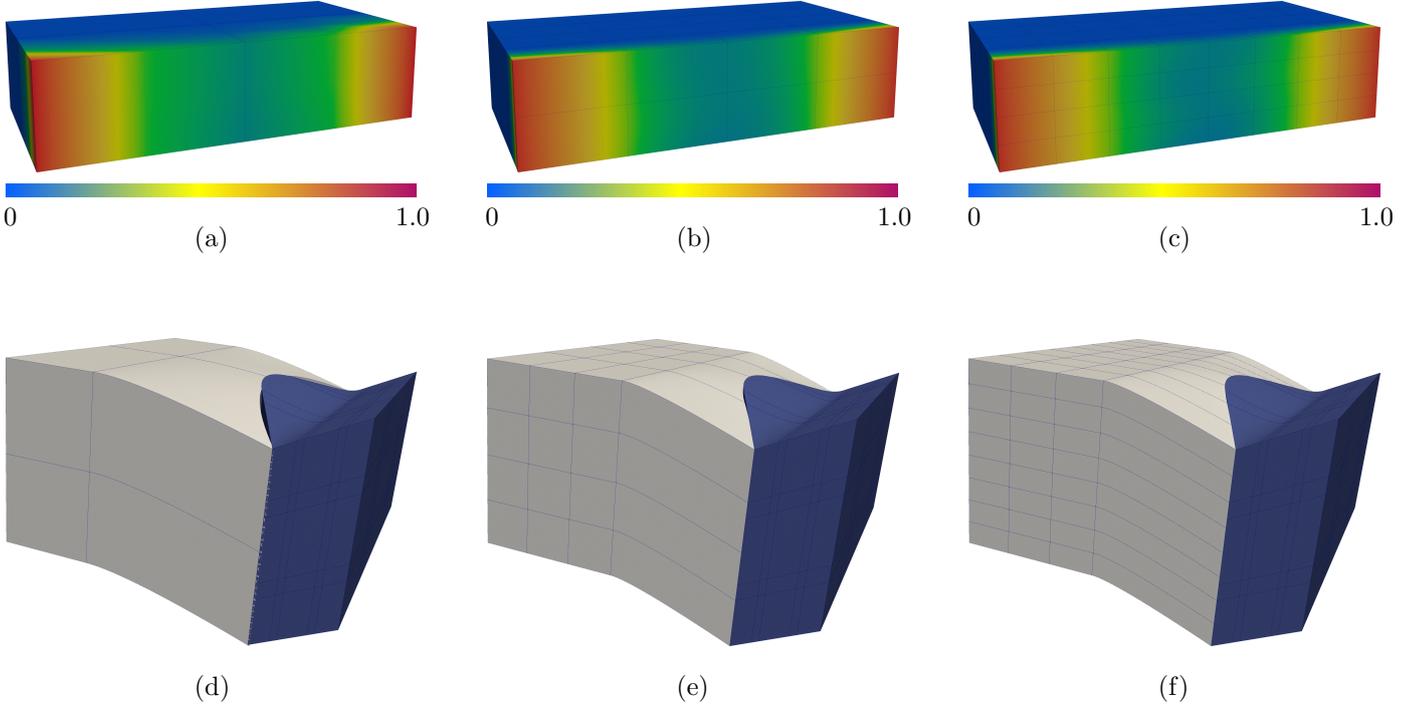
\begin{figure*}[t]\centering
			\mbox {
    		\def\svgwidth{.61\columnwidth}
    		%!TEX root = ../paper.tex
%% Creator: Inkscape 1.1 (c68e22c387, 2021-05-23), www.inkscape.org
%% Accompanies image file 'weak_conf0.pdf' (pdf, eps, ps)
\begingroup%
  \makeatletter%
  \providecommand\color[2][]{%
    \errmessage{(Inkscape) Color is used for the text in Inkscape, but the package 'color.sty' is not loaded}%
    \renewcommand\color[2][]{}%
  }%
  \providecommand\transparent[1]{%
    \errmessage{(Inkscape) Transparency is used (non-zero) for the text in Inkscape, but the package 'transparent.sty' is not loaded}%
    \renewcommand\transparent[1]{}%
  }%
  \providecommand\rotatebox[2]{#2}%
  \newcommand*\fsize{\dimexpr\f@size pt\relax}%
  \newcommand*\lineheight[1]{\fontsize{\fsize}{#1\fsize}\selectfont}%
  \ifx\svgwidth\undefined%
    \setlength{\unitlength}{504.07766799bp}%
    \ifx\svgscale\undefined%
      \relax%
    \else%
      \setlength{\unitlength}{\unitlength * \real{\svgscale}}%
    \fi%
  \else%
    \setlength{\unitlength}{\svgwidth}%
  \fi%
  \global\let\svgwidth\undefined%
  \global\let\svgscale\undefined%
  \makeatother%
  \begin{picture}(1,0.96772816)%
    \lineheight{1}%
    \setlength\tabcolsep{0pt}%
    \put(0,0){\includegraphics[width=\unitlength,page=1]{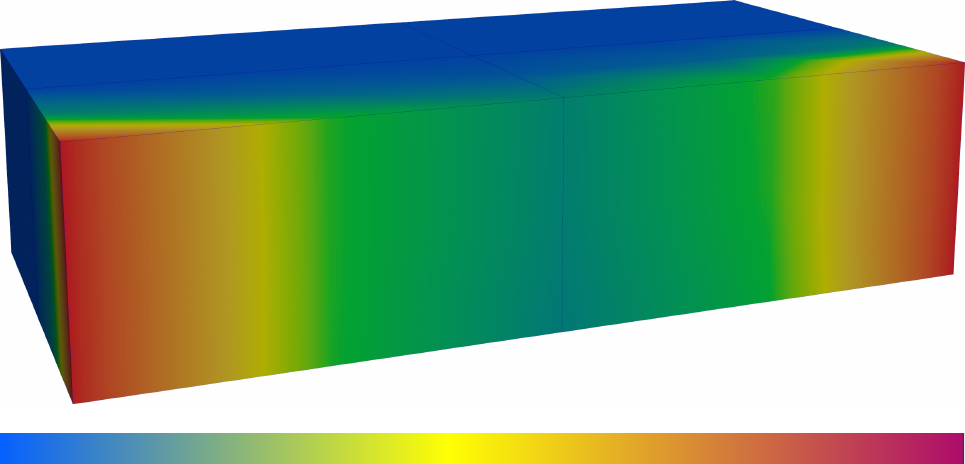}}%
    \put(-0.00450428,-0.07){\color[rgb]{0,0,0}\makebox(0,0)[lt]{\lineheight{1.25}\smash{\begin{tabular}[t]{l}$0$\end{tabular}}}}%
    \put(0.95,-0.07){\color[rgb]{0,0,0}\makebox(0,0)[lt]{\lineheight{1.25}\smash{\begin{tabular}[t]{l}$1.0$\end{tabular}}}}%
    \put(0.46,-0.12){\color[rgb]{0,0,0}\makebox(0,0)[lt]{\smash{\begin{tabular}[t]{l}(a)\end{tabular}}}}%
  \end{picture}%
\endgroup%

    		\qquad
    		\def\svgwidth{.61\columnwidth}
    		%!TEX root = ../paper.tex
%% Creator: Inkscape 1.1 (c68e22c387, 2021-05-23), www.inkscape.org
%% Accompanies image file 'weak_conf1.pdf' (pdf, eps, ps)
\begingroup%
  \makeatletter%
  \providecommand\color[2][]{%
    \errmessage{(Inkscape) Color is used for the text in Inkscape, but the package 'color.sty' is not loaded}%
    \renewcommand\color[2][]{}%
  }%
  \providecommand\transparent[1]{%
    \errmessage{(Inkscape) Transparency is used (non-zero) for the text in Inkscape, but the package 'transparent.sty' is not loaded}%
    \renewcommand\transparent[1]{}%
  }%
  \providecommand\rotatebox[2]{#2}%
  \newcommand*\fsize{\dimexpr\f@size pt\relax}%
  \newcommand*\lineheight[1]{\fontsize{\fsize}{#1\fsize}\selectfont}%
  \ifx\svgwidth\undefined%
    \setlength{\unitlength}{463.13176332bp}%
    \ifx\svgscale\undefined%
      \relax%
    \else%
      \setlength{\unitlength}{\unitlength * \real{\svgscale}}%
    \fi%
  \else%
    \setlength{\unitlength}{\svgwidth}%
  \fi%
  \global\let\svgwidth\undefined%
  \global\let\svgscale\undefined%
  \makeatother%
  \begin{picture}(1,0.97917571)%
    \lineheight{1}%
    \setlength\tabcolsep{0pt}%
    \put(0,0){\includegraphics[width=\unitlength,page=1]{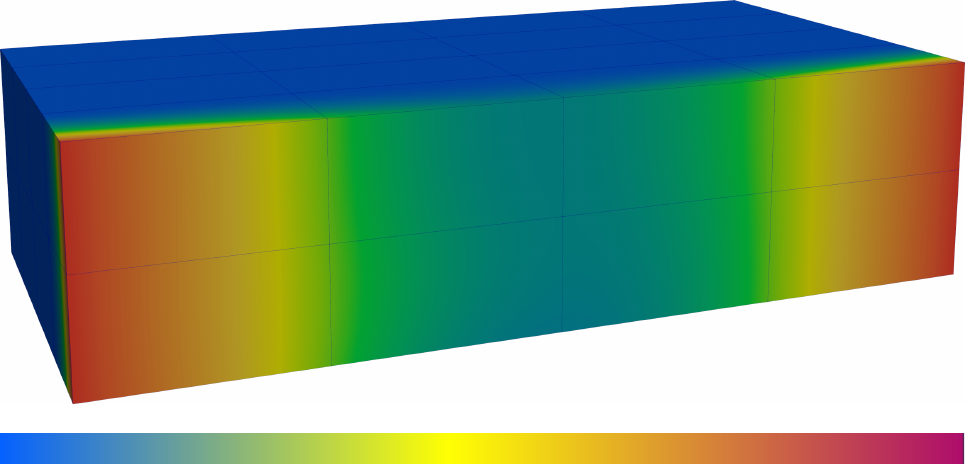}}%
    \put(-0.00450428,-0.07){\color[rgb]{0,0,0}\makebox(0,0)[lt]{\lineheight{1.25}\smash{\begin{tabular}[t]{l}$0$\end{tabular}}}}%
    \put(0.95,-0.07){\color[rgb]{0,0,0}\makebox(0,0)[lt]{\lineheight{1.25}\smash{\begin{tabular}[t]{l}$1.0$\end{tabular}}}}%
    \put(0.46,-0.12){\color[rgb]{0,0,0}\makebox(0,0)[lt]{\smash{\begin{tabular}[t]{l}(b)\end{tabular}}}}%
  \end{picture}%
\endgroup%

    		\qquad
    		\def\svgwidth{.61\columnwidth}
    		%!TEX root = ../paper.tex
%% Creator: Inkscape 1.1 (c68e22c387, 2021-05-23), www.inkscape.org
%% Accompanies image file 'weak_conf2.pdf' (pdf, eps, ps)
\begingroup%
  \makeatletter%
  \providecommand\color[2][]{%
    \errmessage{(Inkscape) Color is used for the text in Inkscape, but the package 'color.sty' is not loaded}%
    \renewcommand\color[2][]{}%
  }%
  \providecommand\transparent[1]{%
    \errmessage{(Inkscape) Transparency is used (non-zero) for the text in Inkscape, but the package 'transparent.sty' is not loaded}%
    \renewcommand\transparent[1]{}%
  }%
  \providecommand\rotatebox[2]{#2}%
  \newcommand*\fsize{\dimexpr\f@size pt\relax}%
  \newcommand*\lineheight[1]{\fontsize{\fsize}{#1\fsize}\selectfont}%
  \ifx\svgwidth\undefined%
    \setlength{\unitlength}{463.13176332bp}%
    \ifx\svgscale\undefined%
      \relax%
    \else%
      \setlength{\unitlength}{\unitlength * \real{\svgscale}}%
    \fi%
  \else%
    \setlength{\unitlength}{\svgwidth}%
  \fi%
  \global\let\svgwidth\undefined%
  \global\let\svgscale\undefined%
  \makeatother%
  \begin{picture}(1,0.97923389)%
    \lineheight{1}%
    \setlength\tabcolsep{0pt}%
    \put(0,0){\includegraphics[width=\unitlength,page=1]{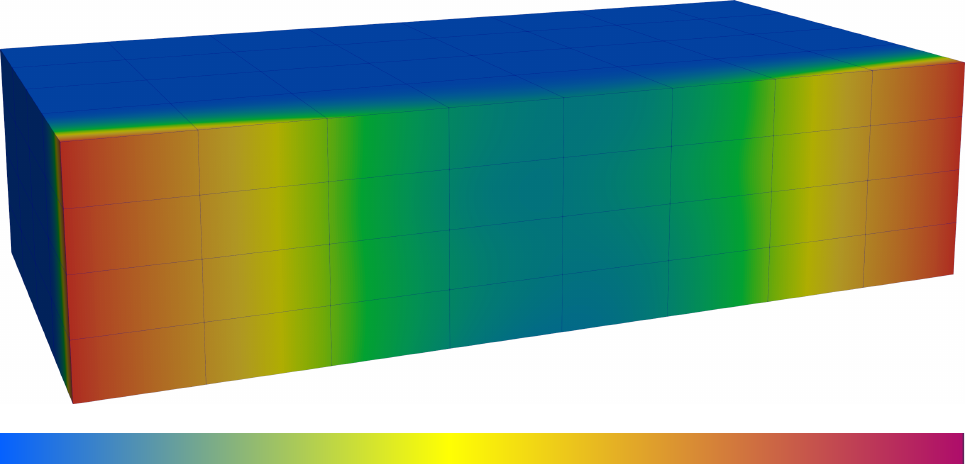}}%
    \put(-0.00450428,-0.07){\color[rgb]{0,0,0}\makebox(0,0)[lt]{\lineheight{1.25}\smash{\begin{tabular}[t]{l}$0$\end{tabular}}}}%
    \put(0.95,-0.07){\color[rgb]{0,0,0}\makebox(0,0)[lt]{\lineheight{1.25}\smash{\begin{tabular}[t]{l}$1.0$\end{tabular}}}}%
    \put(0.46,-0.12){\color[rgb]{0,0,0}\makebox(0,0)[lt]{\smash{\begin{tabular}[t]{l}(c)\end{tabular}}}}%
  \end{picture}%
\endgroup%

   		}\\
   		\vspace{1.5cm}
 			\mbox {
 				\def\svgwidth{.61\columnwidth}
  			%!TEX root = ../paper.tex
%% Creator: Inkscape 1.1 (c68e22c387, 2021-05-23), www.inkscape.org
%% Accompanies image file 'weak_warp0.pdf' (pdf, eps, ps)
\begingroup%
  \makeatletter%
  \providecommand\color[2][]{%
    \errmessage{(Inkscape) Color is used for the text in Inkscape, but the package 'color.sty' is not loaded}%
    \renewcommand\color[2][]{}%
  }%
  \providecommand\transparent[1]{%
    \errmessage{(Inkscape) Transparency is used (non-zero) for the text in Inkscape, but the package 'transparent.sty' is not loaded}%
    \renewcommand\transparent[1]{}%
  }%
  \providecommand\rotatebox[2]{#2}%
  \newcommand*\fsize{\dimexpr\f@size pt\relax}%
  \newcommand*\lineheight[1]{\fontsize{\fsize}{#1\fsize}\selectfont}%
  \ifx\svgwidth\undefined%
    \setlength{\unitlength}{463.13176332bp}%
    \ifx\svgscale\undefined%
      \relax%
    \else%
      \setlength{\unitlength}{\unitlength * \real{\svgscale}}%
    \fi%
  \else%
    \setlength{\unitlength}{\svgwidth}%
  \fi%
  \global\let\svgwidth\undefined%
  \global\let\svgscale\undefined%
  \makeatother%
  \begin{picture}(1,0.8173757)%
    \lineheight{1}%
    \setlength\tabcolsep{0pt}%
    \put(0,0){\includegraphics[width=\unitlength,page=1]{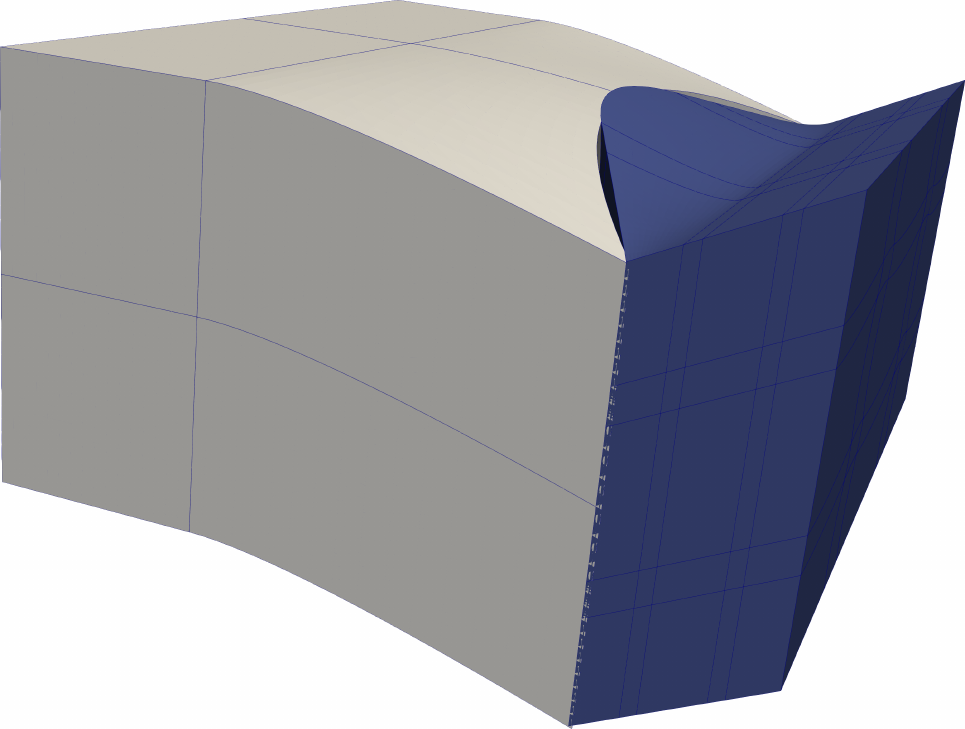}}%
    \put(0.46,-0.12){\color[rgb]{0,0,0}\makebox(0,0)[lt]{\smash{\begin{tabular}[t]{l}(d)\end{tabular}}}}%
  \end{picture}%
\endgroup%

  			\qquad
 				\def\svgwidth{.61\columnwidth}
  			%!TEX root = ../paper.tex
%% Creator: Inkscape 1.1 (c68e22c387, 2021-05-23), www.inkscape.org
%% Accompanies image file 'weak_warp1.pdf' (pdf, eps, ps)
\begingroup%
  \makeatletter%
  \providecommand\color[2][]{%
    \errmessage{(Inkscape) Color is used for the text in Inkscape, but the package 'color.sty' is not loaded}%
    \renewcommand\color[2][]{}%
  }%
  \providecommand\transparent[1]{%
    \errmessage{(Inkscape) Transparency is used (non-zero) for the text in Inkscape, but the package 'transparent.sty' is not loaded}%
    \renewcommand\transparent[1]{}%
  }%
  \providecommand\rotatebox[2]{#2}%
  \newcommand*\fsize{\dimexpr\f@size pt\relax}%
  \newcommand*\lineheight[1]{\fontsize{\fsize}{#1\fsize}\selectfont}%
  \ifx\svgwidth\undefined%
    \setlength{\unitlength}{463.13176332bp}%
    \ifx\svgscale\undefined%
      \relax%
    \else%
      \setlength{\unitlength}{\unitlength * \real{\svgscale}}%
    \fi%
  \else%
    \setlength{\unitlength}{\svgwidth}%
  \fi%
  \global\let\svgwidth\undefined%
  \global\let\svgscale\undefined%
  \makeatother%
  \begin{picture}(1,0.8141757)%
    \lineheight{1}%
    \setlength\tabcolsep{0pt}%
    \put(0,0){\includegraphics[width=\unitlength,page=1]{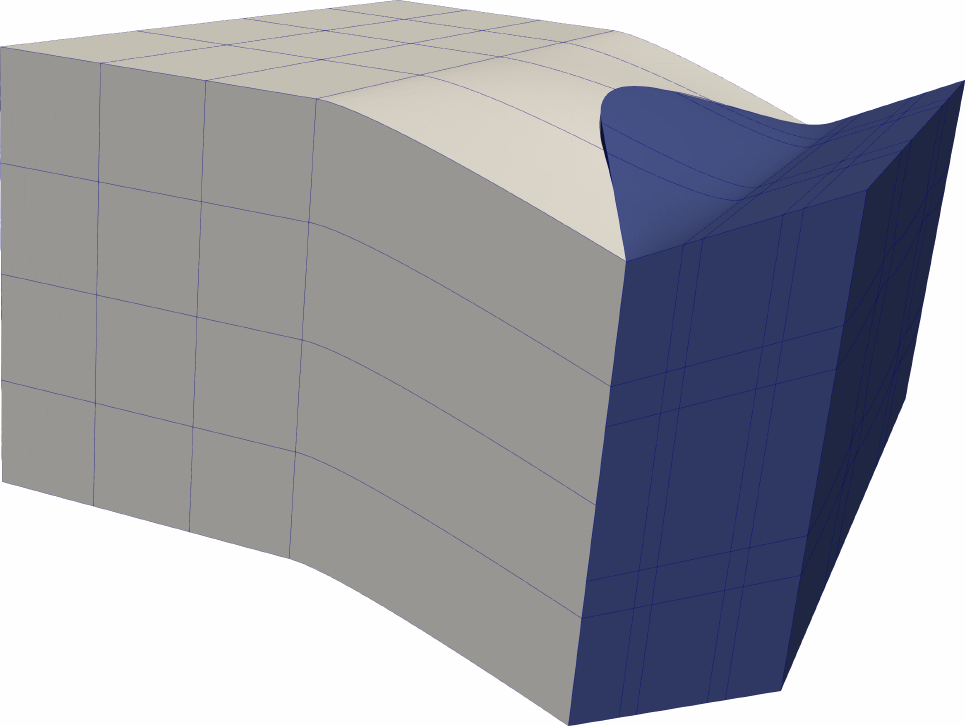}}%
    \put(0.46,-0.12){\color[rgb]{0,0,0}\makebox(0,0)[lt]{\smash{\begin{tabular}[t]{l}(e)\end{tabular}}}}%
  \end{picture}%
\endgroup%

  			\qquad
 				\def\svgwidth{.61\columnwidth}
  			%!TEX root = ../paper.tex
%% Creator: Inkscape 1.1 (c68e22c387, 2021-05-23), www.inkscape.org
%% Accompanies image file 'weak_warp2.pdf' (pdf, eps, ps)
\begingroup%
  \makeatletter%
  \providecommand\color[2][]{%
    \errmessage{(Inkscape) Color is used for the text in Inkscape, but the package 'color.sty' is not loaded}%
    \renewcommand\color[2][]{}%
  }%
  \providecommand\transparent[1]{%
    \errmessage{(Inkscape) Transparency is used (non-zero) for the text in Inkscape, but the package 'transparent.sty' is not loaded}%
    \renewcommand\transparent[1]{}%
  }%
  \providecommand\rotatebox[2]{#2}%
  \newcommand*\fsize{\dimexpr\f@size pt\relax}%
  \newcommand*\lineheight[1]{\fontsize{\fsize}{#1\fsize}\selectfont}%
  \ifx\svgwidth\undefined%
    \setlength{\unitlength}{463.13176332bp}%
    \ifx\svgscale\undefined%
      \relax%
    \else%
      \setlength{\unitlength}{\unitlength * \real{\svgscale}}%
    \fi%
  \else%
    \setlength{\unitlength}{\svgwidth}%
  \fi%
  \global\let\svgwidth\undefined%
  \global\let\svgscale\undefined%
  \makeatother%
  \begin{picture}(1,0.8141757)%
    \lineheight{1}%
    \setlength\tabcolsep{0pt}%
    \put(0,0){\includegraphics[width=\unitlength,page=1]{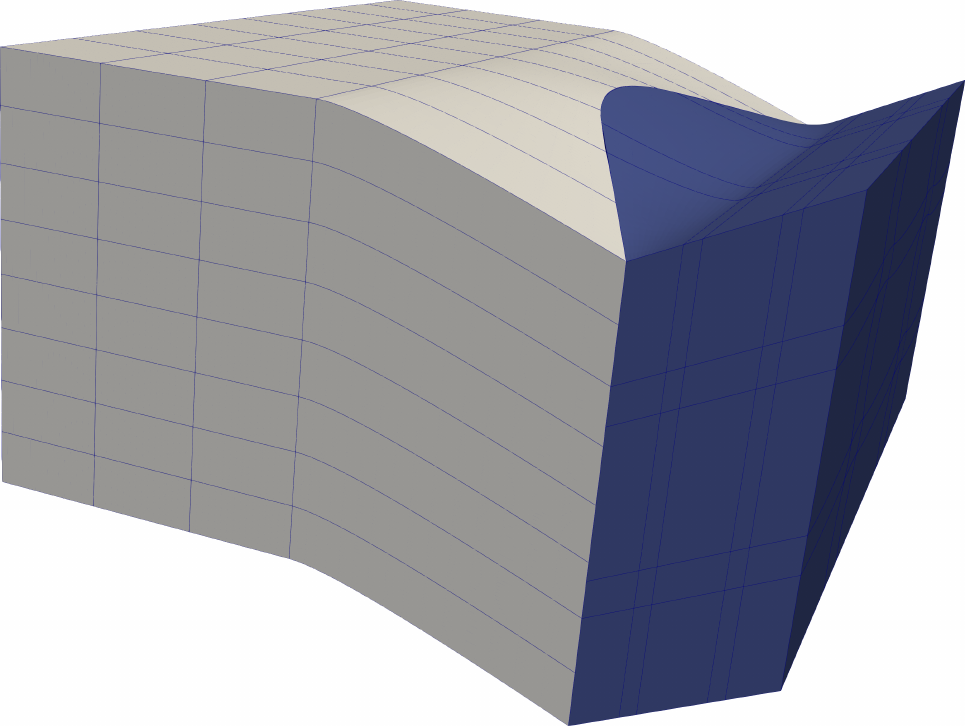}}%
    \put(0.46,-0.12){\color[rgb]{0,0,0}\makebox(0,0)[lt]{\smash{\begin{tabular}[t]{l}(f)\end{tabular}}}}%
  \end{picture}%
\endgroup%

  		}
  		\vspace{0.2cm}
  		\caption{
  			Example described in Section~\ref{sub:weak_continuity}. The blue body in Figure~\ref{fig:example-wc-res}, $\Omega_2$, is
  			deformed as shown in (d), while the solid $\Omega_1$ is deformed accordingly, using the Lee--Lyche--M{\o}rken quasi-interpolant.
  			The influence of the distortion of $\Omega_2$ on $\Omega_1$ is shown in~(a), (b) and (c) for different refinements of the mesh of
  			$\Omega_1$. The corresponding deformations are instead visible in (d), (e) and~(f), respectively.
  		}
  		\label{fig:continuity}
		\end{figure*}

		Figure~\ref{fig:continuity}~(a) shows the influence that a displacement $\delta T_2$ has on $\Omega_1$ by representing $\delta T_1$ and $\delta T_2$ as scalar fields with values between $[0,1]$, where 0 represents no deformation and 1 represents maximal deformation. As expected, the maximal deformation of $\Omega_1$ is localized around the interface $\Gamma_{1,2}$. Nevertheless, when we deform $\Omega_1$ and $\Omega_2$ according to the respective displacements, it is visible a discrepancy between the resulting deformed models, see Figure~\ref{fig:continuity}~(d). In order to improve this result we refine the mesh of $T_1$ in two different steps. In the first step we insert a new knot in the middle of each span for each knot vector of $T_1$ in the three parametric directions, while, in the second, we perform the same operation but only for two parametric directions, see Figures~\ref{fig:continuity}~(b) and~(c), respectively. The results of this procedure are visible in Figures~\ref{fig:continuity}~(e) and~(f), respectively. We notice that, already after the first
		refinement, the behavior of $\delta T_1$ reproduces much closely the distortion $\delta T_2$ and improves even further with the next
		refinement step.
	\subsection{Weak continuity with non-conforming interface} \label{ssub:weak_continuity_not_conf}% (fold)
\begin{figure*}[t]\centering
			\subfloat[] {
    		\def\svgwidth{0.9\columnwidth}
    		%!TEX root = ../notes.tex
%% Creator: Inkscape 1.1 (c68e22c387, 2021-05-23), www.inkscape.org
%% Accompanies image file 'example2_model.pdf' (pdf, eps, ps)
\begingroup%
  \makeatletter%
  \providecommand\color[2][]{%
    \errmessage{(Inkscape) Color is used for the text in Inkscape, but the package 'color.sty' is not loaded}%
    \renewcommand\color[2][]{}%
  }%
  \providecommand\transparent[1]{%
    \errmessage{(Inkscape) Transparency is used (non-zero) for the text in Inkscape, but the package 'transparent.sty' is not loaded}%
    \renewcommand\transparent[1]{}%
  }%
  \providecommand\rotatebox[2]{#2}%
  \newcommand*\fsize{\dimexpr\f@size pt\relax}%
  \newcommand*\lineheight[1]{\fontsize{\fsize}{#1\fsize}\selectfont}%
  \ifx\svgwidth\undefined%
    \setlength{\unitlength}{7410.10772945bp}%
    \ifx\svgscale\undefined%
      \relax%
    \else%
      \setlength{\unitlength}{\unitlength * \real{\svgscale}}%
    \fi%
  \else%
    \setlength{\unitlength}{\svgwidth}%
  \fi%
  \global\let\svgwidth\undefined%
  \global\let\svgscale\undefined%
  \makeatother%
  \begin{picture}(1,0.62576123)%
    \lineheight{1}%
    \setlength\tabcolsep{0pt}%
    \put(0,0){\includegraphics[width=\unitlength,page=1]{img/example2_model.pdf}}%
    \put(0.7,0.04){\color[rgb]{0,0,0}\makebox(0,0)[lt]{\lineheight{1.25}\smash{\begin{tabular}[t]{l}$\Omega_2$\end{tabular}}}}%
    \put(0.3,0.62){\color[rgb]{0,0,0}\makebox(0,0)[lt]{\lineheight{1.25}\smash{\begin{tabular}[t]{l}$\Omega_1$\end{tabular}}}}%
  \end{picture}%
\endgroup%

   		}\quad
   		\subfloat[] {
    		\def\svgwidth{0.9\columnwidth}
    		%!TEX root = ../notes.tex
%% Creator: Inkscape 1.1 (c68e22c387, 2021-05-23), www.inkscape.org
%% Accompanies image file 'example22.pdf' (pdf, eps, ps)
\begingroup%
  \makeatletter%
  \providecommand\color[2][]{%
    \errmessage{(Inkscape) Color is used for the text in Inkscape, but the package 'color.sty' is not loaded}%
    \renewcommand\color[2][]{}%
  }%
  \providecommand\transparent[1]{%
    \errmessage{(Inkscape) Transparency is used (non-zero) for the text in Inkscape, but the package 'transparent.sty' is not loaded}%
    \renewcommand\transparent[1]{}%
  }%
  \providecommand\rotatebox[2]{#2}%
  \newcommand*\fsize{\dimexpr\f@size pt\relax}%
  \newcommand*\lineheight[1]{\fontsize{\fsize}{#1\fsize}\selectfont}%
  \ifx\svgwidth\undefined%
    \setlength{\unitlength}{586.01517932bp}%
    \ifx\svgscale\undefined%
      \relax%
    \else%
      \setlength{\unitlength}{\unitlength * \real{\svgscale}}%
    \fi%
  \else%
    \setlength{\unitlength}{\svgwidth}%
  \fi%
  \global\let\svgwidth\undefined%
  \global\let\svgscale\undefined%
  \makeatother%
  \begin{picture}(1,1.09473946)%
    \lineheight{1}%
    \setlength\tabcolsep{0pt}%
    \put(0,0){\includegraphics[width=\unitlength,page=1]{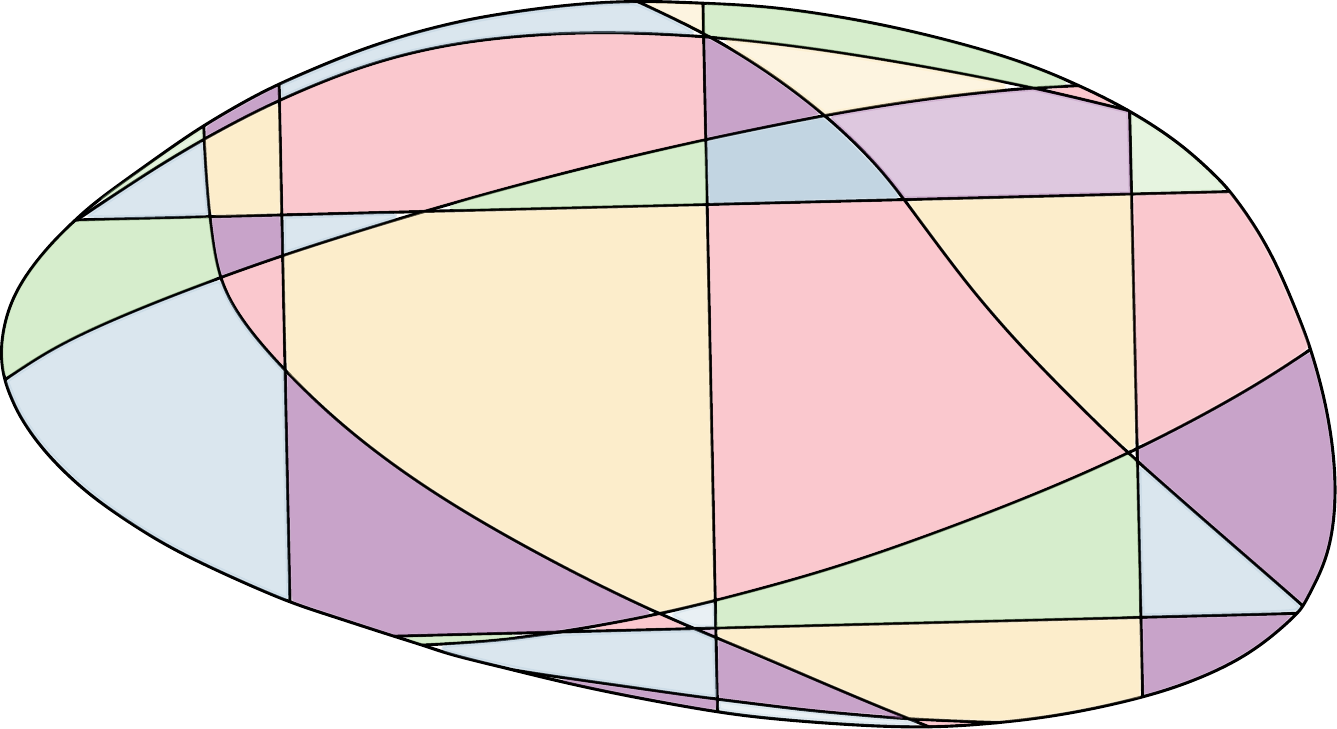}}%
    \put(0.1, 0.53){\color[rgb]{0,0,0}\makebox(0,0)[lt]{\lineheight{1.25}\smash{\begin{tabular}[t]{l}$\hat\Gamma_{1,2}$\end{tabular}}}}%
    % \put(0.15,-0.2){\color[rgb]{0,0,0}\makebox(0,0)[lt]{\smash{\begin{tabular}[t]{l}(b)\end{tabular}}}}%
  \end{picture}%
\endgroup%

   		}
  		\caption{
  			Example described in Section~\ref{ssub:weak_continuity_not_conf}. The original geometries $\Omega_1$ and $\Omega_2$ are
  			shown in~(a). The pull-back curves of the intersection between the isosurfaces of $\Omega_2$ and the interface $\Gamma_{1,2}$ are
  			shown in~(b). The regions to be extracted are marked with different colors.
  		}
  		\label{fig:example-wc-not_conf}
		\end{figure*}

		Let us now consider the situation in Figure~\ref{fig:example-wc-not_conf}. If the interface $\Gamma_{1,2}$ is not the image of a whole
		face of both $\Omega_1$ and $\Omega_2$, we say that the solids are in a non-conforming geometrical setting.

		Enforcing weak continuity constraints in such cases is more complicated due to the fact that the function $\delta T_2$ that
		needs to be projected can have discontinuities along the trimming curve of $\Gamma_{1,2}$. Using the same approach as in 
		Section~\ref{sub:weak_continuity} would therefore results in oscillations in the final results due to the Gibbs
		phenomenon~\cite{Gibbs:1898:FS,Gibbs:1899:FS, Wilbraham:1848:OAC}. Gibbs phenomenon could be avoided by carefully choosing the element
		$K_\ell$ in which to perform the $\mathcal{L}^2$ projection but at the cost of obtaining a `block', pixelized behavior.

		In order to prevent this unpleasant effect, we here propose an approach that enforces weak continuity among the two solids with a
		mass lumping strategy.
	
		Figure~\ref{fig:example-wc-not_conf} shows the setting of this numerical experiment. In Figure~\ref{fig:example-wc-not_conf}~(a), the
		domains $\Omega^*_1$ and $\Omega^*_2$ are shown, together with the isoparametric surfaces of $\Omega^*_2$. These surfaces are
		intersected with the interface $\Gamma_{1,2}$ and produce six intersection curves, whose pull-backs, together with the parametric grid 
		inherited by $\Gamma_{1,2}$ are shown in Figure~\ref{fig:example-wc-not_conf}~(b).

		In order to impose weak continuity constraints in this setting, we introduce the auxiliary function
		\[
			\delta T_{1,2} = 	\begin{cases}
													\delta T_2, &\text{inside } \hat\Gamma_{1,2}\\
													0, 					&\text{otherwise},
												\end{cases}
		\]
		and the domain
		\[
			\Theta = \bigcup_{i\in\Lambda} \text{supp}(B_{i,\dd_1})\subseteq \hat\Gamma_{1,2},
		\]
		where
		\[
			\Lambda = \{i \colon \hat\Gamma_{1,2} \cap \text{supp}(B_{i,\dd_1})\ne\emptyset\}.
		\]
		The spline-based level set function~\cite{Verhoosel:2015:IBG} is then defined as
		\begin{equation}
			\label{eq:convolution}
			\delta T_1 = \sum_{i = 0}^{n_1} B_{i,\dd_1} p_i,
		\end{equation}
		where $n_1$ is the number of control points of $M_1$ and
		\begin{equation}
			\label{eq:convoluted-ctlpt}
			p_i = \frac{\int_{\Theta} B_{i,\dd_1} \delta T_{1,2}}{\int_{\Theta} B_{i,\dd_1}}.
		\end{equation}

		There are several reasons for imposing weak continuity constraints using a spline-based level set approach. Using similar arguments
		as the one proposed in~\cite{Verhoosel:2015:IBG}, it can be shown that $\delta T_1$ satisfies
		\[
			\int_\Theta \delta T_1 = \int_\Theta \delta T_{1,2},
		\]
		and therefore $\delta T_1$ preserves the average value of $\delta T_{1,2}$ over $\Theta$. Moreover it can be shown~\cite{Verhoosel:2015:IBG}
		that~\eqref{eq:convolution} is bounded from above and from below by the maximum and minimum of $\delta T_{1,2}$ and therefore no wild
		oscillations due to the Gibbs phenomenon are to be expected. Finally, the computation of~\eqref{eq:convolution} does not necessitate
		the resolution of a linear system as defined in Section~\ref{sub:weak_continuity} and therefore the new material specification can be
		computed much more efficiently.

		In order to compute the $i$-th control point we first note that the denominator in~\eqref{eq:convoluted-ctlpt} can be computed using
		standard numerical techniques. As for the numerator instead, we write
		\begin{equation}
			\label{eq:blend-numerator}
			\int_{\Theta} B_{i,\dd_1} \delta T_{1,2} = \int_{K_i \cap \hat\Gamma_{1,2}} B_{i,\dd_1} \delta T_2,
		\end{equation}
		with $K_i = \text{supp} (B_{i,\dd_1})$. We notice that the right hand-side integral in~\eqref{eq:blend-numerator} can be computed using
		a similar technique as used in Section~\ref{sub:weak_continuity}.

\begin{figure*}[t]\centering
			\mbox {
    		\def\svgwidth{.61\columnwidth}
    		%!TEX root = ../paper.tex
%% Creator: Inkscape 1.1 (c68e22c387, 2021-05-23), www.inkscape.org
%% Accompanies image file 'weak_conf1.pdf' (pdf, eps, ps)
\begingroup%
  \makeatletter%
  \providecommand\color[2][]{%
    \errmessage{(Inkscape) Color is used for the text in Inkscape, but the package 'color.sty' is not loaded}%
    \renewcommand\color[2][]{}%
  }%
  \providecommand\transparent[1]{%
    \errmessage{(Inkscape) Transparency is used (non-zero) for the text in Inkscape, but the package 'transparent.sty' is not loaded}%
    \renewcommand\transparent[1]{}%
  }%
  \providecommand\rotatebox[2]{#2}%
  \newcommand*\fsize{\dimexpr\f@size pt\relax}%
  \newcommand*\lineheight[1]{\fontsize{\fsize}{#1\fsize}\selectfont}%
  \ifx\svgwidth\undefined%
    \setlength{\unitlength}{463.13176332bp}%
    \ifx\svgscale\undefined%
      \relax%
    \else%
      \setlength{\unitlength}{\unitlength * \real{\svgscale}}%
    \fi%
  \else%
    \setlength{\unitlength}{\svgwidth}%
  \fi%
  \global\let\svgwidth\undefined%
  \global\let\svgscale\undefined%
  \makeatother%
  \begin{picture}(1,0.97917571)%
    \lineheight{1}%
    \setlength\tabcolsep{0pt}%
    \put(0,0){\includegraphics[width=\unitlength,page=1]{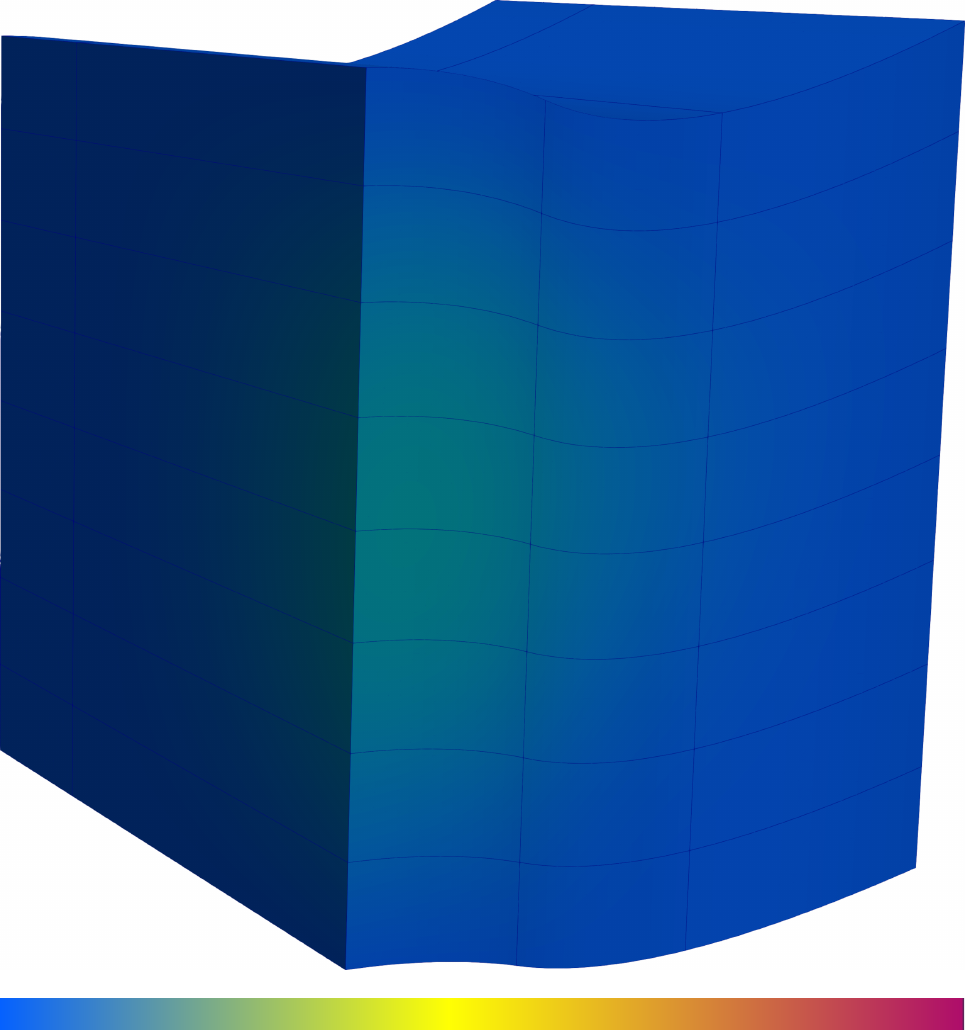}}%
    \put(-0.00450428,-0.07){\color[rgb]{0,0,0}\makebox(0,0)[lt]{\lineheight{1.25}\smash{\begin{tabular}[t]{l}$0$\end{tabular}}}}%
    \put(0.95,-0.07){\color[rgb]{0,0,0}\makebox(0,0)[lt]{\lineheight{1.25}\smash{\begin{tabular}[t]{l}$1.0$\end{tabular}}}}%
    \put(0.46,-0.12){\color[rgb]{0,0,0}\makebox(0,0)[lt]{\smash{\begin{tabular}[t]{l}(a)\end{tabular}}}}%
  \end{picture}%
\endgroup%

    		\qquad
    		\def\svgwidth{.61\columnwidth}
    		%!TEX root = ../paper.tex
%% Creator: Inkscape 1.1 (c68e22c387, 2021-05-23), www.inkscape.org
%% Accompanies image file 'weak_conf1.pdf' (pdf, eps, ps)
\begingroup%
  \makeatletter%
  \providecommand\color[2][]{%
    \errmessage{(Inkscape) Color is used for the text in Inkscape, but the package 'color.sty' is not loaded}%
    \renewcommand\color[2][]{}%
  }%
  \providecommand\transparent[1]{%
    \errmessage{(Inkscape) Transparency is used (non-zero) for the text in Inkscape, but the package 'transparent.sty' is not loaded}%
    \renewcommand\transparent[1]{}%
  }%
  \providecommand\rotatebox[2]{#2}%
  \newcommand*\fsize{\dimexpr\f@size pt\relax}%
  \newcommand*\lineheight[1]{\fontsize{\fsize}{#1\fsize}\selectfont}%
  \ifx\svgwidth\undefined%
    \setlength{\unitlength}{463.13176332bp}%
    \ifx\svgscale\undefined%
      \relax%
    \else%
      \setlength{\unitlength}{\unitlength * \real{\svgscale}}%
    \fi%
  \else%
    \setlength{\unitlength}{\svgwidth}%
  \fi%
  \global\let\svgwidth\undefined%
  \global\let\svgscale\undefined%
  \makeatother%
  \begin{picture}(1,0.97917571)%
    \lineheight{1}%
    \setlength\tabcolsep{0pt}%
    \put(0,0){\includegraphics[width=\unitlength,page=1]{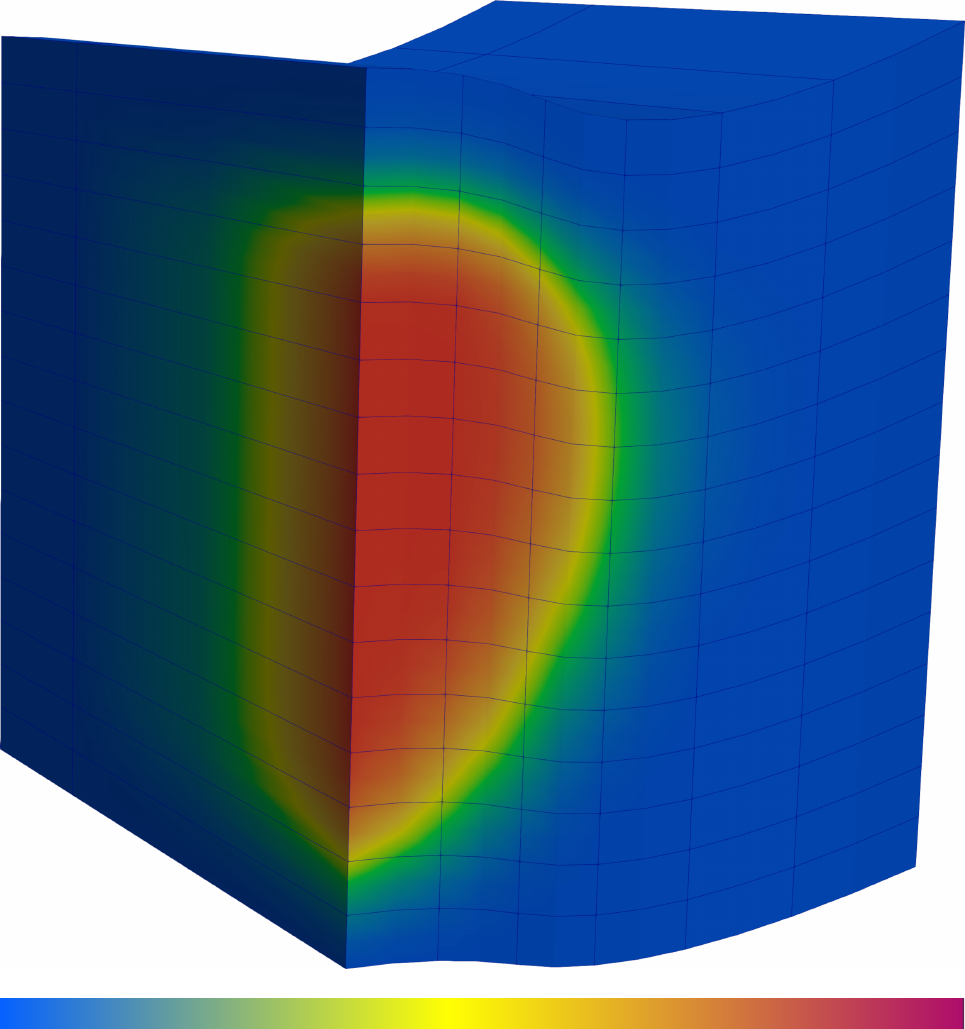}}%
    \put(-0.00450428,-0.07){\color[rgb]{0,0,0}\makebox(0,0)[lt]{\lineheight{1.25}\smash{\begin{tabular}[t]{l}$0$\end{tabular}}}}%
    \put(0.95,-0.07){\color[rgb]{0,0,0}\makebox(0,0)[lt]{\lineheight{1.25}\smash{\begin{tabular}[t]{l}$1.0$\end{tabular}}}}%
    \put(0.46,-0.12){\color[rgb]{0,0,0}\makebox(0,0)[lt]{\smash{\begin{tabular}[t]{l}(b)\end{tabular}}}}%
  \end{picture}%
\endgroup%

    		\qquad
    		\def\svgwidth{.61\columnwidth}
    		%!TEX root = ../paper.tex
%% Creator: Inkscape 1.1 (c68e22c387, 2021-05-23), www.inkscape.org
%% Accompanies image file 'weak_conf1.pdf' (pdf, eps, ps)
\begingroup%
  \makeatletter%
  \providecommand\color[2][]{%
    \errmessage{(Inkscape) Color is used for the text in Inkscape, but the package 'color.sty' is not loaded}%
    \renewcommand\color[2][]{}%
  }%
  \providecommand\transparent[1]{%
    \errmessage{(Inkscape) Transparency is used (non-zero) for the text in Inkscape, but the package 'transparent.sty' is not loaded}%
    \renewcommand\transparent[1]{}%
  }%
  \providecommand\rotatebox[2]{#2}%
  \newcommand*\fsize{\dimexpr\f@size pt\relax}%
  \newcommand*\lineheight[1]{\fontsize{\fsize}{#1\fsize}\selectfont}%
  \ifx\svgwidth\undefined%
    \setlength{\unitlength}{463.13176332bp}%
    \ifx\svgscale\undefined%
      \relax%
    \else%
      \setlength{\unitlength}{\unitlength * \real{\svgscale}}%
    \fi%
  \else%
    \setlength{\unitlength}{\svgwidth}%
  \fi%
  \global\let\svgwidth\undefined%
  \global\let\svgscale\undefined%
  \makeatother%
  \begin{picture}(1,0.97917571)%
    \lineheight{1}%
    \setlength\tabcolsep{0pt}%
    \put(0,0){\includegraphics[width=\unitlength,page=1]{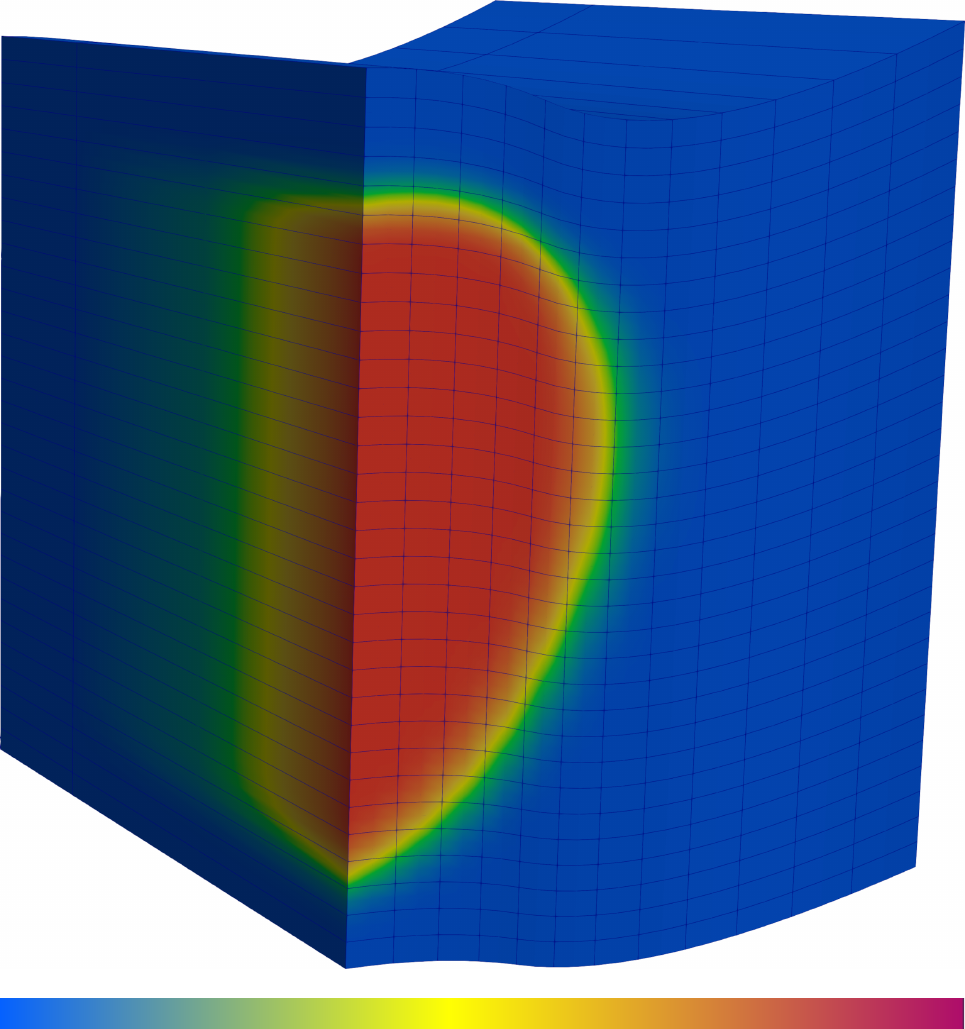}}%
    \put(-0.00450428,-0.07){\color[rgb]{0,0,0}\makebox(0,0)[lt]{\lineheight{1.25}\smash{\begin{tabular}[t]{l}$0$\end{tabular}}}}%
    \put(0.95,-0.07){\color[rgb]{0,0,0}\makebox(0,0)[lt]{\lineheight{1.25}\smash{\begin{tabular}[t]{l}$1.0$\end{tabular}}}}%
    \put(0.46,-0.12){\color[rgb]{0,0,0}\makebox(0,0)[lt]{\smash{\begin{tabular}[t]{l}(c)\end{tabular}}}}%
  \end{picture}%
\endgroup%

   		}\\
   		\vspace{0.1cm}
 			\mbox {
	 			\def\svgwidth{.61\columnwidth}
  			%!TEX root = ../paper.tex
%% Creator: Inkscape 1.1 (c68e22c387, 2021-05-23), www.inkscape.org
%% Accompanies image file 'weak_warp0.pdf' (pdf, eps, ps)
\begingroup%
  \makeatletter%
  \providecommand\color[2][]{%
    \errmessage{(Inkscape) Color is used for the text in Inkscape, but the package 'color.sty' is not loaded}%
    \renewcommand\color[2][]{}%
  }%
  \providecommand\transparent[1]{%
    \errmessage{(Inkscape) Transparency is used (non-zero) for the text in Inkscape, but the package 'transparent.sty' is not loaded}%
    \renewcommand\transparent[1]{}%
  }%
  \providecommand\rotatebox[2]{#2}%
  \newcommand*\fsize{\dimexpr\f@size pt\relax}%
  \newcommand*\lineheight[1]{\fontsize{\fsize}{#1\fsize}\selectfont}%
  \ifx\svgwidth\undefined%
    \setlength{\unitlength}{463.13176332bp}%
    \ifx\svgscale\undefined%
      \relax%
    \else%
      \setlength{\unitlength}{\unitlength * \real{\svgscale}}%
    \fi%
  \else%
    \setlength{\unitlength}{\svgwidth}%
  \fi%
  \global\let\svgwidth\undefined%
  \global\let\svgscale\undefined%
  \makeatother%
  \begin{picture}(1,0.8173757)%
    \lineheight{1}%
    \setlength\tabcolsep{0pt}%
    \put(0,0){\includegraphics[width=\unitlength,page=1]{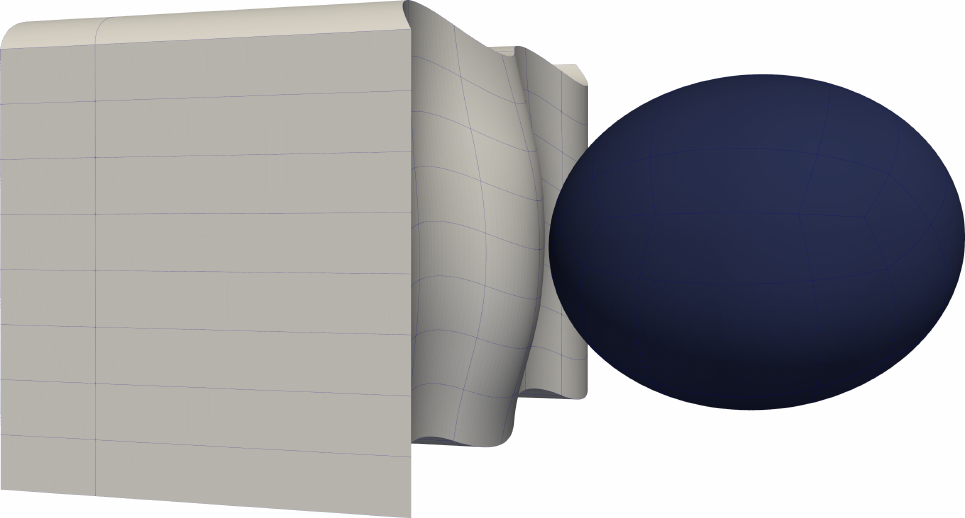}}%
    \put(0.46,-0.12){\color[rgb]{0,0,0}\makebox(0,0)[lt]{\smash{\begin{tabular}[t]{l}(d)\end{tabular}}}}%
  \end{picture}%
\endgroup%

  			\qquad
 				\def\svgwidth{.61\columnwidth}
  			%!TEX root = ../paper.tex
%% Creator: Inkscape 1.1 (c68e22c387, 2021-05-23), www.inkscape.org
%% Accompanies image file 'weak_warp2.pdf' (pdf, eps, ps)
\begingroup%
  \makeatletter%
  \providecommand\color[2][]{%
    \errmessage{(Inkscape) Color is used for the text in Inkscape, but the package 'color.sty' is not loaded}%
    \renewcommand\color[2][]{}%
  }%
  \providecommand\transparent[1]{%
    \errmessage{(Inkscape) Transparency is used (non-zero) for the text in Inkscape, but the package 'transparent.sty' is not loaded}%
    \renewcommand\transparent[1]{}%
  }%
  \providecommand\rotatebox[2]{#2}%
  \newcommand*\fsize{\dimexpr\f@size pt\relax}%
  \newcommand*\lineheight[1]{\fontsize{\fsize}{#1\fsize}\selectfont}%
  \ifx\svgwidth\undefined%
    \setlength{\unitlength}{463.13176332bp}%
    \ifx\svgscale\undefined%
      \relax%
    \else%
      \setlength{\unitlength}{\unitlength * \real{\svgscale}}%
    \fi%
  \else%
    \setlength{\unitlength}{\svgwidth}%
  \fi%
  \global\let\svgwidth\undefined%
  \global\let\svgscale\undefined%
  \makeatother%
  \begin{picture}(1,0.8141757)%
    \lineheight{1}%
    \setlength\tabcolsep{0pt}%
    \put(0,0){\includegraphics[width=\unitlength,page=1]{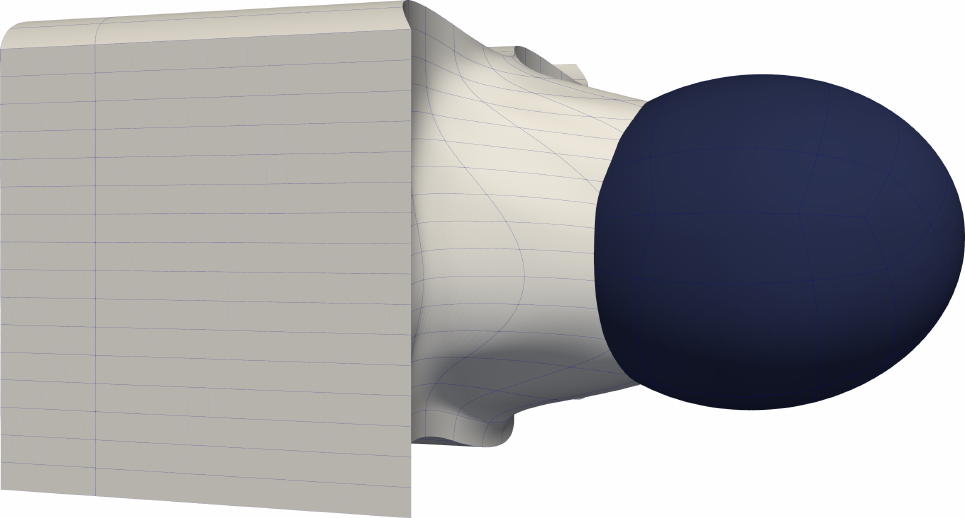}}%
    \put(0.46,-0.12){\color[rgb]{0,0,0}\makebox(0,0)[lt]{\smash{\begin{tabular}[t]{l}(e)\end{tabular}}}}%
  \end{picture}%
\endgroup%

  			\qquad
 				\def\svgwidth{.61\columnwidth}
  			%!TEX root = ../paper.tex
%% Creator: Inkscape 1.1 (c68e22c387, 2021-05-23), www.inkscape.org
%% Accompanies image file 'weak_warp2.pdf' (pdf, eps, ps)
\begingroup%
  \makeatletter%
  \providecommand\color[2][]{%
    \errmessage{(Inkscape) Color is used for the text in Inkscape, but the package 'color.sty' is not loaded}%
    \renewcommand\color[2][]{}%
  }%
  \providecommand\transparent[1]{%
    \errmessage{(Inkscape) Transparency is used (non-zero) for the text in Inkscape, but the package 'transparent.sty' is not loaded}%
    \renewcommand\transparent[1]{}%
  }%
  \providecommand\rotatebox[2]{#2}%
  \newcommand*\fsize{\dimexpr\f@size pt\relax}%
  \newcommand*\lineheight[1]{\fontsize{\fsize}{#1\fsize}\selectfont}%
  \ifx\svgwidth\undefined%
    \setlength{\unitlength}{463.13176332bp}%
    \ifx\svgscale\undefined%
      \relax%
    \else%
      \setlength{\unitlength}{\unitlength * \real{\svgscale}}%
    \fi%
  \else%
    \setlength{\unitlength}{\svgwidth}%
  \fi%
  \global\let\svgwidth\undefined%
  \global\let\svgscale\undefined%
  \makeatother%
  \begin{picture}(1,0.8141757)%
    \lineheight{1}%
    \setlength\tabcolsep{0pt}%
    \put(0,0){\includegraphics[width=\unitlength,page=1]{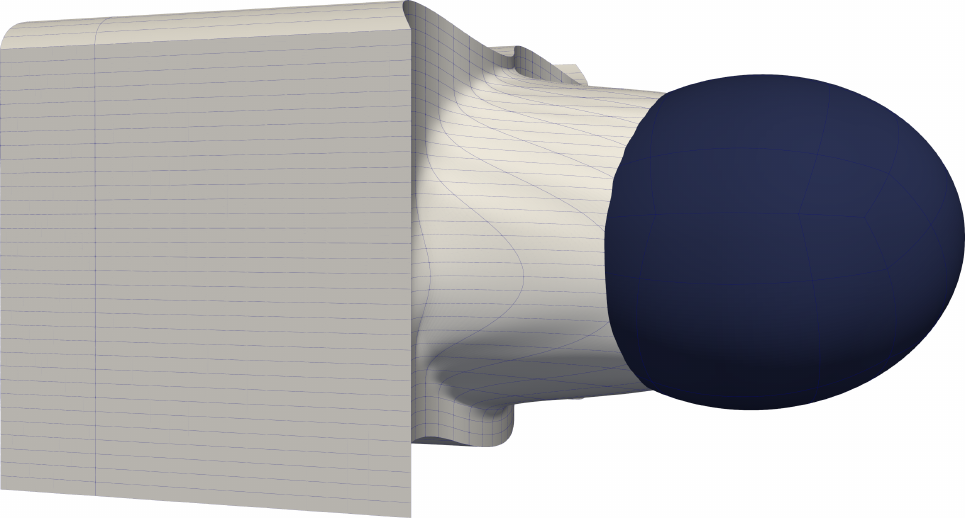}}%
    \put(0.46,-0.12){\color[rgb]{0,0,0}\makebox(0,0)[lt]{\smash{\begin{tabular}[t]{l}(f)\end{tabular}}}}%
  \end{picture}%
\endgroup%

  		}
  		\vspace{0.2cm}
  		\caption{
  			Example described in Section~\ref{ssub:weak_continuity_not_conf}. The blue body in Figure~\ref{fig:example-wc-not_conf},
  			$\Omega_2$, is translated along one of the main axes as shown in (d). The deformation of $\Omega_1$ is computed using the
  			spline-based level set approach. The influence of the distortion of $\Omega_2$ on $\Omega_1$ is shown in~(a), (b), and (c) for
  			different refinements of the mesh of $\Omega_1$. The corresponding deformations are instead visible in (d), (e), and~(f),
  			respectively.
  		}
  		\label{fig:warp_continuity}
		\end{figure*}

		In order to simulate the effect of a deformation of $\Omega_2$, we translate $\Omega_2$ over one of the main directions of a constant quantity. Without loss of generality we assume that
		\[
			\delta T_2 = (0,c,0)^T,\qquad c\in\RR.
		\]
		The effects of this deformation are represented in Figure~\ref{fig:warp_continuity} for the different refinements of $T_1$. As visible
		in Figure~\ref{fig:warp_continuity}~(d), for a very coarse mesh of $\Omega_1$ the convolution-based strategy does not approximate well
		enough the deformation of $\Omega_2$. This is due to the fact that~\eqref{eq:convolution} reproduces the average of the distortion of
		$\delta T_2$ in $\Theta$. From~\eqref{eq:convoluted-ctlpt} it is clear that, in order to get a better approximation of $\delta T_2$,
		the support of the basis functions $K_i$ should be small enough, so to guarantee that the average of $\delta T_2$ on $K_i$ is a good
		approximation of the behavior of $\delta T_2$ over the same domain. In order to improve this result we refine the mesh of $T_1$ in
		two different steps. In
		each step we insert a new knot in the middle of each span for each knot vector of $T_1$ in two of the three parametric directions. This
		is equivalent to splitting each element of $T_1$ into four sub-elements. The results of these refinement steps on the behavior of
		$\delta T_1$ are visible in Figures~\ref{fig:warp_continuity}~(e) and~(f). Notice that, as the mesh gets finer, the influence of the
		distortion gets more localized around the trimming curve of the interface $\Gamma_{1,2}$, as well. 
		Figures~\ref{fig:warp_continuity}~(a)--(c) show instead the influence of the distortion of $\Omega_2$ over the first body, showing that
		only the elements that are closer to the interface are in practice affected by this procedure.

	% % subsection material_blending (end)
	\subsection{Poisson's problem}\label{sub:poisson}
\begin{figure*}[t]\centering
			\mbox {
	 			\def\svgwidth{.60\columnwidth}
  			%!TEX root = ../paper.tex
%% Creator: Inkscape 1.1 (c68e22c387, 2021-05-23), www.inkscape.org
%% Accompanies image file 'numerical_integration_models.pdf' (pdf, eps, ps)
\begingroup%
  \makeatletter%
  \providecommand\color[2][]{%
    \errmessage{(Inkscape) Color is used for the text in Inkscape, but the package 'color.sty' is not loaded}%
    \renewcommand\color[2][]{}%
  }%
  \providecommand\transparent[1]{%
    \errmessage{(Inkscape) Transparency is used (non-zero) for the text in Inkscape, but the package 'transparent.sty' is not loaded}%
    \renewcommand\transparent[1]{}%
  }%
  \providecommand\rotatebox[2]{#2}%
  \newcommand*\fsize{\dimexpr\f@size pt\relax}%
  \newcommand*\lineheight[1]{\fontsize{\fsize}{#1\fsize}\selectfont}%
  \ifx\svgwidth\undefined%
    \setlength{\unitlength}{5661bp}%
    \ifx\svgscale\undefined%
      \relax%
    \else%
      \setlength{\unitlength}{\unitlength * \real{\svgscale}}%
    \fi%
  \else%
    \setlength{\unitlength}{\svgwidth}%
  \fi%
  \global\let\svgwidth\undefined%
  \global\let\svgscale\undefined%
  \makeatother%
  \begin{picture}(1,0.75382535)%
    \lineheight{1}%
    \setlength\tabcolsep{0pt}%
    \put(0,0){\includegraphics[width=\unitlength,page=1]{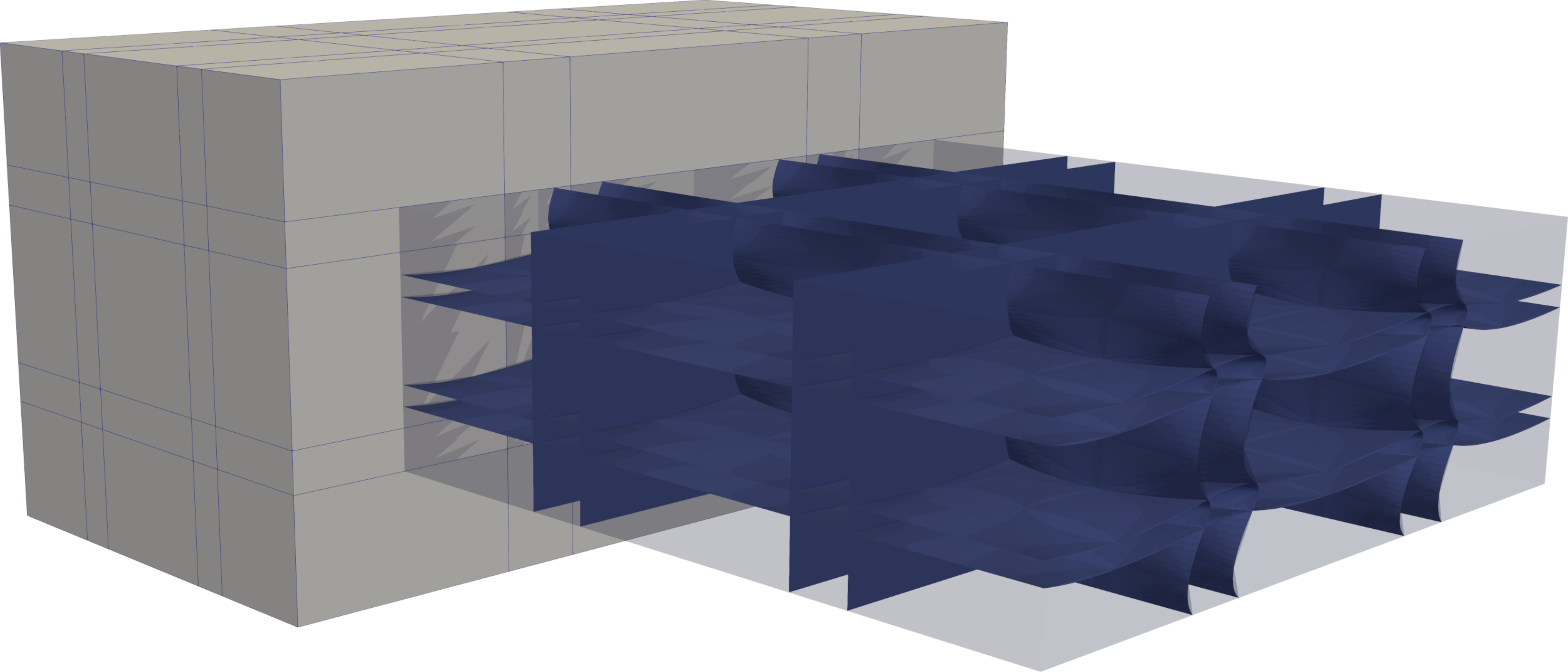}}%
    \put(0.14138577,0.45){\color[rgb]{0,0,0}\makebox(0,0)[lt]{\lineheight{1.25}\smash{\begin{tabular}[t]{l}$\Omega_1^*$\end{tabular}}}}%
    \put(0.46,-0.12){\color[rgb]{0,0,0}\makebox(0,0)[lt]{\smash{\begin{tabular}[t]{l}(a)\end{tabular}}}}%
    \put(0.8,0.0){\color[rgb]{0,0,0}\makebox(0,0)[lt]{\lineheight{1.25}\smash{\begin{tabular}[t]{l}$\Omega_2^*$\end{tabular}}}}%
  \end{picture}%
\endgroup%

  			\qquad
 				\def\svgwidth{.60\columnwidth}
  			%!TEX root = ../paper.tex
%% Creator: Inkscape 1.1 (c68e22c387, 2021-05-23), www.inkscape.org
%% Accompanies image file 'numerical_integration.pdf' (pdf, eps, ps)
\begingroup%
  \makeatletter%
  \providecommand\color[2][]{%
    \errmessage{(Inkscape) Color is used for the text in Inkscape, but the package 'color.sty' is not loaded}%
    \renewcommand\color[2][]{}%
  }%
  \providecommand\transparent[1]{%
    \errmessage{(Inkscape) Transparency is used (non-zero) for the text in Inkscape, but the package 'transparent.sty' is not loaded}%
    \renewcommand\transparent[1]{}%
  }%
  \providecommand\rotatebox[2]{#2}%
  \newcommand*\fsize{\dimexpr\f@size pt\relax}%
  \newcommand*\lineheight[1]{\fontsize{\fsize}{#1\fsize}\selectfont}%
  \ifx\svgwidth\undefined%
    \setlength{\unitlength}{467.03762704bp}%
    \ifx\svgscale\undefined%
      \relax%
    \else%
      \setlength{\unitlength}{\unitlength * \real{\svgscale}}%
    \fi%
  \else%
    \setlength{\unitlength}{\svgwidth}%
  \fi%
  \global\let\svgwidth\undefined%
  \global\let\svgscale\undefined%
  \makeatother%
  \begin{picture}(1,0.95098331)%
    \lineheight{1}%
    \setlength\tabcolsep{0pt}%
    \put(0,0){\includegraphics[width=\unitlength,page=1]{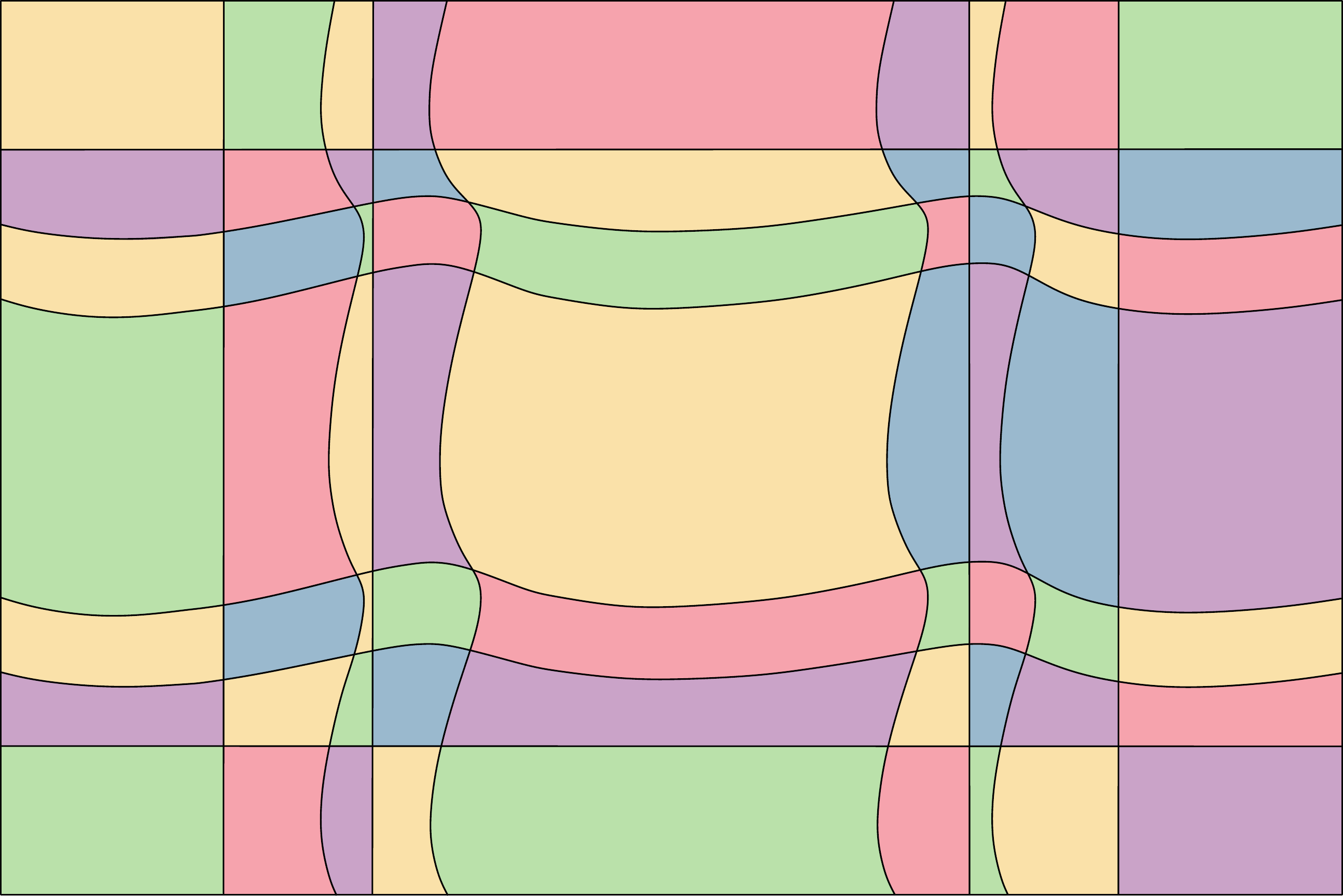}}%
    \put(0.6,0.72){\color[rgb]{0,0,0}\makebox(0,0)[lt]{\lineheight{1.25}\smash{\begin{tabular}[t]{l}$\hat\Gamma_{1,2}$\end{tabular}}}}%
    \put(0.46,-0.12){\color[rgb]{0,0,0}\makebox(0,0)[lt]{\smash{\begin{tabular}[t]{l}(b)\end{tabular}}}}%
  \end{picture}%
\endgroup%

  			\qquad
 				\def\svgwidth{.60\columnwidth}
  			%!TEX root = ../paper.tex
%% Creator: Inkscape 1.1 (c68e22c387, 2021-05-23), www.inkscape.org
%% Accompanies image file 'nd_boundaries.pdf' (pdf, eps, ps)
\begingroup%
  \makeatletter%
  \providecommand\color[2][]{%
    \errmessage{(Inkscape) Color is used for the text in Inkscape, but the package 'color.sty' is not loaded}%
    \renewcommand\color[2][]{}%
  }%
  \providecommand\transparent[1]{%
    \errmessage{(Inkscape) Transparency is used (non-zero) for the text in Inkscape, but the package 'transparent.sty' is not loaded}%
    \renewcommand\transparent[1]{}%
  }%
  \providecommand\rotatebox[2]{#2}%
  \newcommand*\fsize{\dimexpr\f@size pt\relax}%
  \newcommand*\lineheight[1]{\fontsize{\fsize}{#1\fsize}\selectfont}%
  \ifx\svgwidth\undefined%
    \setlength{\unitlength}{700.26198249bp}%
    \ifx\svgscale\undefined%
      \relax%
    \else%
      \setlength{\unitlength}{\unitlength * \real{\svgscale}}%
    \fi%
  \else%
    \setlength{\unitlength}{\svgwidth}%
  \fi%
  \global\let\svgwidth\undefined%
  \global\let\svgscale\undefined%
  \makeatother%
  \begin{picture}(1,0.77887705)%
    \lineheight{1}%
    \setlength\tabcolsep{0pt}%
    \put(0,0){\includegraphics[width=\unitlength,page=1]{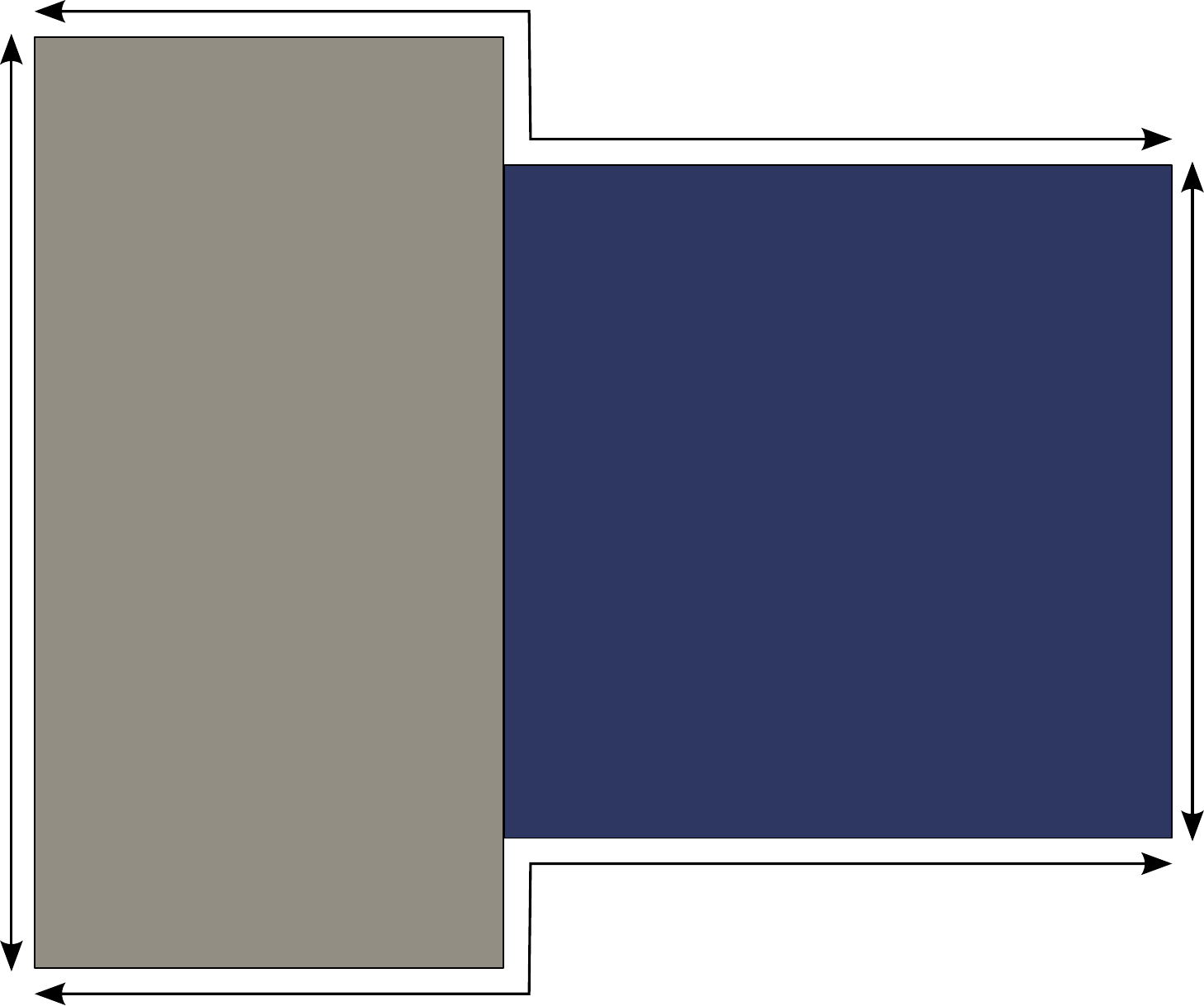}}%
    \put(-0.12,0.4){\color[rgb]{0,0,0}\makebox(0,0)[lt]{\smash{\begin{tabular}[t]{l}$\Gamma_D$\end{tabular}}}}%
    \put(1.02,0.4){\color[rgb]{0,0,0}\makebox(0,0)[lt]{\smash{\begin{tabular}[t]{l}$\Gamma_D$\end{tabular}}}}%
    \put(0.27,0.4){\color[rgb]{0,0,0}\makebox(0,0)[lt]{\smash{\begin{tabular}[t]{l}$\Gamma_{1,2}$\end{tabular}}}}%
    \put(0.55,-0.02){\color[rgb]{0,0,0}\makebox(0,0)[lt]{\smash{\begin{tabular}[t]{l}$\Gamma_N$\end{tabular}}}}%
    \put(0.55,0.8){\color[rgb]{0,0,0}\makebox(0,0)[lt]{\smash{\begin{tabular}[t]{l}$\Gamma_N$\end{tabular}}}}%
    \put(0.46,-0.12){\color[rgb]{0,0,0}\makebox(0,0)[lt]{\smash{\begin{tabular}[t]{l}(c)\end{tabular}}}}%
  \end{picture}%
\endgroup%

  		}
  		\vspace{0.5cm}
  		\caption{
  			Example described in Section~\ref{sub:poisson}. Two volumetric boxes $\Omega^*_1$ and $\Omega^*_2$ in a contact position are
  		 	shown in (a). The intersection curves between the knot surfaces of $T_2$ and $\Gamma_{1,2}$ are pulled-back in the 
  		 	parametric space $\hat\Gamma_{1,2}$ and, together with the knot lines of the interface, form the curvilinear drawing in~(b). 
  		 	The regions to be extracted are marked with different colors. Finally, the definition of the Dirichlet and 
  		 	Neumann boundaries are shown in~(c).
  		}
  		\label{fig:example-contact}
		\end{figure*}

		In this section we use our algorithm to solve the Poisson's equation in a contact problem context. Let us consider the geometric
		setting as in Figure~\ref{fig:example-contact}~(a).
		The two boxes $\Omega_1$ and $\Omega_2$ touch along the interface but do not intersect and are parameterized by two trivariate B-splines
		$T_1\in\mathbb{S}_{\dd_1}$ and $T_2\in\mathbb{S}_{\dd_2}$, where we assume the degrees in $\dd_1$ to satisfy $d^{(1)}_\alpha\ge2$,
		$\alpha = u, v, w$. In this geometric setting, the
		strong form of the Poisson's problem can be stated as follows. Given a domain $\Omega = \Omega_1 \cup \Omega_2$, find
		$u\colon \Omega\rightarrow \RR$ such that
		\begin{equation}
			\label{eq:poisson_problem}
			\begin{aligned}
				-\Delta u_i &= f\qquad \text{in } \Omega_i,\quad i = 1, 2,&\\
				\partial_{\nn_1} u_1 + \partial_{\nn_2} u_2 &= 0\qquad \text{on } \Gamma_{1,2},&\\
				u_1 - u_2 &= 0\qquad \text{on } \Gamma_{1,2},&\\
				u &= 0\qquad \text{on } \Gamma_D,&\\
				\nabla u\cdot\nn &= 0\qquad \text{on } \Gamma_N,&
			\end{aligned}
		\end{equation}
		where $u_i = \restr{u}{\Omega_i}$, $i = 1, 2$, $f \colon \Omega \rightarrow \RR$ and $\nn$ is the outward normal of $\partial \Omega$.
		We assume that
		\begin{equation*}
			\begin{aligned}
				\bar{\Gamma}_D \cup \bar{\Gamma}_N &= \partial \Omega\\
				\Gamma_D \cap \Gamma_N &= \partial \emptyset
			\end{aligned}
		\end{equation*}
		The definition of $\Gamma_D$ and $\Gamma_N$ for the considered volumetric
		model is shown in Figure~\ref{fig:example-contact}~(c) for a horizontal cross section of $\Omega$. The upper and lower faces of
		both $\Omega_1$ and $\Omega_2$ are considered as Neumann boundaries.

		In order to present the weak formulation of~\eqref{eq:poisson_problem}, we first introduce the function spaces
		\begin{align*}
			\vcal_h = &\{v_{1,h}\in\bbs_{\dd_1,\tl_1}\colon \restr{v_{1,h}}{\Gamma_D \cap \partial\Omega_1} = 0\}\oplus\\
								&\{v_{2,h}\in\bbs_{\dd_2,\tl_2}\colon \restr{v_{2,h}}{\Gamma_D \cap \partial\Omega_2} = 0\}
		\end{align*}
		and
		\[
			\Lambda_h = \{\lambda_h\in tr(\bbs_{\dd_1 - 2,\tl_1^\prime})\},
		\]
		where we use the notation $\dd - k = (d_u - k, d_v - k, d_w - k)$ and we denote with $tr(\,\cdot\,)$ the trace operator over
		$\hat\Gamma_{1,2}$ and with $\tl_1^\prime$ a knot vector obtained from $\tl_1$ by removing the first and the last two knots in each
		parametric direction.
		The function space $\Lambda_h$ is where the Lagrange multipliers for the resolution of the Poisson problem are
		going to be defined. This particular choice of $\Lambda_h$ guarantees the inf-sup stability of the Lagrange
		multipliers~\cite{Brivadis:2015:IMM} and is the reason for the hypothesis on the degrees of the trivariate $T_1$.

		Hence, the discrete, weak form of Poisson's problem can be stated as follows. Find $u_h\in\vcal_h$ and $\lambda_h\in\Lambda_h$ such that
		\begin{equation}
			\label{eq:weak_form}
			\begin{aligned}
				a(u_h, v_h) + b(v_h, \lambda_h) &= f(v_h)\qquad &v_h\in \vcal_h\\
				b(u_h,\mu_h) &= 0\qquad &\mu_h\in\Lambda_h,
			\end{aligned}
		\end{equation}
		where
		\begin{align*}
			a(u, v) &= \int_{\Omega} \nabla u \cdot \nabla v,\\
			b(v, \mu) &= \int_{\Gamma_{1, 2}} [\,v\,] \mu,\\
			f(v) &= \int_{\Omega} f v,
		\end{align*}
		and $[\,u\,] = u_1 - u_2$ is the usual jump operator from $\Omega_1$ to $\Omega_2$ over $\Gamma_{1,2}$.

		Denoting with $u_{i,h} = \restr{u_h}{\Omega_i}$, $i = 1, 2$, $b(u_h, \mu_h)$ reads
		\[
			b(u_h, \mu_h) = \int_{\Gamma_{1,2}} \mu_h u_{1,h} - \int_{\Gamma_{1,2}} \mu_h u_{2,h},
		\]
		and therefore, it appears clear that, in order to obtain a good approximation for $b(u_h, \mu_h)$ it is needed to compute the mesh
		intersection between the meshes inherited by $\mu_h$ and $u_{2,h}$.

\begin{figure*}[t]\centering
 			\mbox {
	 			\def\svgwidth{.61\columnwidth}
  			%!TEX root = ../paper.tex
%% Creator: Inkscape 1.1 (c68e22c387, 2021-05-23), www.inkscape.org
%% Accompanies image file 'weak_warp0.pdf' (pdf, eps, ps)
\begingroup%
  \makeatletter%
  \providecommand\color[2][]{%
    \errmessage{(Inkscape) Color is used for the text in Inkscape, but the package 'color.sty' is not loaded}%
    \renewcommand\color[2][]{}%
  }%
  \providecommand\transparent[1]{%
    \errmessage{(Inkscape) Transparency is used (non-zero) for the text in Inkscape, but the package 'transparent.sty' is not loaded}%
    \renewcommand\transparent[1]{}%
  }%
  \providecommand\rotatebox[2]{#2}%
  \newcommand*\fsize{\dimexpr\f@size pt\relax}%
  \newcommand*\lineheight[1]{\fontsize{\fsize}{#1\fsize}\selectfont}%
  \ifx\svgwidth\undefined%
    \setlength{\unitlength}{463.13176332bp}%
    \ifx\svgscale\undefined%
      \relax%
    \else%
      \setlength{\unitlength}{\unitlength * \real{\svgscale}}%
    \fi%
  \else%
    \setlength{\unitlength}{\svgwidth}%
  \fi%
  \global\let\svgwidth\undefined%
  \global\let\svgscale\undefined%
  \makeatother%
  \begin{picture}(1,0.8173757)%
    \lineheight{1}%
    \setlength\tabcolsep{0pt}%
    \put(0,0){\includegraphics[width=\unitlength,page=1]{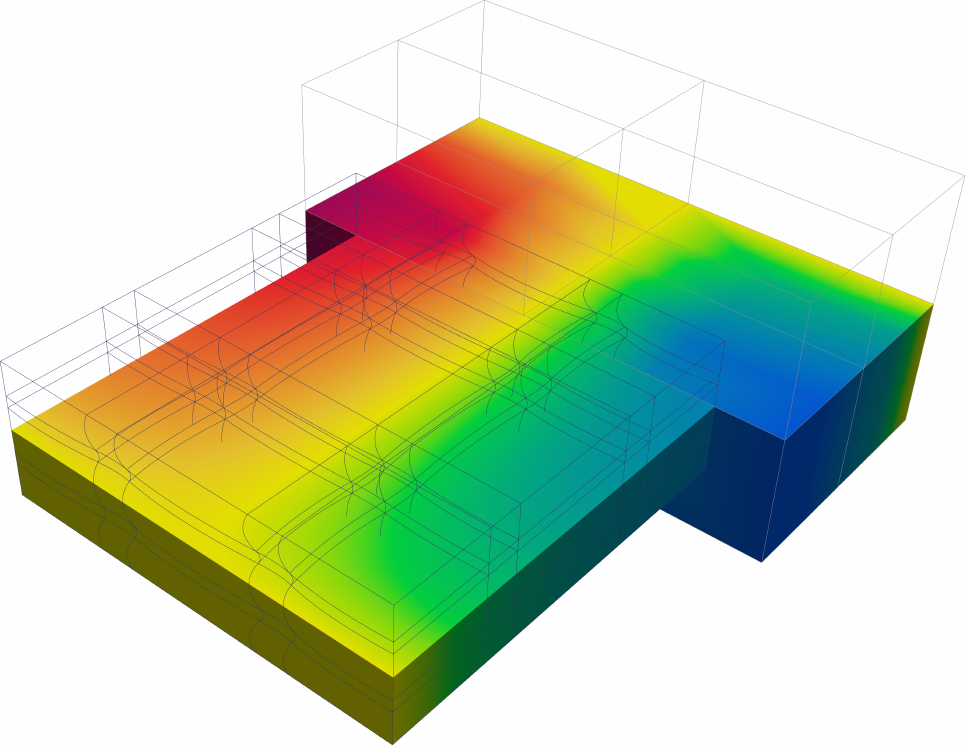}}%
    \put(0.46,-0.12){\color[rgb]{0,0,0}\makebox(0,0)[lt]{\smash{\begin{tabular}[t]{l}(a)\end{tabular}}}}%
  \end{picture}%
\endgroup%

  			\qquad
 				\def\svgwidth{.61\columnwidth}
  			%!TEX root = ../paper.tex
%% Creator: Inkscape 1.1 (c68e22c387, 2021-05-23), www.inkscape.org
%% Accompanies image file 'weak_warp0.pdf' (pdf, eps, ps)
\begingroup%
  \makeatletter%
  \providecommand\color[2][]{%
    \errmessage{(Inkscape) Color is used for the text in Inkscape, but the package 'color.sty' is not loaded}%
    \renewcommand\color[2][]{}%
  }%
  \providecommand\transparent[1]{%
    \errmessage{(Inkscape) Transparency is used (non-zero) for the text in Inkscape, but the package 'transparent.sty' is not loaded}%
    \renewcommand\transparent[1]{}%
  }%
  \providecommand\rotatebox[2]{#2}%
  \newcommand*\fsize{\dimexpr\f@size pt\relax}%
  \newcommand*\lineheight[1]{\fontsize{\fsize}{#1\fsize}\selectfont}%
  \ifx\svgwidth\undefined%
    \setlength{\unitlength}{463.13176332bp}%
    \ifx\svgscale\undefined%
      \relax%
    \else%
      \setlength{\unitlength}{\unitlength * \real{\svgscale}}%
    \fi%
  \else%
    \setlength{\unitlength}{\svgwidth}%
  \fi%
  \global\let\svgwidth\undefined%
  \global\let\svgscale\undefined%
  \makeatother%
  \begin{picture}(1,0.8173757)%
    \lineheight{1}%
    \setlength\tabcolsep{0pt}%
    \put(0,0){\includegraphics[width=\unitlength,page=1]{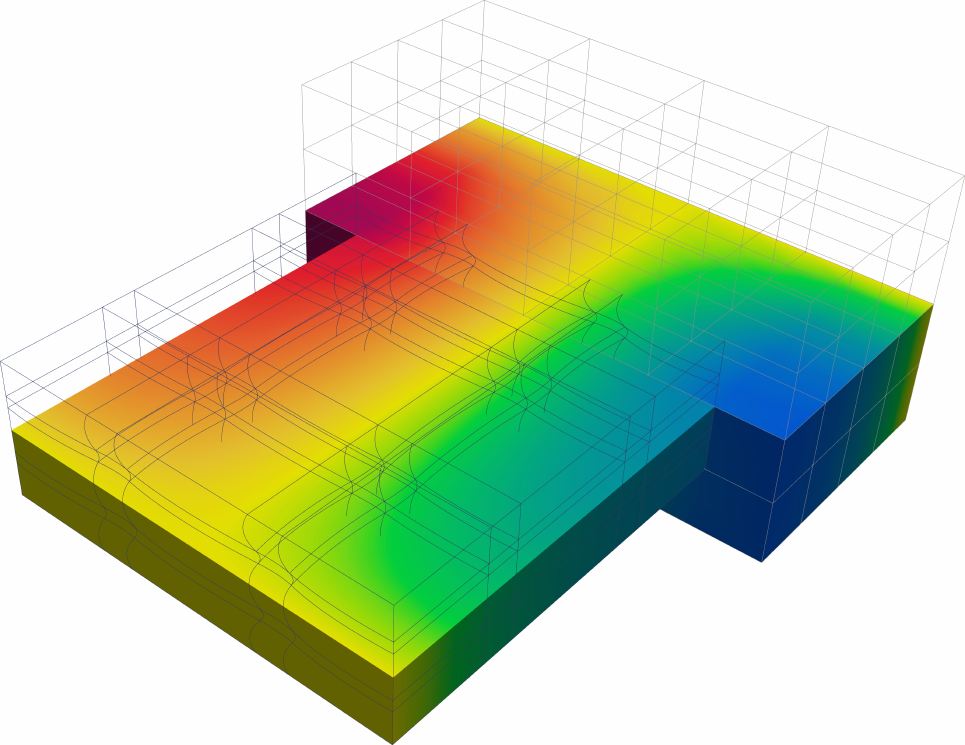}}%
    \put(0.46,-0.12){\color[rgb]{0,0,0}\makebox(0,0)[lt]{\smash{\begin{tabular}[t]{l}(b)\end{tabular}}}}%
  \end{picture}%
\endgroup%

  			\qquad
 				\def\svgwidth{.61\columnwidth}
  			%!TEX root = ../paper.tex
%% Creator: Inkscape 1.1 (c68e22c387, 2021-05-23), www.inkscape.org
%% Accompanies image file 'weak_warp0.pdf' (pdf, eps, ps)
\begingroup%
  \makeatletter%
  \providecommand\color[2][]{%
    \errmessage{(Inkscape) Color is used for the text in Inkscape, but the package 'color.sty' is not loaded}%
    \renewcommand\color[2][]{}%
  }%
  \providecommand\transparent[1]{%
    \errmessage{(Inkscape) Transparency is used (non-zero) for the text in Inkscape, but the package 'transparent.sty' is not loaded}%
    \renewcommand\transparent[1]{}%
  }%
  \providecommand\rotatebox[2]{#2}%
  \newcommand*\fsize{\dimexpr\f@size pt\relax}%
  \newcommand*\lineheight[1]{\fontsize{\fsize}{#1\fsize}\selectfont}%
  \ifx\svgwidth\undefined%
    \setlength{\unitlength}{463.13176332bp}%
    \ifx\svgscale\undefined%
      \relax%
    \else%
      \setlength{\unitlength}{\unitlength * \real{\svgscale}}%
    \fi%
  \else%
    \setlength{\unitlength}{\svgwidth}%
  \fi%
  \global\let\svgwidth\undefined%
  \global\let\svgscale\undefined%
  \makeatother%
  \begin{picture}(1,0.8173757)%
    \lineheight{1}%
    \setlength\tabcolsep{0pt}%
    \put(0,0){\includegraphics[width=\unitlength,page=1]{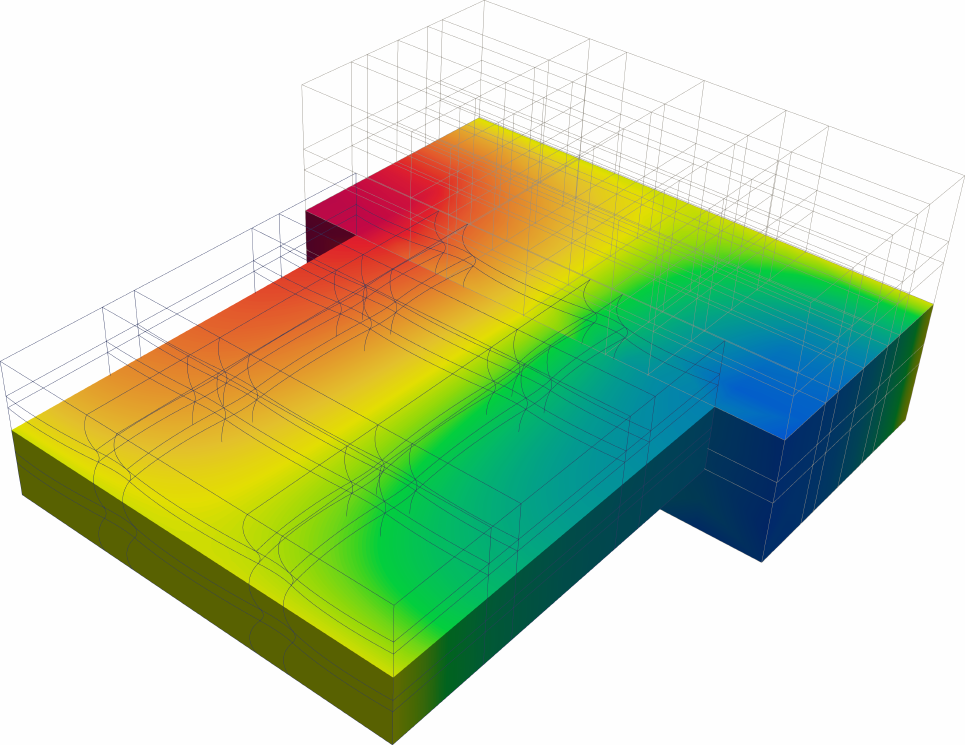}}%
    \put(0.46,-0.12){\color[rgb]{0,0,0}\makebox(0,0)[lt]{\smash{\begin{tabular}[t]{l}(c)\end{tabular}}}}%
  \end{picture}%
\endgroup%

  		}
  		\vspace{0.2cm}
  		\caption{Example described in Section~\ref{sub:poisson}. In~(a), (b), and~(c), the obtained solution of the Poisson's problem
  						 in~\eqref{eq:poisson_problem} for different refinements of the mesh of $\Omega^*_1$ are shown.}
  		\label{fig:contact_warp}
		\end{figure*}
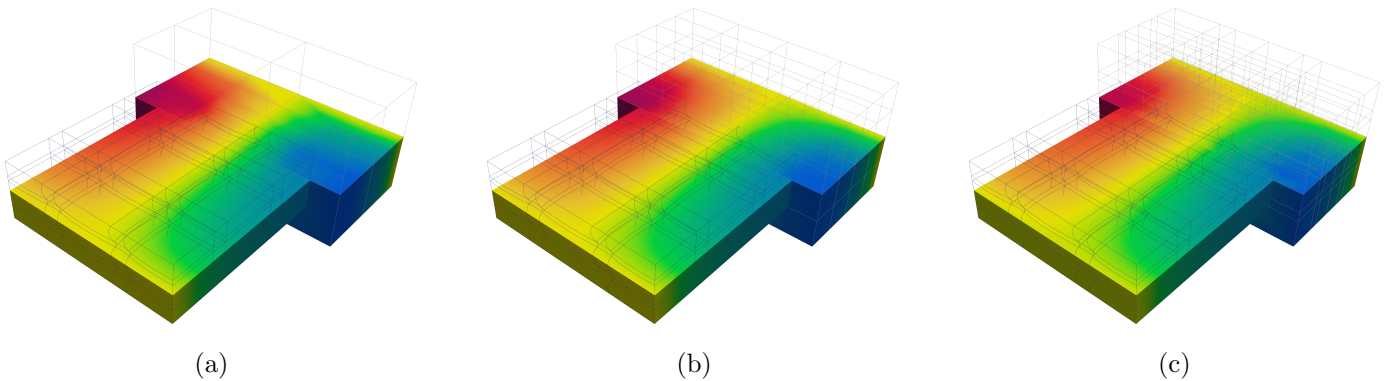

		To this end, we use our region extraction algorithm in order to identify the regions visible in Figure~\ref{fig:example-contact}~(b).
		In order to visualize the solution of the Poisson's problem, we associate to $\Omega_2$ the sinusoidal source
		\[
			f(x, y, z) = {\Big(\frac{\pi}{10}\Big)}^6 \sin\Big(\frac{\pi}{10} x\Big)
		\]
		and we enforce equality over the interface $\Gamma_{1,2}$ with the Lagrange multipliers method. Figure~\ref{fig:contact_warp} shows
		$u_1$ and $u_2$ according to three different refinements of the mesh of $\Omega_1$.

% section numerical_results (end)
\section{Conclusion}\label{sec:conclusion}
	In this work we have presented a novel region extraction algorithm for curvilinear drawings that allows to easily identify the regions 
	bounded by a set of planar curves. The algorithm has shown to be a powerful tool to precisely compute integrals over the interfaces of
	solids of piecewise polynomials defined in different meshes. In the literature, this kind of integrals are computed by simply
	approximating the regions in which the splines are represented by polynomials. These approximations polluted the quality of the
	integration and, in all our numerical tests, resulted in a lower precision than the one shown in
	Section~\ref{sub:computation_of_integrals_over_planar_domains}. The precise computation of integrals of this type come in handy in many
	practical operations.

	Three of such operations have been discussed in Section~\ref{sec:applications}. In the context of weak continuity enforcement, different 
	results can be achieved using different projectors for the Lee-Lyche-M{\o}rken quasi-interpolant or a different quasi-interpolant 
	altogether. Different choices could be, for example, the use of the standard $\mathcal{L}^2$ projector or yet the use of the
	projector proposed in~\cite{Zou:2018:IBD}. Both of these choices would take advantage of our algorithm for a precise computation of the 
	integrals. Other choices of the local knot interval $K_\ell$ different than the one taken in Section~\ref{sub:weak_continuity} can also
	influence the behavior of the result. Different approaches can be taken according to the needs of the specific application (locality,
	reproduction of the method etc.) and it is therefore our idea that the best results can be achieved only using application tailored
	methods.

	Among the applications discussed, the enforcement of the weak continuity is probably the more interesting from a geometric modeling point 
	of view. In a recent work~\cite{Masalha:2021:HPT}, Masalha and colleagues proposed a procedure for creating heterogeneous parametric 
	trivariate fillets. Despite the interesting geometric algorithms proposed, the obtained fillets are not connected to the original objects 
	and therefore new fillets need to be created every time a deformation applies to the input objects. The weak continuity constraints 
	proposed in this work can overcome this issue, creating a unique multipatch geometry for the entire volumetric model.

	Among the possible applications of Algorithm~\ref{alg:region-extraction} that have not been presented in this work, there are several 
	ones that are more tailored to the context of IGA, such as mortar methods and contact problems. These are promising research areas that
	we are currently investigating and will be the arguments of forthcoming works.

\section{Acknowledgment}
	P. Antolin and A. Buffa are partially supported by the ERC AdG project CHANGE n. 694515. P. Antolin, A. Buffa and E. Cirillo are
	partially supported by the SNSF through the project “Design-through-Analysis (of PDEs): the litmus test” n. 40B2-0 187094 (BRIDGE
	Discovery 2019). These supports are gratefully acknowledged.

%%%%%%%%%%%%%%%%%%%%%%%%%%%%%%%%%%%%%%%%%%%%%%%%%%%%%%%%%%%%%%%%%%%%%%%%%%%%%
\bibliographystyle{model1-num-names}
\bibliography{jabbrv,references}

\end{document}